\documentclass[11pt, final]{article}
\usepackage{amsmath,amssymb,amsthm}
\usepackage{bm}
\usepackage{graphicx}
\usepackage{fullpage}
\usepackage{color}
\usepackage{amssymb}
\usepackage{amscd}
\usepackage[mathscr]{eucal}
\usepackage{enumitem}
\usepackage{amsfonts}
\usepackage{float}
\usepackage{mathtools}
\usepackage{accents}
\usepackage{comment}

\theoremstyle{plain}
\newtheorem{Theorem}{Theorem}[section]
\newtheorem{Lemma}[Theorem]{Lemma}
\newtheorem{Proposition}[Theorem]{Proposition}

\theoremstyle{definition}
\newtheorem{Definition}[Theorem]{Definition}

\newcommand{\hl}[1]{#1}
\newcommand\highlighting[1]{#1}

\newcommand{\revv}[1]{}

\newcommand{\revision}[1]{%
  \begingroup
  #1%
  \endgroup
}

\newcommand{\bH}{{\bf H}}

\newcommand\calJ{\mathcal{J}}
\newcommand\calZ{\mathcal{Z}}
\newcommand\calP{\mathcal{P}}
\newcommand\bR{\mathbb{R}}

\newcommand\bN{\mathbb{N}}
\newcommand{\supp}{{\rm supp}}
\newcommand{\interior}[1]{\mathring{#1}}
\renewcommand{\interior}[1]{\accentset{\circ}{#1}}
\newcommand{\apT}{{\tilde{\tau}}}

\newcommand{\vep}{\varepsilon}

\newcommand{\mbold}[1]{\mathbf{#1}}

\newcommand{\abs}[1]{\lvert #1 \rvert}

\newcommand{\vnorm}[1]{\left\lVert #1 \right\rVert}

\newcommand{\floor}[1]{\lfloor #1 \rfloor}

\newcommand{\ceil}[1]{\lceil #1 \rceil}

\newcommand{\vset}[1]{\left\{ #1 \right\}}

\newcommand{\wh}[1]{\widehat{#1}}
\newcommand{\wt}[1]{\widetilde{#1}}
\newcommand{\mbb}[1]{\mathbb{#1}}

\newcommand{\mscr}[1]{\mathscr{#1}}
\newcommand{\mcal}[1]{\mathcal{#1}}

\newcommand{\eps}{\varepsilon}

\newcommand{\Ster}{S_{\infty}}
\newcommand{\hter}{h_{\infty}}

\title{Large-Deviation Analysis for Stochastic Models of Bacterial Evolution}

\author{
Robert Azencott\thanks{Department of Mathematics, University of Houston, Houston, TX, razencot@cougarnet.uh.edu}, \hspace*{2mm}
Brett Geiger\thanks{Department of Mathematical Sciences, High Point University, NC},
\hspace*{2mm}
Ilya Timofeyev\thanks{Department of Mathematics, University of Houston, Houston, TX, itimofey@cougarnet.uh.edu}, \hspace*{2mm}
}

\begin{document}

\maketitle

\begin{abstract}
Radical shifts in the genetic composition of large cell populations are \emph{\highlighting{rare events}} with quite low probabilities that direct numerical simulations generally fail to evaluate accurately. 
In this paper, we develop a theoretical large-deviation framework for a class of Markov chains modeling the genetic evolution of bacteria such as \emph{E. coli} \revision{in ``locked-box'' laboratory experiments.} In particular, we develop the cost function for discrete-time Markov chains that describe the daily evolution of histograms of bacterial populations.  \revision{We obtain explicit formulas for the cost function for interior histograms.} We also develop explicit formulas that can be used to numerically quantify the most likely evolutionary trajectories connecting an initial histogram and the target histogram.
\end{abstract}


\medskip


\section{Introduction}\label{mainobj}

We focus this theoretical study on stochastic models for the genetic evolution of bacterial populations. 
Stochastic modeling of evolutionary dynamics has been an active area of research for several decades (e.g.,~\cite{durrett2008,darwin2011,meleard2015,meleard2016,komarova2005}), and we utilize tools from stochastic analysis to study rare events that correspond to the emergence of non-dominant 
genotypes in long-term \revision{Lenski} evolutionary experiments \revision{of \emph{E. coli} (e.g.,~\cite{cooper1,heger,len94a,Fox2015,Barrick2013,Barrick2010,Cooper2001,Deatherage2017,Woods2011,Levy2015,Plank1979}).}
In particular, 
we recast a large class of such models as discrete-time Markov chains in the large-dimensional space of population histograms for bacteria with multiple genotypes.
These models, sometimes called ``locked-box'' models,  have been developed to describe, in the context of laboratory experiments, the daily evolution 
of finite-size bacterial populations~\cite{azencoop,rice,heger}. On day $n$, the bacterial population is characterized by the histogram $H_n$ of all genotype frequencies. Due to daily selection, in typical laboratory contexts the size $N$ of the population remains roughly constant.  Bacterial genotypes are characterized by their 
fitness, or equivalently by their growth factors, which remain fixed throughout the evolution of the bacterial colony. Random mutations are assumed to be roughly Poissonian and to occur at a fixed very small rate $m$. For a mutating cell of genotype $i$  the probability $q_{i,j}$ of generating a cell of genotype $j$ is assumed to be constant in time. Recall that many estimation techniques based on intensive simulations and experimental data have been implemented and tested (see~\cite{azencoop,heger,Gordo2013,Illinworth2012,Simon2013,Ross1996}) to evaluate the mutation rate $m$ and the growth factors (or selective advantages) of genotypes of interest.

\revision{The mathematical model considered here is similar
to the stochastic evolutionary model studied by R. Azencott et al. in~\cite{azencoop} to analyze the long-term laboratory experiments of T. Cooper on genetic evolution of \emph{Escherichia coli} bacteria. 
Multiple other publications  (e.g.,~\cite{cooper1,heger,len94a,Fox2015,Barrick2013,Barrick2010,Cooper2001,Deatherage2017,Woods2011,Levy2015,Plank1979}) have focused on the same theme, following  the well-known “locked-box”  scheme of  the Lenski experiments on \emph{E. coli}. 
In these experiments, on day $n$, the current cell population has size $N$ and grows freely until nutrient exhaustion. One then extracts (by dilution, for instance) a random sample of approximately $N$ cells, which constitutes the next-day population. 
``Locked-box'' models of genetic evolution for such experiments implement a succession of  ``daily'' cycles comprising three steps: growth phase, mutations, and random selection of a sub-sample of fixed size $N$.
In~\cite{azencoop}, the main challenge was to develop and test new accurate estimators of the main model parameters, namely,  the  mutation rate and selective advantage of beneficial mutations in \emph{E. coli}. To this end, genetically homogeneous subpopulations were tagged with biological markers. After parameter estimation,  the goodness of fit between the locked-box models and the experimental data was validated via intensive simulations.

Motivated by this  practical validation of  locked-box models, we undertook in  the present paper  the mathematical task of developing implementable  algorithms to study rare genetic events of high interest in Lenski-type experiments,   with an emphasis on the rare fixation of  genotypes that do not have maximum fitness. 
Recall that fixation of a specific genotype means that its frequency quickly becomes unusually large, or even overwhelming, in the cell population. 
Fixation of a genotype with moderate fitness has a very small probability, which vanishes at an exponential rate for a large cell population size $N$; hence, we naturally focus our present study on developing  a rigorous mathematical large deviations theory for the analysis
of rare genetic events in locked-box models.

Our emphasis is on developing numerically implementable algorithms to actually compute exponential rates of decay for probabilities of rare genetic events. Our first key theoretical tool is to determine the large-deviation functionals controlling rare process trajectories in the path space of all sequences of population histograms. We actually obtain fairly explicit theoretical formulas for these large-deviation functionals, enabling the numerical analysis of fixation events as well as the numerical computation of the most likely evolutionary histogram trajectory linking an initial population histogram $H$ to a terminal histogram $G$. The present paper focuses on rigorously proving the theoretical large-deviation results and the associated explicit formulas needed for the numerical computation of the most likely evolutionary paths leading to genotype fixation in locked-box bacterial experiments.

The numerical applicability of rare-event analysis has not often been exploited in concrete models of cell populations. Previous large-deviation results for stochastic population evolution have involved theoretical asymptotic studies, such as~\cite{champ01,champ02,champ6,champ8}. These papers have studied very general asexual population evolution in the space of phenotypic trait vectors. In such models, Darwinian evolution of asexual populations is driven by birth and death rates, which are themselves dependent on phenotypic traits. These traits are approximately transmitted to offspring, with rare but important variations due to gene mutations. Competition for limited resources forces a permanent or roughly periodic selection. Mutations are assumed to be rare enough such that, most of the time, only one currently dominant trait vector can coexist with the trait vectors of emerging mutants.

Here, we consider evolutionary contexts quite different from~\cite{champ01,champ02,champ6,champ8} since we focus on  locked-box models that combine two random steps: mutations and daily dilution. Thus, observed rare events can be triggered by unexpected combinations of these two steps, and it is important to understand the joint impact of these steps on experiments in bacterial evolutionary dynamics 
\cite{Barrick2013,Desai2007,heger,Levy2015,Peng2017}.
Our large-deviation study focuses on  specific  discrete Markov chains which live in the Euclidean  space of cell-population histograms and  describe the genetic evolution of bacterial populations in locked-box experiments. Our goal here is to rigorously develop  formulas to quantify rare genetic events such as  fixation of genes with moderate fitness, with the aim of implementing them numerically (as in our companion paper~\cite{azencott2023rare}).

We develop an explicit rate function defined for all paths connecting an initial population histogram $G$ to a target histogram $H$ in a fixed number of time steps. Minimizing the rate function over the set $S(G,H)$ of all histogram paths linking $G$ to $H$ then yields formulas to quantify the log-probability of rare events, for instance, when the target histogram $H$ describes the fixation of a non-dominant genotype. We also develop an explicit second-order reverse recurrence equation satisfied by optimal paths minimizing the rate function over $S(G,H)$. These recurrence equations essentially solve the optimization problem of finding the most likely path maximizing the probability of fixation for any specific genotype.

Extensive numerical results based on our theoretical developments are presented in our companion paper~\cite{azencott2023rare}, where we discuss the numerical implementation of our algorithms for computing the most likely evolutionary trajectory linking two population histograms $H$ and $G$, and explore how to handle numerical complexity when the number of genotypes increases. In particular, Ref.~\cite{azencott2023rare} extensively analyzes models with $g = 3, 4, 5$ genotypes, as well as discussing how to optimize the search for the optimal penultimate point in the backward search algorithm presented in Section~\ref{s:numtraj}. The complexity of the backward search algorithm arises from the necessity to search for the optimal penultimate point in the path of histograms. In~\cite{azencott2023rare}, we develop, implement, and test accelerated algorithms for this numerical search. 
In the locked-box evolutionary models parametrized by the experimental results derived in~\cite{azencoop} from the T. Cooper experiments,  we also implemented  and tested (see~\cite{su_thesis,azencott2023rare})  a gradient descent algorithm for computing the most likely evolutionary path linking any two given population histograms.
The Ph.D. thesis~\cite{su_thesis} contains a detailed discussion about the model and numerical techniques,
including a numerical validation of the rate functional derived in this paper using the importance sampling algorithms.
The numerical techniques implemented in~\cite{su_thesis,azencott2023rare} allowed us to compute the most likely paths for numbers of genotypes  $g = 7, 8, 9, 10$  using a multi-core implementation. 
This demonstrates that our approach can successfully handle relatively high-dimensional problems.
In addition, we believe that, with a high-performance computing implementation, our optimization algorithm could be extended to handle models with 20--30 genotypes.
}

\revision{The locked-box models considered here assume a constant population size $N$ at the beginning of each day, Poisson mutations with constant rates, fixed growth factors, and no interaction with the environment. These assumptions are common in studies of locked-box models motivated by the Lenski experiments. Mathematically, some parameter variations can be accounted for by introducing additional genotypes. However, as discussed below, the practical computational limit of our current implementation is approximately 10--20~genotypes. Nevertheless, models with 10 genotypes can already provide valuable insight into the mechanisms of genetic evolution in locked-box experiments. We also believe that extending the model to variable mutation rates is possible, provided that the expected number of mutations does not vary significantly from one day to the next.
Extending the present framework to models that describe interactions with the environment represents one of the most challenging directions for future research. In this context, in many epidemic models
(see, e.g.,~\cite{assaf2017wkb,gomez2023markovian,papageorgiou2025enhanced}), the environment can directly affect the transition structure. This makes the transition mechanism state-dependent, time-dependent, or constrained, and the analysis of the one-step cost function becomes more difficult. In particular, it is unlikely that the one-step cost function can be obtained in an explicit form. One possible extension is a deterministic time-dependent environment where the mutation matrix and growth factors may depend on time. In that case, it should be possible to define a one-step cost function and extend the analysis to this case. Random environments are much more challenging and extending LDT analysis to such cases would probably require a substantially modified framework.}

\subsection{\highlighting{Novelty}
 and Scientific Contribution}
\revision{This  paper addresses  several mathematical challenges in order to rigorously apply large deviations theory (LDT) in a Markov chain framework to characterize rare events for the stochastic locked-box models of bacterial genetic evolution. 
The mathematical developments required here go beyond classical LDT results,
which certainly do not directly provide the estimates and explicit formulas required for this setting. 
For sets $A$ of evolution trajectories, typical LDT results are of the form 
$N^{-1} \log P(A) = -\Lambda(A) + o(1)$ (with $o(1)\to 0$ as $N \to \infty$), 
where $N$ is the population size and $\Lambda(A) \ge 0$ is derived from an LDT rate function $\lambda(\text{path})$  defined on all process paths in the space of population histograms.}

\revision{However, we need to obtain \emph{explicit uniform bounds} 
 on the speed at which $o(1)$ tends to zero for large $N$.
Indeed,  the rate function $\lambda(\text{path})$ is computed by minimizing the sum of two interacting rate functions, $\lambda_{mut}$ and $\lambda_{sel}$, corresponding to mutation and selection, respectively. Uniform speeds of convergence for $o(1)$ are essential for proving that this minimization yields the correct path-space rate function $\lambda$  (see Theorem~\ref{thtrajLD}).
Another delicate technical point is the derivation of local continuity estimates for the rate function $\lambda$. These estimates require, in particular, explicit H\"{o}lder continuity constants for the mutation rate function $\lambda_{mut}$ (see Proposition~\ref{propholdermut}). 
A further difficulty is that genotype frequencies are rational numbers. Thus, for fixed large $N$, the population histograms evolve within a very large finite discretization  of a convex Euclidean set. 
This discrete structure must be handled carefully, since the most likely histogram path realizing a rare event is computed in continuous space.
Finally, we obtain a nontrivial but fully explicit form of the rate function (see \eqref{def:lambda} and  \eqref{formulaCHG}) for generic paths in the space of histograms. This is an important practical aspect of this work, allowing us to develop an implementable computational strategy to search for the most likely path linking any prescribed pair of initial and terminal histograms. This goal is discussed briefly in Section~\ref{s:numtraj}, and its successful numerical implementation is presented in another publication~\cite{azencott2023rare}.
Overall, our results  provide a  significant analytical and practical advance in the concrete modeling and analysis of rare genetic events in locked-box experiments on large bacterial populations.}

\revision{Several major results in this paper are restricted to interior-population histograms, where every genotype has strictly positive frequency. This assumption is essential for deriving the explicit large-deviation rate function and the optimization algorithm using the reverse recurrence relation. A rigorous extension of the LDT analysis to boundary histograms, where one or more genotype frequencies are exactly zero, would require substantial theoretical advancements. 
In Section~\ref{sec:boundary}, we present one particular boundary case in which we treat the emergence of a single genotype and obtain an explicit form of the one-step cost function. This suggests that extensions to more general boundary cases may be possible. However, the reverse optimization algorithm in Section~\ref{s:numtraj} would need to be modified accordingly, since most likely paths may involve several intermediate boundary histograms.
We also would also like to point out that the lower bounds imposed on genotype frequencies may be chosen to be sufficiently small in many situations, so that the explicit formulas developed here provide intuition and insights into evolutionary scenarios involving boundary histograms.}

\subsection{\highlighting{Manuscript Organization}}
In Section~\ref{stochmodel}, we describe the Markov chain dynamics that serves as the stochastic model of bacterial evolution.  In Section~\ref{LDcycles}, we establish a large-deviation framework for these Markovian successions of daily cycles, in order to derive an explicit formula for the one-step cost function given by \eqref{formulaCHG} in Theorem~\ref{thCHG}.  These results are applied in Section~\ref{LDPathTheory} to formulate a large deviations theory for evolutionary trajectories taking values in the space of population histograms.
The main application of our large-deviation framework is considered in Section~\ref{mostlikely}, where we develop theoretical results for computing the most likely path connecting an initial population state  $G$ 
 to a fixed given terminal state. We conclude our main results in Section~\ref{s:numtraj} by deriving an explicit reverse recurrence relation that can be used to numerically compute the most likely path connecting the initial $G$ and the final $H$ in the space of histograms.

%%%%%%%%%%%%%%%%%%%%%%%%%%%%%%%%%%%%%%%%%%%%%
\section{Stochastic Model for Bacterial Evolution Experiments}\label{stochmodel}
To model the main features of random bacterial evolution, we focus on a class of Markov chains often used in this context~\cite{azencoop,rice,heger}. We assume that all possible cell genotypes belong to a finite \hl{set} 
$\{1, 2, \cdots, g \}$. 
Cells of genotype $j$ are called \emph{j-cells} here and have a fixed growth factor $F_j > 0$. We always order genotypes by increasing fitness so that $F_1 < F_2<\cdots < F_g$. The genotype $g$ with the highest fitness is called \hl{dominant}. 
We also denote the vector of growth factors as $F \coloneqq \left[ F_1,\cdots, F_g \right]$.

When a $j$-cell divides, its genotype is typically inherited by the two daughter cells, unless a rare random genotype mutation from $j$ to $k \neq j$ occurs. Mutation occurrences have an approximately Poisson distribution with a very small mutation rate, $m$, with typical values in the range $10^{-9} \le m \le 10^{-6}$. 

%%%%%%%%%%%%%%%%
\begin{Definition}\label{Qjk}
\hl{When a} 
mutation occurs during the division of a given  $j$-cell, the conditional probability that a mutant daughter cell will be a $k$-cell is denoted $q_{j,k} \ge 0$. The transition mutation matrix $Q \in \bR^g \times \bR^g$, with entries $q_{j,k} \ge 0$, naturally verifies  $\sum_k q_{j,k} =1$ and $q_{j,j} =0$. 
\end{Definition}

\begin{Definition}\label{P}
The population size $N$ and the process parameters $\calP=\{m, g, F_1,\cdots, F_g, Q\}$
define our stochastic evolutionary model for finite-size bacterial populations.
\end{Definition}
In all proofs throughout the manuscript, we assume that these
parameters are fixed. We will also assume a large population size $N \gg 1$ and a small mutation rate $m \ll 1$ and utilize expansions to obtain the leading-order terms in $N$ and $m$ in some of the proofs. For typical bacterial population experiments, one has $m<10^{-6}$ and $N> 10^5$.

We detail below the three successive phases (growth, mutation, random selection) implemented in 
each daily evolutionary cycle and we define the space of population histograms which quantify the concentrations of $j$-cells for each $j \in \vset{1,\dots,g}$.  Then, we describe the Markov transition kernel associated with each daily cycle.

%%%%%%%%%%%%%%%%%%%%%%%%%%%%%%%%%%%%%%%%%%%%%
\subsection{The Space of Genetic Histograms} \label{histograms}
Each bacterial population is described by a histogram of genotype frequencies $H \coloneqq \left[ H(1),\cdots, H(g) \right]$. We often denote by $H_n \coloneqq \left[ H_n(1),\cdots, H_n(g) \right]$ the histogram of bacterial frequencies on the $n$-th day.
Next, we introduce some basic definitions.
%%%%%%%%%%%%%%%%%
\begin{Definition} \label{Nrational}
A matrix $A$ will be called \emph{$N$-rational} if for a positive integer $N$ all the coefficients of $N A$ are non-negative integers.  
\end{Definition}
%%%%%%%%%%%%%%%%%
\begin{Definition} \label{calH}
Denote by  $\mathcal{H} =\{H \in \bR^g\}$ the set of all possible population histograms. Each histogram $H$ is a vector of length $g$ such that $0 \le H(j) \le 1$ and $\sum_j H(j) = 1$. Note that $\mathcal{H} \subset \bR^g$ is compact and convex.
\end{Definition}
%%%%%%%%%%%%%%%%
\begin{Definition} \label{calHN}
The subset $\mathcal{H}_N$ of $\mathcal{H}$ is the set of $N$-rational 
vectors $H \in \mathcal{H}$. Note that $\mathcal{H}_N$ is finite, with $card(\mathcal{H}_N) \leq (N+1)^g$.

In a cell population of size $N$, denote as $N_j$  the number of  $j$-cells, and  $H(j) = N_j/N$ their concentration. This  ``genetic histogram'' $H = [ H(1),\dots,H(g) ] $  then belongs to  $\mathcal{H}_N$.
\end{Definition}
 
The sets $\mathcal{H}$ and $\mathcal{H}_N$
are endowed with the \emph{$L_{\infty}$-distance}
$$
\vnorm{H - G} = \max_j\,{ \abs{H(j) - G(j)}} \quad \text{for all } \; H,G\in\mathcal{H}.
$$ 

%%%%%%%%%%%%%%%%%
\begin{Definition}
\label{def:suppmin}
The ``boundary'' of $\mathcal{H}$ is the set of histograms $H$ for which at least one $H(j)=0$.  To quantify  closeness to the boundary for any $H \in \mathcal{H}$ we define the support $\supp(H)$ and the essential minimum $b(H) > 0$  by 
\begin{equation}\label{supp(H)b(H)}
\supp (H) = \vset{j \, | \,H(j)>0} \quad \text{and} \quad b(H) = \min_{j \in \supp (H)} \; H(j) .
\end{equation}
\end{Definition}
Note that $b(H) \geq 1/N$ for all $H \in \mathcal{H}_N$.  
\revision{The definition above distinguishes interior histograms with $b(H)>0$ from boundary histograms with $b(H)=0$. For interior histograms, it is assumed that all genotypes are present in the population, although some frequencies may be very small, e.g., $b(H)=N^{-1}$. 
For boundary histograms, one or more genotypes are absent from the population. Several of our key results are developed for interior histograms; however, since we can consider interior histograms with $b(H)=N^{-1}$, these results provide useful insight into evolutionary scenarios involving the emergence of new genotypes.}
Genetic evolution during $T$ days is then described by the histogram \hl{trajectory }
$\bH \coloneqq \{H_1, H_2,\cdots, H_T\}$, where $H_n \in \mcal{H}_N$ is the population histogram on day $n$ for $1 \leq n \leq T$; each histogream $H_n$ describes the distribution of $g$ genotypes on day $n$.
We now present the three phases of each ``daily'' cycle.

Our ``one-day'' time unit refers to the fixed duration of each single evolutionary cycle. In laboratory experiments, single-cycle durations are roughly constant but may be shorter than 24 h in real time. In the T. Cooper experiments, for instance (e.g.,~\cite{Cooper2001}), the growth phase typically takes 8--12 h, after which the nutrients are 
exhausted and cells become dormant.

%%%%%%%%%%%%%%%%%%%%%%%%%%%%%%%%%%%%%%%%%%%%%
\subsection{Path Space of Histogram Sequences}
Genetic evolution is modeled as a sequence of daily cycles indexed by day $n$. To simplify our evolution model, we assume that all daily random mutations occur nearly simultaneously, after the growth phase. One can develop analogous models where random mutations can occur at any time during the growth phase (see, for instance,~\cite{schoneman_thesis}); but these models involve stochastic differential equations driven by Poisson processes, so that rigorous large-deviation frameworks similar to the theory developed here require more complicated mathematical proofs. We consider that the simplified models studied here are sufficient for quantitative analysis of genetic evolutionary pathways in realistic experimental setups.

We now describe the stochastic dynamics of each daily cycle.
At the beginning of day $n$, the population $pop_n$ is always of size $N$ and is identified by its genetic histogram $H_n \in \mbb{R}^g$, recording the frequencies $H_n(j)$ of $j$-cells in population $pop_n$. 
The day $n$ cycle initiated by $pop_n$ generates $pop_{n+1}$ (and hence its genetic histogram $H_{n+1}$) in three successive phases. 

{\textbf{Phase 1:}
}  Purely deterministic growth, with growth factor $F_j$ for the $j$-cell colony.

{\textbf{Phase 2:}}  Random mutations occur simultaneously at the end of the growth phase.

{\textbf{Phase 3:}}  After growth and mutations, a random sub-sample of fixed size $N$ is extracted from the current population, and will constitute $pop_{n+1}$. 

\begin{figure}[H]
    \centerline{\includegraphics[width=0.8\textwidth]{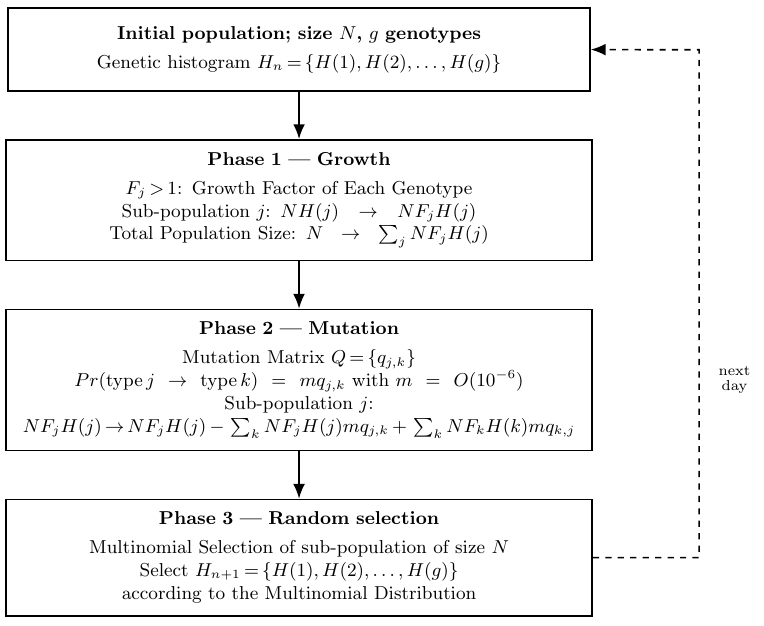} }
    \caption{\revision{Schematic representation of the Markov chain dynamics for histograms $H_n \to H_{n+1}$ during the three phases in the daily cycle. }}
    \label{fig1}
\end{figure}

The three phases are schematically illustrated in Figure \ref{fig1}.

%%%%%%%%%%%%%%%%%%%%%%%%%%%%%%%%%%%%%%%%%%%%%%
\subsection{Phase 1---Deterministic Growth} 

In actual experiments on bacterial evolution (see, e.g.,~\cite{cooper1,heger}), the daily multiplicative growth factors $F_j$ for observable genotypes are typically in the range $[20,300]$. These growth factors can be computed as $F_j = \exp(\Delta  \times a_j)$, where $\Delta $ is the duration of the growth phase and $a_j > 0$ is the ``selective advantage''  of $j$-cells.
Detectable selective advantages over the ancestor genotype are typically larger than 0.01 \revision{(e.g.,~\cite{azencoop})}.

On day $n$, the initial population $pop_n$ has genetic histogram $H_n$, and the $j$-cell colony has initial size $N H_n(j)$. The $j$-cell colony then grows according to its growth factor, $F_j$, and reaches the size  $N F_j H_n(j)$. At the end of the growth phase, the population reaches a much larger  size, 
$
\ceil{N \langle F, H \rangle} = \ceil{ \sum_j N F_j H(j) }
$, 
where
$\langle F, H \rangle$ is the inner product of  vectors $F,H$ and $ \ceil{u}$ is the smallest integer $\geq u$.
The population genetic histogram after growth, $\varphi$, 
is then given by 
$
\varphi_j(H) = \ceil{ N F_j H(j) }  /  \ceil{ N \langle F, H \rangle}$.  
We naturally approximate $\varphi(H)$ by $\Phi(H) \in \mathcal{H} $, which we define as
\begin{equation} \label{phi}
\Phi_j(H) = F_j H(j) / \langle F, H \rangle.
\end{equation} 
\hl{For }
$N>20$ and all $H \in \mathcal{H}$, the approximation error is
\begin{equation} \label{intervalphi}
\vnorm{\Phi(H) - \varphi(H)} \leq  \frac{4 F_g }{N F_1}.
\end{equation}

%%%%%%%%%%%%%%%%%%%%%%%%%%%%%%%%%%%%%%%%%%%%%%
\subsection{Phase 2---Random Mutations} 

As outlined earlier, in our simplified model, we consider that on day $n$ all random mutations occur {\hl{simultaneously,}} right after the growth phase. 

\subsubsection{Random-Mutation Matrices and Mutation Rates} \label{lincons} 
\vspace{6pt}
\begin{Definition} \label{Rn}
During the mutation phase of day $n$, a random number $R_n(j,k)$ of $j$-cells mutate into $k$-cells. All the $R_n(j,k)$ are non-positive integers and $R_n(j,j)=0$ for all $j$. Denote by $R_n$ the
corresponding $g \times g$ random matrix.
\end{Definition}

On day $n$, at the end of the growth phase, each  $j$-cell colony has just reached the finite size $F_j N H_n(j)$. This forces the random numbers $R_n(j,k)$ of $j$-cells mutating into another genotype $k\neq j$ to verify a set of natural linear inequalities, which we explicitly develop below, before describing technically the joint probability distribution of the $R_n(j,k)$.

%%%%%%%%%%%%%%%%
\begin{Definition} \label{mutM}
The mean emergence rate $m$ of mutants is assumed to have the same very small value for all $j$-cell colonies. For bacterial populations, one typically has $10^{-9} \leq m \leq 10^{-6}$. Recall that $Q =\{q_{j,k}\}$ is the fixed conditional mutation transition matrix defined in \eqref{Qjk}. So for a large $N$, we want the conditional expectation $E(R_n(j,k) | H_n )$ to be very close to $m F_j N H_n(j) q(j,k).$
\end{Definition}

%%%%%%%%%%%%%%%%%%%%%%%%%%%%%%%%%%%%%%%%%%%%%%
\subsubsection{Mutation Constraint Set $K(H)$} 
Given $H_n$, the total number of mutants $\sum_k R_n(j,k)$ emerging from the $j$-cell colony must be smaller than $N F_j H_n(j)$, which is the number of $j$-cells after growth. This imposes $g$ linear constraints on $R_n$, namely, 
$$
\sum_k R_n(j,k) < N F_j H_n(j) \;\; \text{whenever}\;\; H_n(j) >0.
$$
\hl{Hence,} the random matrix  $r_n= R_n/N$ must belong to a convex set of matrices $K(H_n)$, which we now describe.

%%%%%%%%%%%%%%%%
\begin{Definition}[Constraint Sets]
\label{KNconstr}
Let $\calZ$ be the set of all $N$-rational $g \times g$ matrices, and note that $R_n / N \in \calZ$. For each histogram $H$ and each $j$, define $K(j, H)$ as the set of all $g \times g$ matrices $r$ with non-negative coefficients, verifying for all $j,k$
\begin{equation} \label{Kj}
\begin{cases}
\sum\limits_{k=1}^g  r_{j,k} <  F_j H(j) & \text{if } \quad H(j) > 0;\\
r_{j,k}  =  0, & \text{if }\quad H(j) = 0 \text{ or } q_{j,k} = 0.
\end{cases}
\end{equation}
\hl{Also, define} 
 the convex set of matrices $K(H)  = \bigcap_{j=1}^g \; K(j, H)$ and let $K_N(H) = \calZ \cap K(H)$ be the subset of all N-rational matrices in $K(H)$.
\end{Definition}
The matrices $r \in K(H)$ are potential values for all the random matrices $r_n =R/N$, which are $N$-rational and must, therefore, belong to $K(H_n)$. 
\revision{The constraints in \eqref{Kj} reflect the requirement that the number of mutants emerging through the mutation $j \to k$ is smaller than the total number of cells of genotype $j$.}

Note that $K_N(H) \subset K(H)$ by definition.  Due to \eqref{Kj}, we have $ r_{j,k} \leq F_g$ and $card( K_N(H) )\leq [ (1+F_g)(N+1) ]^{g^2}$ for all $H\in \mathcal{H}$, $r \in K(H)$, $(j,k)$, and $N$. The next technical lemma quantifies how accurately $K_N(H)$ approximates  $K(H)$ for a large $N$.

%%%%%%%%%%%%%%%%
\begin{Lemma} \label{density}
Fix $a>0$. For $N > \frac{g^2}{a F_1}$, any $H \in \mathcal{H}$ with $b(H) \geq a$, and any $r \in K(H)$, there is an $N$-rational matrix $s = \{ s_{j,k} \}$ such that 
\begin{equation}\label{ss}
s \in K_N(H) \subset K(H), \qquad 
\supp(s) = \supp(r), \qquad|| s - r || \leq g/N. 
\end{equation}
\end{Lemma}
%%%%%%%%%%%%%%%%%
\begin{proof} See  Appendix \ref{sec:ap2}.
\end{proof}

\revision{Since the population at the beginning of each day is assumed to have size $N$, we need to work with $N$-rational histograms. However, for analytical purposes it is sometimes necessary to consider histograms with arbitrary real-valued frequencies. Thus, we need to estimate how close an $N$-rational histogram can be to a given histogram with real-valued frequencies.}
We build our stochastic mutation model to enforce 
$
P( R_n/N \in K(H) \; | \; H_n = H ) = P( R_n/N \in K_N(H) \; | \; H_n = H ) = 1
$
We also want to ensure that for large $N$, the conditional joint distribution of the $R_n(j,k)$,  given $H_n$, will become extremely close to a product of  independent Poisson distributions with means $m\, q_{j,k} N F_j H_n(j)$. To this end, we first introduce the companion matrices $Z_n$.
%

%%%%%%%%%%%%%%%%%%%%%%%%%%%%%%%%%%%%%%%%%%%%%%
\subsubsection{The Poisson Companion Matrices $Z$} 
\label{compMat}
Given $H_n = H$, let $Z$ be a matrix of independent random variables $Z_{j,k}$ having Poisson distributions with respective means 
$
\mbb{E}[ Z_{j,k} \; | \; H_n= H ] = m q_{j,k} N F_j H(j).
$ 
We show that the conditional probability $P( Z/N \in K(H) \; | \; H_n = H ) $ tends to 1 very fast for large $N$ (see Theorem~\ref{companionTH}). So, for all histograms $H$ and matrices $z$ with non-negative coefficients, we define the conditional distribution of $R_n$, given $H_n=H$, by 
$
P( R_n/N = z \; | \; H_n = H )  = P(Z/N = z \; | \; Z/N \in K(H)). 
$
Since $Z/N$ is $N$-rational, this forces $R_n/N$ to be $N$-rational, with
$$
P( R_n/N \in K_N(H) \; | \; H_n = H ) = P( R_n/N \in K(H) \; | \; H_n = H ) = 1.
$$

\hl{We now} complete the analysis of the mutation phase by computing the population histogram at the end of Phase 2.

%%%%%%%%%%%%%%%%%%%%%%%%%%%%%%%%%%%%%%%%%%%%%%
\subsubsection{Population Histogram After Mutations}
\label{JnMut}
\revision{Recall that at the end of the growth phase the population becomes $N F_j H_n(j)$, although some of these values may not be $N$-rational.}
Accounting for all the random mutations that occurred during Phase 2, each $j$-cell colony has lost $\sum_k R_n(j,k)$ outgoing mutants and gained $\sum_k R_n(k,j)$ incoming mutants. But the total population size has not changed since the end of Phase 1 and is given by $\ceil{N \langle F, H_n \rangle}$. Hence, the random number  of $j$-cells at the end of Phase 2 is given by
\begin{equation} \label{Wn}
 \ceil{N F_j H_n(j)} - \sum_k R_n(j,k) + \sum_k R_n(k,j).
\end{equation}
\hl{Hence}, given $H_n= H \in \mathcal{H}_N$ and $R_n/N = r \in K_N(H)$, the population histogram $J_n$ at the end of Phase 2 is a deterministic function $J_n= \calJ(H_n, R_n/N, N)$ of $H_n, R_n/N$ and $N$. The histogram-valued function $(H, r, N) \to \calJ(H,r,N)$ does not depend on $n$ and is defined for all $H \in \mathcal{H}$ and all $r \in K(H)$. We can rewrite \eqref{Wn} to get
\begin{equation} \label{WWn}
J_n(j) = \calJ_j (H,r,N) =\frac{ \ceil{ N F_j H(j)} }{ \ceil{ N \langle F, H \rangle }} + \; \frac{N}{\ceil{ N \langle F, H \rangle } } \left[- \sum_k  r_{j,k} + \sum_k  r_{k,j} \right].
\end{equation}

\hl{One has} $\left|\frac{\ceil{ N F_j H(j)}}{\ceil{ N \langle F, H \rangle }} - \frac{F_j H(j)}{\langle F, H \rangle} \right| \leq \frac{4 F_g}{N F_1}$ due to \eqref{intervalphi}. Therefore, given $H_n= H \in \mathcal{H}_N$ and $R_n/N = r \in K_N(H)$, the histogram $J_n(j) =\calJ_j(H,r,N) $ is well approximated by
\begin{equation} \label{psi}
\Psi_j(H, r) = \frac{1}{\langle F, H \rangle} \left(F_j H(j) - \sum_k  r_{j,k} + \sum_k  r_{k,j}\right) .
\end{equation}

\hl{The function} $\Psi_j(H, r)$ is well-defined for all $H \in \mathcal{H}$ and $r \in K(H)$. Moreover, one has $\Psi_j(H, r) \geq 0$ and $\sum_j  \Psi_j(H, r) =1$, so that $\Psi(H, r)$ is always a histogram belonging to $\mathcal{H}$.

One can directly verify  that when $N > 20,$
\begin{equation} \label{JJn}
\vnorm{\calJ(H, r, N) - \Psi(H, r)} \leq \frac{13 F_g}{N F_1}
\end{equation}
for all $H \in \mathcal{H}$ and  $r \in K(H).$ In particular, with probability 1 one has
\begin{equation} 
\label{Jn}
\vnorm{J_n - \Psi(H_n , R_n/N)} \leq \frac{13 F_g}{N F_1}.
\end{equation}
\hl{Since the histogram} $J_n$ is  $ \tilde{N}$-rational with 
$\tilde{N} = \ceil{N \langle F, H \rangle}$, its essential minimum $b(J_n)$ must verify 
\begin{equation} \label{bJn}
b(J_n) \geq \frac{1}{\tilde{N}} \geq \frac{1}{N \langle F, H_n \rangle} \geq \frac{1}{N F_g}. 
\end{equation}
\hl{For} $H \in \mathcal{H}$ and $r \in K(H)$, Equations~\eqref{Kj} and \eqref{psi} yield, for each $j$, 
\begin{equation} \label{zeroPsi}
\Psi_j(H, r) = 0 \quad \text{if and only if} \; \; H(j) = r_{j,k} = r_{k,j}= 0 \;\; \text{for all} \; k.
\end{equation}

\hl{The relations} \eqref{Wn} and $R_n/N \in K_N(H_n)$ imply $J_n(j) = 0$ if and only if $H_n(j) = R_n(k,j) = 0$ for all $k$. Similarly, for all $N$-rational $H$ and $r \in K(H)$, Equation~\eqref{WWn} shows that $\calJ_j(H, r, N ) = 0$ if and only if $H(j) = r_{k,j} = 0$ for all $k$. Hence, one has $\supp(J_n) = \supp( \Psi(H_n , R_n/N) )$ for all $H_n$ and $R_n/N \in K(H_n)$.  In addition, $\supp(J(H, r, N))  =  \supp( \Psi(H, r) )$ for all $N$-rational $H$ and $r \in K(H)$.
Equation~\eqref{zeroPsi} proves that for $H, H' \in \mathcal{H}$, $r \in K(H),$ and $r' \in K(H')$,

\begin{equation} \label{supppsiHh}
\begin{aligned}
& \supp(\Psi(H,r) ) = \supp(\Psi(H', r') \; \; \text{whenever} \;\; \\
& \supp(H) = \supp(H') \text{ and } \supp(r) = \supp(r').
\end{aligned}
\end{equation}

We conclude the study of Phase 2 with a couple of relations that will be useful later on.  For all $H$ and $r \in K(H)$, the partial derivatives of $\Psi(H,r)$ are given by 
\begin{equation} \label{difPsi}
\begin{cases}
\frac{\partial}{\partial H(i)} \Psi_j  =  \frac{- F_i F_j H(j)}{\langle F, H \rangle^2} + \frac{1_{\{ i=j \}} F_j}{\langle F, H\rangle}, \\
\frac{\partial}{\partial r_{j,k}} \Psi_j  =  \frac{- F_j H(j)}{\langle F, H \rangle}, \\
 \frac{\partial}{\partial r_{k,j}} \Psi_j  =  \frac{F_j H(j)}{\langle F, H \rangle},
\end{cases}
\end{equation}
where the indicator $1_{\{i=j\}}$ equals $1$ if $i=j$ and $0$ otherwise.  Since $F_j H(j) \leq  \langle F, H \rangle$ and $\langle F, H\rangle  \geq  F_1$, we have 
\begin{equation} \label{lipPsi}
\vnorm{\Psi(H',r') - \Psi(H,r)} \leq \frac{3 g F_g}{F_1} \left(\vnorm{ H' - H } + \vnorm{ r' - r }  \right)
\end{equation}
for all $H', H \in \mathcal{H}$, $r \in K(H)$, and $r' \in K(H')$.

%%%%%%%%%%%%%%%%%%%%%%%%%%%%%%%%%%%%%%%%%%%%%%
\subsection{Phase 3---Random Selection} 
At the end of Phase 2 on day $n$, the current population $POP_n$ has a large  size $\tilde{N} \simeq N \langle F,H \rangle$ and histogram $J_n \simeq \Psi(H_n , R_n/N)$, with $\Psi$ given by \eqref{psi}. During Phase 3, one extracts from $POP_n$ a random sample of fixed size $N$. This sample  becomes the new initial population $pop_{n+1}$ on day $n+1$ and has genetic histogram $H_{n+1}$. Phase 3 is thus a simplified emulation of natural selection. 
%%%%
%%%%
%%%%
%%%%

The multinomial distribution $\mu_{N,J}(V)$ parameterized by $N$ and histogram $J$ is defined for all vectors $V \in \mbb{N}^g$ having  integer coordinates $V(j) \geq 0$ such that  $\sum_j V(j) = N$ by 
\begin{equation} \label{multinomial}
\mu_{N, J}(V) = N! \; \prod_{j \in \supp(V)} \; \frac{J(j)^{V(j)}}{V(j)!}. 
\end{equation}

\hl{When} $J_n(j) =0$, no mutant of genotype $j$ is present before selection, so that $H_{n+1}(j) = 0$. Hence, $\supp(H_{n+1}) \subset \supp(J_n)$ with probability 1. For $G \in \mathcal{H}_N$, all coordinates of $V = N G$ are integers, and one has
\begin{equation}\label{multinomial2}
P( H_{n+1} = G \; | \; H_n ; R_n ) = P(H_{n+1} = G \; | \; J_n) = 
\begin{cases}
\mu_{N, J_n}(N G),&  \supp(G) \subset \supp(J_n), \\
 0, & \text{otherwise.}
\end{cases}
\end{equation}

\hl{The multinomi}al distribution \eqref{multinomial} has mean $N J$. Hence, Equation~\eqref{multinomial2} implies
\linebreak $
\mbb{E}[ H_{n+1} \, | \, H_n , R_n ]  =  \mbb{E}[ H_{n+1} \, | \, J_n] = J_n. 
$
%

%%%%%%%%%%%%%%%%%%%%%%%%%%%%%%%%%%%%%%%%%%%%%%
\subsection{Markov Chain Dynamics in the Space of Histograms} \label{TransProb}
At the completion of the three phases during a daily cycle, the population size always returns to the large but \emph{\hl{fixed size}  
 $N$}. The cycle on day $n$ thus induces a stochastic transition $H_n \to H_{n+1} $ in the space of genetic histograms $\mathcal{H}\subset \mbb{R}^g$. The succession of daily
 cycles just described above generates a \emph{\hl{time-homogeneous Markov chain}} 
$\{H_n\to H_{n+1}\}$ on the state space $\mathcal{H} $ of all histograms. However, for each fixed population size $N$, the actual state space of the Markov chain $\{H_n\to H_{n+1}\}$ is the finite set $\mathcal{H}_N$ of $N$-rational histograms, which has size $card(\mathcal{H}_N) \leq (N+1)^{g}$.
The Markov chain' {\hl{transition kernel}} 
$\mathcal{Q}_N (H, G) = P(H_{n+1} = G \; | \; H_n = H)$ for $H, G \in \mathcal{H}_N$ is given by the finite sum 
\begin{equation}\label{kernel}
\mathcal{Q}_N (H, G) = \sum_{r \; \in \; K_N(H) } \; P( R_n/N = r \; | \; H_n = H ) \; P( H_{n+1}= G \; | \; R_n/N = r, H_n=H ).
\end{equation}
\hl{Recall th}at $card(K_N(H)) \leq [ (1+F_g)(N+1) ]^{g^2}$. 
This Markov chain generates random histogram trajectories $ \mbold{H} = [H_1\, H_2,\cdots ,H_T]$ of arbitrary duration $T$. In the following sections, we develop an explicit large deviations theory for histogram trajectories. Since realistic experiments on bacterial evolution involve large populations with \mbox{$5\times 10^5 \le N \le 10^8$}, our rare-event study naturally focuses on asymptotic results for $N \to \infty$.

%%%%%%%%%%%%%%%%%%%%%%%%%%%%%%%%%%%%%%%%%%%%%%
\section{Large Deviations for Daily Cycles}
\label{LDcycles}
For the time-homogeneous Markov chain $\{H_n\to H_{n+1}\}$, the distribution of histogram trajectories is essentially determined by its transition kernel $\mathcal{Q}_N $, given in \eqref{kernel}. A rare-event analysis for random histogram trajectories must hence  start with a large-deviation analysis of the daily transition kernel.  In this section, we carry out a large-deviation analysis for random mutations and for random selection, concluding with an explicit large-deviation framework for daily transitions.  
The distribution of our daily random-mutation matrices $R_n$ is well approximated by products of Poisson distributions, and daily random selections follow multinomial distributions. For large $N$ we combine  precise large-deviation results for Poisson and multinomial distributions   to derive accurate   approximations of the daily transition kernel.  Technically however, to obtain smooth formulas one has to  accurately quantify   how to relax the $N$-rationality conditions for the random histograms $H_n, J_n$ introduced above.

%%%%%%%%%%%%%%%%
\begin{Definition}
\label{balls}
For any histogram $H,$ define the ball $V_N(H)$ and the $N$-rational ball $B_N(H) \subset V_N(H)$ by $V_N(H)  =   \vset{ H' \in \mathcal{H} \;: \; \vnorm{H' - H} \leq \frac{2}{3 N}} $ and
$B_N(H)  =   \mcal{H}_N \cap V_N(H)$. 
\end{Definition} 
Note that $card(B_N(H)) \leq 2^g$. The following lemma  characterizes ``uniformity'' within  $B_N(H)$ and helps to relax $N$-rationality.

%%%%%%%%%%%%%%%%
\begin{Lemma} \label{basicBN}
Fix any histogram $H$. For $N$ large enough, all histograms $H'$ in the  N-rational ball  $B_N(H)$ have the  same support as $H$, and an  essential minimum $b(H')$ larger than $b(H)/2$. More precisely, these results hold as soon as  $N > \frac{2}{b(H)}$.
\end{Lemma}
%%%%%%%%%%%%%%%%
\begin{proof} 
For $H' \in B_N(H)$ and any $k \in \supp(H')$, one has $H'(k) \geq 1/N$ since $H'$ is $N$-rational. This implies 
\[
H(k) \geq H'(k) - 2/(3 N) > 1/(3 N) >0.
\]
\hl{Hence,} $k \in \supp(H)$ and therefore $\supp(H') \subset \supp(H)$.
Since $\vnorm{ H' - H} \leq \frac{2}{3 N}$, one has for $N > \frac{2}{b(H)}$ and any $j \in \supp(H)$, 
$$
H'(j) \geq H(j) - \frac{2}{3 N} \geq b(H) - \frac{2}{3 N} > \frac{2 b(H)}3.
$$
\hl{This inequality}  yields $b(H') \geq 2 b(H)/3$, which proves the lemma.
\end{proof}

%%%%%%%%%%%%%%%%%%%%%%%%%%%%%%%%%%%%%%%%%%%%%%%
\subsection{Random-Mutation Matrices and Products of Poisson Distributions}
\label{subs3.1}
Given the initial histogram $H_n =H$ for day $n$,  the conditional distribution of the random-mutation matrix $R_n$ was defined in two steps (see Section~\ref{compMat}). First, one defines a companion random matrix $Z_n$ of independent random variables $Z_n(j,k)$, where each $Z_n(j,k)$  is Poisson-distributed with mean $m q_{j,k} N F_j H_n(j)$. Then, for  $r \in K_N(H)$ one defines  $P(R_n/N = r \; | \; H_n=H  )$ by the ratio 
$$
\frac
{P(Z_n/N = r \; | \;H_n=H)}
{P (Z_n/N \in K(H) \; | \;H_n=H )}.
$$

\hl{In this section} we prove Theorem~\ref{companionTH}, which states that as $N\to +\infty$, the random-mutation matrix $R_n$ and its Poisson companion $Z_n$ (see Section~\ref{compMat}) have nearly identical conditional distributions given  $H_n = H$. A key point is to prove that the conditional probability $P( Z_n/N \in K(H) \; | \; H_n = H ) $ tends to 1 as exponential speed $\simeq d(H)^N$, where $0 <d(H) <1 $ is given by the following definition.
%
%%%%%%%%%%%%%%%%
\begin{Definition} 
\label{d(H)}
The essential minimum $b(H) >0$ of any histogram $H$ is given by Equation~\eqref{supp(H)b(H)}. 
Define the \emph{decay coefficient} $0<d(H)<1$ of any histogram $H$ by $d(H) = \exp( - \delta F_1  b(H) )$, where $\delta = \log(1/m) -1 > 0$.

For $H \in \mathcal{H}_N$, one has  $b(H) \geq 1/N$ and hence $d(H)^N \leq e^{- \delta \; F_1} \ll 1$ for all $H \in \mcal{H}_N$ and all $N$. 
In actual laboratory experiments~\cite{len94a,azencoop,Gordo2012,Desai2007}, one typically has $F_1 > 100$ and  mutation rates $m \leq 10^{-6},$ so that 
$\delta F_1 > 12800$ and hence $d(H)^N < 10^{- 5500}$, which is practically zero.
\end{Definition}

%
%%%%%%%%%%%%%%%%
\begin{Theorem} \label{companionTH}
For large $N$, the random-mutation matrix $R_n$ on day $n$ and its Poissonian companion $Z_n$  have nearly identical conditional distributions given $H_n = H$. Fix any histogram $H^*$ and the N-rational ball $B_N(H^*)$. Let $c = 2 +\log(g) / F_1 $. 
Then, provided $N > \frac{c}{b(H^*)} $, one has for all $n$, all $H \in B_N(H^*)$, and all matrices $z \in K(H),$
\begin{equation} \label{control}
1 \leq  \frac{P( R_n/N = z \; | \; H_n = H )} {P( Z_n/N = z \; | \; H_n = H )} \leq 1 + 2 g\,  d(H^*)^{N/2}.
\end{equation}
\end{Theorem}
This result can be reformulated in terms of the conditional density  $f_N(z, H) =\frac{P( R_n/N = z \; | \; H_n = H )}{P( Z_n/N = z \; | \; H_n = H )}$ of $R_n/N$ with respect to $Z_n/N$. 

Fix any $0< a < 1$. Then, $f_N(z, H) \to 1$ at exponential speed as $N \to \infty$ and the convergence rate is uniform for $b(H^*) \geq a$, $ H \in B_N(H^*)$, and $z \in K(H)$.  Therefore, the conditional joint distribution of the mutation matrix $[R_n(j,k)]$, given $H_n \in B_N(H^*)$, becomes, for large $N,$ practically equal to the product of  Poisson distributions with respective means $m q_{j,k} N F_j H_n(j)$.

To prove Theorem \eqref{companionTH}, we  first need  two lemmas.

%%%%%%%%%%%%%%%%
\begin{Lemma}[Poisson Large Deviations] 
\label{poissonLD}
Fix any $u >0,$ and let $X$ be a random variable having a Poisson distribution with mean $N u$. For any $v > u > w > 0$ and for any integer $N,$ the following inequalities hold:
\begin{align} 
P( X \geq N v )  &\leq  \exp (- N \left[ u + v \log (v / u) - v) \right]), \label{cramerN1}\\
P( X \leq N w )  &\leq  \exp (- N \left[ u + w \log (w / u) - w) \right]). \label{cramerN2}
\end{align}
\hl{Moreover, }for any $v \geq 0$ such that $N v$ is an integer, one has 
\begin{equation}\label{XNv}
\frac{1}{N}\log{P( X= N v )} = - (u + v \log(v/u) - v) + o(1) 
\end{equation}
with $| o(1) | \leq 2 \log(N)/N $.
\end{Lemma}
%%%%%%%%%%%%%%%%%
\begin{proof}
Express the Poisson random variable $X$ as 
$X = \sum_{s=1}^N  X_s$, where $X_1,\cdots, X_N$ are i.i.d Poisson random variables with the same  mean $u >0$. Then,  for all $N$ and $v > u$ (see, e.g.,~\cite{azenLDT2012,dembo}), the empirical mean $X/N$ of the $X_i$ must  verify the large-deviation inequality
$
P ( X/N \geq  v) \leq \exp(- N \lambda(v))
$,
where $\lambda(v)$ is the large-deviation rate for the Poisson distribution, and is given by the well-known formula $\lambda(v) = u + v \log (v / u) - v$ for all $v >u >0$. This proves \eqref{cramerN1} and a similar argument proves \eqref{cramerN2}.

The classical proof of the Stirling formula  for factorials can be easily modified to yield the following uniform inequality:
\begin{equation}\label{stirling}
\abs{\log(N!) - N \log(N/e)} \leq 2 \log(N) \quad \text{for all  } N \ge 1.
\end{equation}
\hl{For any } $v \geq 0$ such that $N v$ is an integer, one has $P(X= N v) = e^{- N u} (N u)^{N v} / (N v) !$, so that
$$
\frac{1}{N}\log{P(X= N v)} = - u + v \log(N u) - \frac{1}{N} \log ((N v) !).
$$
\hl{For} $v > 0$, \eqref{stirling} implies
$$
\frac{1}{N} \log ((N v)!)= v \log(N) + v \log(v) -v + o(1)
$$
with $\abs{o(1)} \leq 2\log{(N)} / N$. The last two equations prove \eqref{XNv} for $v >0$. Finally, \eqref{XNv} is trivially true for $v =0$ with the convention $0 \log(0) = 0$. This concludes the lemma.
\end{proof}
%

%%%%%%%%%%%%%%%%
\begin{Lemma} \label{sumexp} 
Fix a positive sequence $\epsilon(N)$ such that $\epsilon(N) \to 0$ as $N \to \infty$. Fix any set $E$ and any function $\lambda(s) \geq 0$ defined for all $s \in E$. Fix $c > 0,\beta > 0,$ and let $E_N$ be a finite subset of $E$ with $card(E_N) \leq c N^\beta$. Consider fast-vanishing exponentials $p_N(s) > 0$ indexed by $s \in E_N$ such that 
$
\frac{1}{N} \log{p_N(s)} = - \lambda(s) + o_s(1)
$, 
where $\abs{o_s(1)} \leq \vep(N)$.  Define $\Lambda(E_N) = \inf_{s \in E_N} \lambda(s) $.  
The sum $p(E_N) = \sum_{ s \in E_N } p_N(s) $ satisfies
$$
\frac{1}{N}\log {p(E_N)} = - \Lambda(E_N) + o(1) 
$$
for all $N$ with $|o(1)| \leq \epsilon(N) + \beta \log(N) / N + \log(c) /N.$ 
\end{Lemma}
%%%%%%%%%%%%%%%%
\begin{proof}
Select $s(N) \in E$ such that $\Lambda(E_N) = \lambda(s(N))$. This yields the lower bound 
\begin{equation} \label{low}
\frac{1}{N}\log {p(E_N)} \geq \frac{1}{N}\log {p(s(N)} = - \lambda(s(N)) + o_{s(N)}(1) \geq \; - \Lambda(E_N) - \epsilon(N).
\end{equation}
\hl{We have} $ p_N(s) \leq \exp (  - N [ \Lambda(E_N) + \epsilon(N) ]  )$ for all $s \in E_N$ by definition of $\Lambda(E_N) $.  This gives $ p(E_N) \leq \exp (  - N [ \Lambda(E_N) + \epsilon(N) ] ),$ which implies
$$
\frac{1}{N}\log {p(E_N)} \leq - \Lambda(E_N) + \epsilon(N) + \log(card(E_N)) / N.
$$
\hl{Combining th}is upper bound with \eqref{low} concludes the proof.
\end{proof}
%%%%%%%%%%%%%%%%

%%%%%%%%%%%%%%%%%
\begin{proof}[Proof of Theorem \ref{companionTH}]
\hl{Consider} 
any $N$-rational histogram $H$. Given $H_n = H$, let $Z_n$ be a matrix of independent Poisson random variables $Z_n(j,k)$ with means $m N q{j,k} F_j H(j)$. The sums $S(j) = \sum_{k= 1}^g  Z_n(j,k) $ then have Poisson distributions with respective means $\mbb{E}[S(j)] = N m F_j H(j)$.
First, we show that, for each $j,$
\begin{equation} \label{controlCj}
P (Z_n \in N  K(j, H) \; | \; H_n = H ) \geq 1- d(H)^N 
\end{equation}
holds for all $N$ and $H$.  For this, we separately consider the two cases $H(j) > 0$ and $H(j) = 0$.\smallskip\\
{\bf \hl{Case 1:}
} Suppose $H(j) > 0$.  For any $s > m$, apply \eqref{cramerN1} to $X = S(j)$, with $v = s F_j H(j)$ and $u = m q(j,k) F_j H(j)$, to obtain
\begin{equation} \label{tail.j}
P( S(j) \geq s N F_j H(j) \; | \; H_n = H) \leq \exp \left( - s N F_j H(j)[\log(1/m) -1] \right).
\end{equation}

\hl{Since }$ \log(1/m) -1 = \delta >0$, the inequality \eqref{tail.j} implies 
\begin{equation} \label {tail.j2}
P( S(j) \geq s N F_j H(j) \; | \; H_n = H) \leq \exp \left(- \delta s  N  F_j H(j) \right) \leq d(H)^{s N}
\end{equation}
for all $N, H$ and $s > m $.  By definition \eqref{Kj} of $K(j,H),$ 
$$P (Z_n / N \in K(j, H) \; | \; H_n = H ) = 1 - P(S(j) \geq N F_j H(j) \; | \; H_n = H).$$
\hl{Equation }\eqref{tail.j2}, with $s=1$, implies \eqref{controlCj} for all $N$ and $H$.
 \smallskip\\
{\bf \hl{Case 2:}} Suppose $ H(j) = 0.$ This implies  $P(Z_n(j,k)=0 \:| \:H_n= H ) = 1$ for all $k.$ Therefore, $S(j) = 0$ and $P( Z_n \in N K(j, H_n) \; | \; H_n = H ) =1$, which trivially satisfies \eqref{controlCj}.
\smallskip\\
Now, since $K(H) = \bigcap_{j = 1}^g  K(j, H)$ and \eqref{controlCj} holds for all $N$ and $H,$ then
\begin{equation} \label{control2}
P (Z_n / N \in K(H_n) \; | \; H_n = H) \geq 1- g \,d(H)^N
\end{equation}
also holds for all $N$ and $H$. The constraint $ N > N(H) = \log(2 g) / (  12  F_1 b(H) ) $ forces $g \,d(H)^N < 1 / 2$, so that $1/(1 - g \, d(H)^N) \leq 1+2 g \, d(H)^N $. For $N > N(H)$, Equation~\eqref{control2} then implies
\begin{equation} \label{controlC}
1 \leq \frac{1}{P( Z_n / N \in K(H) \; | \; H_n = H)} \leq 1 + 2 g \, d(H_n)^N.
\end{equation}

\hl{For all matrices} $z \in K(H)$, we have
$$
P (R_n/N = z \; | \; H_n = H) = \frac{P ( Z_n / N = z \; | \; H_n =H )} {P ( Z_n/N \in  K (H) \; | \; H_n = H)}.
$$
\hl{For} $H_n= H$, $z \in K(H)$, and $N > N(H)$, Equation~\eqref{controlC} yields

\begin{equation} \label{compar}
\begin{aligned}
P ( Z_n/N = z \; | \; H_n =H) & \leq P ( R_n/N = z \; | \; H_n = H)  \\
& \leq \left(1+ 2 g \, d(H)^N\right)  P ( Z_n/N = z \; | \; H_n = H ). 
\end{aligned}
\end{equation}

Consider any fixed $H^* \in \mathcal{H}$ and any $H_n = H \in V_N(H^*)$, which is equivalent to $H \in B_N(H^*)$ since $H_n$ is $N$-rational. By Lemma \ref{basicBN}, for $N > \frac{2}{b(H^*)}$, 
$\supp(H) = \supp(H^*) $ and $b(H) \geq b(H^*)/2$, so that $d(H)< d(H^*)^{1/2}$.  Let $N_1 = \frac{c}{b(H^*)}$, with $c = 2 + \frac{\log(g)}{F_1}$. Then, $N > N_1$ will force $N > N(H)$ provided $H \in B_N(H^*)$. Equation~\eqref{compar} then implies, for all $N$-rational $H \in B_N(H^*)$, $z \in K_N (H)$, and $N > N_1$,
\begin{equation} \label{N3}
1 \leq  \frac{P( R_n/N = z \; | \; H_n = H )} {P( Z_n/N = z \; | \; H_n = H )} \leq 1 + 2 g \, d(H^*)^{N/2}, 
\end{equation}
which proves \eqref{control}. The uniform convergence in the statement of the theorem is an easy consequence of \eqref{control}.
\end{proof}

%%%%%%%%%%%%%%%%%%%%%%%%%%%%%%%%%%%%%%%%%%%%%%%
\subsection{Large Deviations for Random-Mutation Matrices $R_n$}
\label{subs3.2}
In this section, we develop explicit   large-deviation formulas for  the daily matrices $R_n$ of random mutations. We have shown that  with probability tending to 1 at exponential speed for large $N$, 
the random variables $R_n(j,k)$ are conditionally independent  given $H_n=H$, and that their   conditional distributions $\pi_{H, i,j}$ are Poissonian, with means depending on $H$ and $j,k$. 
This indicates that, given $H_n=H$, the conditional large-deviation rate function of $R_n$ should be computed by summing over all $j,k$ the Poisson rate functionals of $\pi_{H, i,j}$. However, proving this requires a few meticulous technical steps.

Since the non-zero coefficients of $R_n$ are all of order $N$, the key matrices to obtain large-deviation results are the matrices  $r_n=R_n/N$, which belong to the following convex cone~$\mcal{R}$. 
%%%%%%%%%%%%%%%%
\begin{Definition}
The random matrices $r_n=R_n/N$ always belong to the convex cone $\mcal{R}$  of all $g\times g$ matrices $r$ such that   $r_{j,k} \geq 0$ and $r_{j,j} = 0$ for all $j,k$. For $r\in \mcal{R}$, denote $\vnorm{ r } = \max_{j,k} \abs{ r_{j,k} } $ the $L_{\infty}$ norm of $r$.

In $\mcal{R}$, determining  relevant ``bounded'' sets $S$  of matrices $r$ involves controlling the \emph{smallest positive} $r_{j,k}$ over $r \in S$.   So for $r\in \mcal{R}$, we define  the \emph{support} $\supp(r)$ and the \emph{essential minimum} $b(r) >0$ of $r$ by
$$
\supp(r) = \{ (j,k) \; | \; r_{j,k} >0 \} \qquad \text{and} \qquad b(r) = \min_{(j,k) \, \in \, \supp(r)} \; r_{j,k}. 
$$
\hl{As} above for histograms, denote by $V_N(r)$  the $L_{\infty}$ ball with center $r\in \mcal{R}$ and  radius $\frac{2}{3 N}$. Let $B_N(r) \subset V_N(r)$ be the set of $N$-rational matrices in $V_N(r)$. Note that $card( B_N(r) ) \leq 2^{g^2}$.
\end{Definition}

Given $H_n = H$, we showed above that the random matrices $ R_n / N$ always belong  to the set $ K_N(H)$  of all N-rational matrices within the convex set $K(H) \subset \mcal{R}$ defined by Definition~\ref{KNconstr}. Note that for $r \in K(H)$ and $(j,k) \in \supp(r),$ one must have $H(j) >0$ and $q_{j,k} > 0$.  

We can now define the adequate large-deviation rate function controlling  the probabilities of rare random-mutation matrices. 
%%%%%%%%%%%%%%%%%
\begin{Definition}[Rate function for random-mutation matrices]
\label{def:mut}
 For any histogram $H\in \mcal{H}$, any matrix $r \in K(H),$ and any $j,k$ such that $H(j) q_{j,k} >0,$ denote by $\pi(H,j,k)$ the Poisson distribution with mean $m\, q_{j,k} F_j H(j)$. Note that for large $N$ $\pi(H,j,k)$ will be extremely close to the distribution of $R_n(j,k)$. Well-known large-deviation formulas for the Poisson distribution show that  the rate function $ L_{j,k} (H,r) \geq 0$ for $\pi(H,j,k)$ is   given by 
 \begin{equation} \label{Ljk}
L_{j,k}(H,r)  = 
m \,q_{j,k} F_j H(j) + r_{j,k} \log\left( \frac{r_{j,k}} { e \,m \,q_{j,k} F_j H(j) }\right) \;\; \text{when} \;\; H(j)\, q(j,k) >0.  
\end{equation}
\hl{When} $H(j)\, q(j,k) = 0$, we  conventionally set $L_{j,k}(H,r)  = 0$.
For Poisson distributions such as $\pi(H,j,k)$, the rate function  $ L_{j,k}(H,r) $ is a \emph{strictly convex} function of $r_{j,k}$ for $r_{j,k} >0$, and one has $L_{j,k}(H,r)=0$ if and only if $r_{j,k}= m\,q_{j,k} F_j H(j)$. Denote by $L(r,H)$ the $g \times g$ matrix  $L_{j,k}(H,r)$.

We show further on  that large deviations for  the random-mutation matrices $r_n= R_n/N$ are controlled by the \emph{rate function} $mut(H,r)$,  given by
\begin{equation}\label{mutcramer}
\begin{aligned}
mut(H,r) & = \sum_{j,k} \; L_{j,k}(H,r) \\
& = \sum_{(j,k) \; | \; m\, q_{j,k} H(j) > 0} 
m\, q_{j,k} F_j H(j) + r_{j,k} \log r_{j,k} - r_{j,k} \log(e\, m\, q_{j,k} F_j H(j) ) .
\end{aligned}
\end{equation}
\end{Definition}
\begin{Proposition} \label{mutproperties}
 For all histograms $H$, the mutation rate function $mut(H,r)$ is a finite, non-negative, \emph{continuous, strictly convex} function of $r \in K(H)$. Moreover, $mut(H,r)$ is a continuous, convex function of $H$ for fixed $r$. For any fixed $0<a<1$, there is an explicit  constant $c(a)$ depending  only on $a$ and on the process parameters $g$, $F$, $Q$, $m$, such that   
\begin{equation} \label{mutbound}
mut(H,r) \leq c(a) \;\; \text{for all $H$ with $b(H) \geq a $ and all $r \in K(H),$}  
\end{equation} 
\hl{Note also tha}t  $mut(H,r) = 0$ if and only if $r_{j,k} = m\, q_{j,k} F_j H(j)$ for all $j,k$.
\end{Proposition} 
\begin{proof}
As with all Poisson rate functions,  $L_{j,k}$ must be  strictly convex in $r_{j,k}$. This proves  strict convexity in $r$    for $mut(H,r)$. Each $L_{j,k}$ is also clearly convex in $H$. This implies convexity in $H$ for  $mut(H,r)$. To prove \eqref{mutbound}, note the elementary inequalities 
\begin{align*} 
 x \,\abs{ \log x } \leq \mcal{D}(v) \; \text{for all } \;\; 0 \leq x \leq v;\\
 \abs{ \log x } \leq \mcal{E}(u,v) \; \text{for all } \;\; 0 < u \leq x \leq v.
\end{align*}
where $\mcal{D}(v)=1/e + \log(1+v),$ and $\mcal{E}(u,v)=\log(1+1/u)+\log(1+v)$.

For $H(j)\, q_{j,k} > 0$, $r \in K(H)$, and $b(H) \ge a$, one has 
$$
r_{j,k} \leq F_j H(j) \leq F_g\;\text{and} \;\; 0 < a \, m\, b(Q) \leq m\,q_{j,k} F_j H(j)) \leq F_g m \vnorm{Q}.
$$  

Combining the last three inequalities yields 
\begin{equation}
\begin{aligned}
\abs{ r_{j,k} \log(r_{j,k}) } & \leq\mcal{D}(F_g), \\
\abs{ r_{j,k} \log(e\, m\, q_{j,k} F_j H(j)) } & \leq F_g \abs{\log(e\, m\, q_{j,k}) }+ F_j H(j) \abs{ \log(F_j H(j)) } \\
& \leq \mcal{E}\left(a \, m\,b(Q) , m\,F_g  ||Q|| \right) + \mcal{D}(F_g).
\end{aligned}
\end{equation}
 \hl{This conclud}es the proof due to formula \eqref{mutcramer}.
\end{proof}

To obtain uniform H\"{o}lder continuity bounds for $mut(H,r)$, the histogram $H$ must belong to one of the  compact sets  $\mathcal{H}(a)$ defined next.  

%%%%%%%%%%%%%%%%
\begin{Definition}
For any fixed  $0 < a < 1,$ define the compact set of histograms $\mcal{H}(a) \subset \mcal{H}$ by
\begin{equation} \label{HIST(a)}
\mcal{H}(a) = \vset{H \in \mathcal{H} \; | \; b(H) \geq a}. 
\end{equation}
\end{Definition}
\begin{Proposition} \label{propholdermut}
Fix $0 < a < 1$ and any H\"{o}lder exponent $0 < \alpha < 1$. For all histograms $H', H \in \mcal{H}(a)$ with $\supp(H') = \supp(H)$, and all matrices $r\in K(H)$, $r'\in K(H')$, the mutation rate function $mut(H,r)$ verifies  
\begin{equation}\label{holdermut}
\abs{ mut(r', H') - mut(r, H)} \leq c \left( \vnorm{r' - r}^{\alpha} + \vnorm{ H'- H}\right)
\quad \text{with~~} c = \frac{ 14 g^2 F_g \log(F_g) }{a (1- \alpha)}.
\end{equation} 
\end{Proposition}
%%%%%%%%%%%%%%%%
\begin{proof}
The proof requires several meticulous inequalities and is presented in  \linebreak \mbox{Appendix~\ref{sec:ap3}}.
\end{proof}

%%%%%%%%%%%%%%%%

We now conclude the analysis of random mutations with a key asymptotic large-deviation result for random-mutation matrices.
%%%%%%%%%%%%%%%%%
\begin{Proposition} \label{LDmut}
Let $R_n$ be the random matrix of mutations on day $n$. Let $mut(H,r)$ be the mutations rate function defined by \eqref{mutcramer}. Fix $a >0$ and the parameters $\calP$ of the bacterial evolution model (Definition~\ref{P}). Then, there is a constant $N_0 = N_0(a, \calP)$ such that for $N > N_0$, the large-deviation formula
\begin{equation} \label{LDmut1}
\frac{1}{N} \log P( R_n / N = r  \; | \;  H_n = H )  =  - mut(H,r) + o(N),
\end{equation}
with  $\abs{ o(N) }  \leq  4 g^2  \log N / N $ will hold uniformly for all $H\in \mcal{H}(a) \cap \mcal{H}_N$ and $r \in K_N(H).$ 
\end{Proposition}
%%%%%%%%%%%%%%%%%
%%%%%%%%%%%%%%%%%%%%%%%
\begin{proof}
The proof of this proposition is presented in Appendix \ref{sec:ap3}.

\end{proof}
%%%%%%%%
%%%%%

%%%%%%%%%%%%%%%%%%%%%%%%%%%%%%%%%%%%%%%%%%%%%%%
\subsection{Large Deviations for Random Selection} \label{LDsel}
After daily growth and mutations, the current population $POP_n$ has a random size much larger than $N$, and the day $n$ cycle ends with the random selection of a sample of size $N$. This sample becomes $pop_{n+1}$ and has a histogram $H_{n+1}$. Given $H_n=H$ and $R_n=R$, $POP_n$ has large size $\eta(H) N$, and the conditional distribution $mult_{H,R}(N)$ of $N H_{n+1}$ is  multinomial. Since both population and sample have very large sizes of order $N$, one cannot directly apply  the classical large deviations for  ``standard'' multinomial distributions. Moreover, the histogram of $POP_n$ given $H_n, R_n$ is itself random. Therefore, in this section we develop an explicit rate function for large deviations of the conditional multinomial  $mult_{H,R}(N)$, and demonstrate that this functional involves the well-known Kullback--Leibler divergence. These results then have  to be transposed into precise rate function formulas for conditional large deviations of $H_{n+1}$ given $H_n=H$ and $R_n=R$. We also analyze the regularity of these conditional rate functions when  $H$ and $R/N$ vary, which is an important technical point to prove large-deviation results for the daily transition from $H_n$ to $H_{n+1}$.

%%%%%%%%%%%%%%%%%
\begin{Definition}[Kulback--Leibler divergence]
\label{def:kl}
The classical {Kullback--Leibler divergence} between two histograms $G$ and $J$ is defined by 
\begin{equation} \label{KL}
KL(G,J)  = \sum\limits_{j \in \supp(G)}  G(j) \log \frac{G(j)}{J(j)} \;\; \text{ when }\;\; \supp(G) \subset \supp(J)
\end{equation}
and by $KL(G,J) = +\infty$ otherwise. 
\end{Definition}
Note that $KL(G,J)\geq 0$, and that $KL(G,J) = 0$ if and only if $G = J.$ 
When $ \supp(G) \subset \supp(J), $ the function $KL(G,J)$ is finite, continuous, and is  strictly convex  in $J$ for fixed $G$, as well as strictly convex in $G$ for fixed $J$, since $KL(G,J)$ has the following  partial derivatives  for all $i,j \in \supp(G) \subset \supp(J)$:
\begin{align}
\partial_{G(i)} KL(G,J)  &=  1 + \log (\frac{G(i)}{J(i)}) \;\; 
\text{ and } \;\;\; 
\partial_{J(i)}  KL(G,J)  =  - \frac{G(i)}{J(i)} 
\label{partialKL},\\
\partial_{G(i)} \partial_{G(j)}  KL(G,J)  &= \left(\frac{1_{\vset{ i  =  j }}}{G(j)}\right)
\partial_{J(i)} \partial_{J(j)} KL(G,J)  = 1_{\vset{ i  =  j }} \frac{G(i)}{J(i)^2}.
\end{align}

On day $n$, we denote the random  population histogram right after growth and mutations (Phases 1 and 2) as $J_n$. For any two histograms $G$ and $J$ such that $\supp(G) \subset \supp(J)$, Equation~\eqref{multinomial2} shows that 
\begin{equation}\label{temp1}
 P(H_{n+1} = G \; | \; J_n = J) = \mu_{N, J}(N G)
\end{equation}
where $\mu_{N,J}(V)$ is the  multinomial distribution defined  by \eqref{multinomial}.
For a large $N$ we now derive an asymptotic expression of $\frac{1}{N} \log(\mu_{N,J}(NG))$.
%%%%%%%%%%%%%%%%
\begin{Proposition} \label{LDmultinomial}
For any histogram  $J\in \mcal{H}$, and any $N$-rational histogram $G\in \mcal{H}$ such that $\supp(G) \subset \supp(J)$,  the multinomial distribution $\mu_{N,J}$ defined by \eqref{multinomial} verifies
\begin{equation} \label{KLlogmu}
\frac{1}{N} \log(\mu_{N,J}(N G) )= - KL(G,J) + o(1) 
\end{equation}
with \emph{uniform} remainder $\abs{ o(1) } \leq 2 (g+1) \log N /N. $
\end{Proposition}
%%%%%%%%%%%%%%%%
\begin{proof}
The proof is presented in Appendix
\ref{sec:ap3}.
\end{proof}
We can now prove the main large-deviation evaluation for random selections. Combined with Proposition~\ref{LDmut}, this result will be used to obtain the large-deviation estimate for the full daily transition kernel $P(H_{n+1} =G \;|\; H_n= H)$.
%%%%%%%%%%%%%%%%
\begin{Proposition} \label{LDKL}
Fix $0< a <1$, and set $N^* = 2 / a + \log(F_g)$. On day $n$, let $J_n $ be the population histogram at the end of the mutation phase (Phase 2). For any $G$ with $b(G) > a$, any  $G'$ in the N-rational ball $ B_N(G)$, and any $N > N^*$, one has the large-deviation estimate
\begin{equation} \label{LDselect1}
\frac{1}{N} \log P ( H_{n+1} = G' \; | \;  H_n  ,  R_n  ) = - KL( G, J_n ) + o(1) 
\end{equation}
with uniform remainder $\abs{ o(1) } \leq (5g + 1) \log N / N$.
\end{Proposition} 
%%%%%%%%%%%%%%%%
\begin{proof}
The proof is presented in Appendix
\ref{sec:ap3}.
\end{proof}

%%%%%%%%%%%%%%%%%%%%%%%%%%%%%%%%%%%%%%%%%%%%%%%
\subsection{Large-Deviation Estimates for the One-Step Transition Kernel}
\label{subs3.4}
For each random phase of any day $n$  cycle, namely, mutations and selection,  we have so far proved precise \emph{conditional} large-deviation estimates. We can now combine these two results to develop a large-deviation analysis for the one-step transition kernel $\mathcal{Q}_N (H, G) = P(H_{n+1} = G \; | \; H_n = H)$.  In this section, we prove that, in the state space $\mathcal{H}$  of population histograms, the Markov transition kernel $\mathcal{Q}_N (H, G)$ converges to a deterministic kernel at exponentially fast speed as $N \to \infty$, and we explicitly compute  the corresponding rate function  $C(H,G) \ge 0$ as the limit $ = \lim_{N\to \infty} \frac{1}{N}\log(\mathcal{Q}_N(H, G))$.

Equations~\eqref{mutcramer} and \eqref{LDselect1} provide two explicitly computed ``partial'' rate functions, namely, the rate function $mut(H,r) \geq 0$ to approximate  $- \frac{1}{N} \log(P(R_n/N = r \; | \;H_n = H))$ and the Kullback--Leibler rate function $KL(G, J) \geq 0$ to approximate $- \frac{1}{N} \log(P(H_{n+1} = G \; | \; J_n = J))$ 
(where $J_n$ is the histogram after the growth and mutation phases).
The deterministic function $\Psi(H,r)$ given by \eqref{psi} provides the approximation  $J_n \approx \Psi(H_n,R_n/N)$ with accuracy $\approx 1/N$. We thus expect $KL(G, \Psi(H,r))$ to be the rate function approximating $- \frac{1}{N} \log(P(H_{n+1} = G \; | \; H_n = H, R_n/N = r )$.

We will define an explicit  composite transition rate $\tau(H,r,G) \ge 0$, which  for large $N$ provides an  estimate of the form 
\[
\tau(H,r,G) \approx - \frac{1}{N}\log(P(\{H_{n+1} = G\} \, \& \,\{R_n/N = r\} \; | \; H_n = H)).
\]
\hl{Hence, fo}r  fixed $H,G$, the probability  $P(\{H_{n+1} = G\} \, \& \,\{R_n/N = r\} \; | \; H_n = H))$ will be nearly maximized over $r$ when the random-mutation matrix $R_n$ is close to $N r^*$, where $ C(H,G) = \tau(H,r^*,G)$ is the minimum of $\tau(H,r,G)$ over $r$. Later on, we show that $C(H,G)$ is the rate function approximating $-\frac{1}{N}\log(P(H_{n+1} = G \;|\; H_n =H))$.

We now present a rigorous development of this large-deviation analysis. 
%%%%%%%%%%%%%%%%
\begin{Definition}[Composite Transition Rate]
To control large deviations for the composite daily transition $(H_n = H) \to (R_n/N = r) \to (H_{n+1} = G),$ we introduce the \emph{composite transition rate} $\tau(H,r,G) \geq 0$, defined by
\begin{equation} \label{tauHrG}
\tau(H,r,G) = mut(H,r) + KL(G, \Psi(H,r)) 
\end{equation}
for $H, G \in \mcal{H}$ and $r \in K(H).$ 
%\hfill \qedsymbol{}
\end{Definition}

Recall that $mut(H,r)$ is a finite continuous convex function of $r \in K(H)$ (see Proposition~\ref{mutproperties}), and  $KL(G,J)$ is a continuous convex function of $J \in \mathcal{H}$, with strict convexity whenever $KL(G,J)$ is finite (see def. \ref{def:kl});  moreover, due to \eqref{psi}, $J = \Psi(H,r) $ is an affine function of $r$, so that   $KL(G,\Psi(H,r))$ is continuous convex in $r$ on the compact convex closure of $K(H)$.

Hence,  $\tau(H,r,G) \geq 0$ is continuous and convex in $r$ on the compact convex closure of $K(H)$, with strict convexity in $r$ whenever $\tau(H,r,G)$ is finite. Similarly, one verifies that whenever the function  $\tau(H,r,G)$ is finite, it is convex in $H$ for fixed $(r,G)$, and convex in $G$  for fixed $(H,r)$.

Since $mut(H,r)$ is finite, one can have    $\tau(H,r,G) = + \infty$ if and only if  $ KL(G, \Psi(H,r)) = + \infty$, which is equivalent to the existence of a genotype $j$ such that $G(j) > 0$ and $ \Psi_j(H,r) = 0$. Due to \eqref{zeroPsi}, we see that for $r \in K(H),$ one has $\tau(H,r,G) = + \infty$ if and only if there is a genotype $j$ such that $G(j) > 0$ and $H(j) = r_{k,j} = r_{j,k}= 0$ for all $k$.

%%%%%%%%%%%%%%%%%
\begin{Definition}[Feasible Transitions and One-Step Cost]
\label{def:onestep}
By construction, for any $H, G \in \mcal{H}$, the transition kernel $P(H_{n+1} = G \; | \; H_n = H)$ is strictly positive if and only if for any $j$ such that $G(j) > 0$ and $H(j) =0$, one can find a $k$ such that $H(k) \, q_{k,j} >0$. We then say that $(H \to G)$ is a \emph{feasible transition}.  Let $\tau(H,r,G) \geq 0$ be the composite transition rate defined by \eqref{tauHrG}. For any $H,G \in \mathcal(H)$, define the {one-step cost function} $C(H,G) \geq 0$ by 
\begin{equation} \label{costHG}
C(H,G) = \min_{r \in K(H)} \tau(H,r,G) = \min_{r \in K(H)} \left[ mut(H,r) + KL(G, \Psi(H,r) ) \right].
\end{equation}
\end{Definition}
{\hl{An explicit }expression of the one-step cost $C(H,G)$ is developed later in this section in Theorem~\ref{thCHG}. }

%%%%%%%%%%%%%%

\begin{Lemma} [Finiteness and convexity of~$C(H,G)$]
\label{finitecost}
For all  $H, G\in \mcal{H}$, the one-step cost function $C(H,G)$  is finite if and only if $(H \to G)$ is a feasible transition.

When $C(H,G)$ is finite, there is a unique $\rho$ in the closure of $K(H)$ such that $C(H, G) = \tau(H, \rho, G)$; moreover,  $C(H,G)$ is then convex in $H$ for fixed $G$ and convex in $G$ for fixed $H$.
\end{Lemma}
%%%%%%%%%%%%%%%%
\begin{proof} 
By definition, $C(H,G)$ is finite if and only if there is at least one $r \in K(H)$ with $\tau(H,r,G) $ finite.  As seen above, this occurs if and only if for each $j$ such that $G(j) >0$ and $H(j) = 0,$ there is a $k$ with $r_{k,j} >0$.  For any $(j,k)$ verifying such a condition, one must also have $m\, q_{k,j} >0$ due to \eqref{Kj} so that $(H \to G)$ is a feasible transition.

Since $\tau(H,\rho,G)$ is continuous in $\rho$ on $\overline{K(H)}$, there is at least one $\rho \in \overline{K(H)}$ such that $C(H, G) = \tau(H, \rho, G)$.
When $C(H,G)$ is finite, this $\rho$ is unique in the compact, convex set $\overline{K(H)}$ since $\tau(H, \rho, G)$ is strictly convex in $\rho$ whenever $\tau(H, \rho, G)$ is finite; the asserted convexity properties of $C(H,G)$ follow then from the separate convexity in $H$ and $G$ of $\tau(H, \rho, G)$.
\end{proof}
%%%%%%%
A key large-deviation result (see Theorem \ref{transLD} 
below) will  relate $C(H,G)$ to the one-step Markov transition kernel $\mcal{Q}(H,G)$  by a large-deviation estimate of the form $1/N log(\mcal{Q}(H,G)) \approx -\, C(H,G)$, valid for large $N$. The proof of Theorem \ref{transLD}  requires precise  H\"{o}lder continuity estimates for $KL(H,G)$ and $\tau(H,r,G)$, which we now state in the next two technical results.

%%%%%%%
%%%%%%%%%%%%%%%%
\begin{Lemma} [Continuity estimates for $KL(G,J)$]
\label{lipKL}
For  $G,J\in \mcal{H}$, then $KL(G,J)$ is finite if and only if $\supp(G) \subset \supp(J)$, and one then has
\begin{equation} \label{betabound}
KL(G,J) \leq \log \beta(G,J) \leq \frac{KL(G,J) + \log(g)}{b(G)} \;\text{where}\;\;\beta(G,J) = \max_{k \in \supp(G)}  \frac{1}{J(k)}.
\end{equation}
\begin{enumerate}[label = (\roman*),leftmargin=*,align=left]
\item For $G, J, J'\in \mathcal{H}$ with  both $KL(G,J)$ and $KL(G,J')$ finite, we have  
\begin{equation}\label{lipKLJI}
\abs{KL(G, J) - KL(G, J')}  \leq  c \vnorm{ J - J' } 
\end{equation}
with $ c = g^{1 + 1/b(G)} e^{\kappa / b(G)}$ and  $\kappa = \min\left( KL(G, J),  KL(G, J') \right)$.
\item
Fix $0 < a \leq 1/2$. Consider any  $G, G', J \in \mathcal{H}$ with  
$KL(G,J)$ and $KL(G',J)$ finite, and verifying $b(G)\geq a$ as well as  $b(G') \geq a$.  Then, one has
\begin{equation} \label{holderKLG}
\abs{KL(G',J) - KL(G,J)} \leq c \vnorm{ G' - G } ^{1/2}
\end{equation} 
with $c = \frac{g}{a} \left(2 + \log g  + \max \vset{ KL(G,J), KL(G',J)}\right)$. 
If, moreover, $\supp(G) = \supp(G')$, one has also the Lipschitz continuity 
\begin{equation} \label{lipKLG}
\abs{KL(G',J) - KL(G,J)} \leq c_1 \vnorm{ G' - G } 
\end{equation} 
with $c_1 = \frac{g}{a} \left(2 + \log g + \min\left(KL(G,J) , KL(G',J) \right)\right)$.
\end{enumerate}
\end{Lemma}
%%%%%%%%%%%%%%%%
\begin{proof}
The proof of this technical lemma is given in the appendix (see Appendix  \ref{sec:ap3}).
\end{proof}

%%%%%%%%%%%%%%%%
\begin{Proposition} [H\"{o}lder continuity for the composite transition cost $\tau(H,r,G)$]
\label{holdertau} 
Fix $A >0, a >0,$ and the process parameters $\calP$ (Definition~\ref{P}). For all $H', H, G\in \mathcal{H}$ with
$\min(b(H),b(H'),b(G)) > a$, 
and matrices $r' \in K(H')$ and $r \in K(H)$ verifying
\begin{equation} \label{hyp}
\supp(H') = \supp(H),\;\; \vnorm{ r' - r } \leq 1,\;\; \tau(H,r,G) \leq A,\;\; \tau(H',r',G)  < \infty 
\end{equation} 
one then has the bound
\begin{equation} \label{eqholdertau}
\abs{\tau(H',r',G) - \tau(H,r,G)} \leq \eta \left[ \vnorm{ r' - r }^{1/2} + \vnorm{ H'- H } \right].
\end{equation}
where the constant $\eta =\eta(A, a, \calP)$ is given by
\begin{equation} \label{eta}
\eta  =  3 g^{2+1/a} e^{A/a} F_g / F_1 + 20 F_g \log(F_g) /a.
\end{equation}
\end{Proposition}
%%%%%%%%%%%%%%%%%
\begin{proof}
Proof of this proposition is given in Appendix \ref{sec:ap3}.
\end{proof}

%%%%%%%%
We first state a key  large-deviation result for the one-step transition kernel controlling daily transitions.  The following theorem will be essential for developing the large-deviation formalism for random histogram trajectories in Section~\ref{LDPathTheory}.

%%%%%%%%%%%%%%%%%
\begin{Theorem} \label{transLD}
Fix any $0< a < 1 < d$ and the process  parameters $\calP$ (Definition~\ref{P}). One-step large deviations for the Markov chain $H_n$ are then controlled by two constants, $c\equiv c(d, a, \calP)$ and $N_0 \equiv N_0(d, a, \calP)$, as follows.  Consider any $N$-rational histograms $H, G \in \mcal{H}$ 
with $\min(b(H), b(G))>a$ and transition cost $C(H,G) \leq d$. Then, the one-step transition kernel $\mathcal{Q}(H,G)$ has a uniform large-deviation approximation valid for all $N >N_0$:  
\begin{equation} \label{onestepLD}
\frac{1}{N} \log \mathcal{Q}(H, G) = -  C(H,G) + o(1)
\end{equation}
with  $\abs{ o(1) } \leq c / \sqrt{N}.$
\end{Theorem}

\begin{proof}
 We sequentially construct constants $c_i>0, N_i >0$, which only depend on the fixed positive $(d, a)$ and the process parameters $\calP$. 

On day $n$, given $H_n= H $ and $R_n/N = r \in K_N(H)$, the population histogram $J_n$ after growth and mutations is a deterministic function $J_n = J(H, r)$ defined by \eqref{WWn}. Recall that $J(H,r)$ can be precisely estimated by $I = \Psi(H,r)$ (see Equation~\eqref{Jn}). 
To transition from  histograms $H_n=H$ to $H_{n+1}=G$, we define a  set  $A_N$ of ``acceptable'' mutation matrices $r_n=R_n/N$ by
\begin{equation} \label{A_N}
A_N = A_N(H,G) = \vset{ r \in K_N(H) \; | \; KL( G, J )< \infty \text{ and }  KL(G, I)< \infty}.
\end{equation}

\hl{Introduce the }approximate transition cost  $\apT(H,r,G) = mut(r, H) + KL( G, J(H,r) )$, which we expect to be very close to the explicit composite transition cost $\tau(H, r,G)$ 
(defined in Definition~\ref{tauHrG}).
We now give a technical but more tractable expression for the transition kernel.

Indeed, there  are  constants $c_0,N_0$ such that, for $N > N_0$ and for any $N$-rational triple $(H,r,G)$ verifying $H,G \in \mcal{H}$, $r \in K_N(H)$, $\min(b(H),b(G)) > a$, then the Markov transition kernel $\mathcal{Q}(H,G)$ (see Definition~\ref{kernel}) verifies 
\begin{equation} \label{trans4}
\mathcal{Q}(H, G)  =  \sum_{r \in  A_N} \exp \left(-N\,[\,\apT(H,r,G) + o(1) \,]\right) 
\end{equation}
with $\abs{ o(1) } \leq c_0  \log N / N.$

To prove \eqref{trans4}, apply \eqref{JJn} and the discussion in Section~\ref{JnMut} to $J=J(H,r)$ and $I=\Psi(H,r)$ to conclude that for $N > N_1 = 20 + 2 F_g / a$ and $c=  13 F_g / F_1$, one has 
\begin{equation} \label{compareJI}
\supp(J)= \supp(I) \;\text{~and~} \; \vnorm{J - I} \leq  c / N. 
\end{equation}
\hl{This forces} either $ \supp(G) \subset \supp(I) = \supp(J) $, so that $KL(G,I)$ and $KL(G,J)$ are both finite, or $\supp(G) \not\subset \supp(I) = \supp(J) $, so that $KL(G,I) = KL(G,J) =  \infty$. Define $\kappa = \kappa(H,r,G) = \min \vset{KL(G,J) , KL(G,I)}$. Set $c_1 = g^{1+1/a}$ and $c_2= c_1 \, c$. Combine \eqref{lipKLJI} with \eqref{compareJI} to obtain, for $r \in A_N$
\begin{equation} \label{**} 
\abs{KL( G, J ) - KL( G, I)}  \leq  c_1  e^{\kappa / a} \vnorm{ J - I } \leq c_2/ N 
\end{equation} 
which implies 
\begin{equation} \label{***} 
\abs{\apT(H,r,G) - \tau(H,r,G)}  \leq c_2/ N.
\end{equation}

\hl{The transition kernel} $\mathcal{Q}$ verifies by construction
\begin{equation} \label{trans0}
\mathcal{Q}(H, G)  =  \sum_{r \in K_N(H)} P( R_n/N = r \; | \; H_n = H )\,  P( H_{n+1} = G \; | \; H_n = H, \, R_n/N = r ).
\end{equation}

\hl{Since} $\vnorm{ r } \leq F_g,$ for all  matrices $r \in K_N(H)$, one must have $card(K_N(H)) \leq [ (N+1) F_g ]^{g^2}. $   From \eqref{LDselect1} we get, for $N > N_2 = N_1 + 2 / a + \log F_g$, 
\begin{equation} \label{trans2}
P(H_{n+1} = G \; | \; H_n = H, \, R_n/N = r) = \exp\left( - N  [ KL(G, J(H, r) )+ o_1(1) ]\right)
\end{equation}
with $\abs{o_1(1)} \leq (5 g+1) \log N / N$.  In \eqref{trans2}, the right-hand side is $0$ unless $r \in A_N$. Hence, the sum in \eqref{trans0} can be restricted to $r \in A_N$. Then, \eqref{LDmut1} yields, for $N > N_3 = N_2 + [ 4 + 2 \log(g) ] / a$,
\begin{equation} \label{trans3}
P( R_n / N = r \; | \;  H_n = H ) = \exp\left(- N [ mut(r, H) + o_2(1)]\right) 
\end{equation} 
with $\abs{o_2(1)} \leq 4 g^2  \log N / N$. For $r\in A_N$ and $N > N_3$, substitute \eqref{trans2} and \eqref{trans3} into \eqref{trans0} to finally prove \eqref{trans4} with $o(1) = o_1(1) + o_2(1),$ provided one sets $N_0=N_3$ and $c_0= 5 g + 1 + 4 g^2$.

For $H,G,r$ as above, and   $N >N_0$, we have now shown that 
\begin{align} 
\supp(J) = \supp(I) \;& \;  \text{~and~} \;  \vnorm{  J - I } \leq  c/ N \; \text{~for~} \; r \in K_N(H), \label{recall1}\\
\abs{KL( G, J ) - KL( G, I)} \;& \leq \; c_1 e^{\kappa / a} \vnorm{ J - I } \; \text{~for~} \; r \in A_N, \label{recall2}
\end{align} 
and we have proved Equation~\eqref{trans4}, namely,
$$
\mathcal{Q}(H, G)  =  \sum_{r  \in  A_N}  \exp (- N\, [\, \apT(H,r,G) + o(1) ]),
$$
with $\abs{o(1)} \leq c_0  \log N / N$.

Assume that  $C(H,G) < d$. Let $U_N = \exp(- N [ T(H,r,G) + o(1) ])$. Partition $A_N$ into two subsets, $\mathcal{M}$ and $\mathcal{R}$, defined by 
\begin{equation} \label{MM} 
\mathcal{M} = \vset{ r \in A_N \; | \; KL( G, J) > 4 d} \qquad \text{and} \qquad  \mathcal{R} = \vset{ r \in A_N \; | \; KL( G, J) \leq 4 d}
\end{equation} 
so that $\mathcal{Q}(H, G)  =  S(\mathcal{R})+ S(\mathcal{M})$, with $S(\mathcal{R})  =  \sum_{r  \in  \mathcal{R} }  U_N$ and $S(\mathcal{M})  =  \sum_{r  \in  \mathcal{M} } U_N.$ For $N > N_4 = N_0 + (c_0/ 2 d)^2$  and  $r \in \mathcal{M}$, the definitions of $U_N$ and $T(H,r,G)$ yield 
$$
- \frac{1}{N}\log(U_N) \geq KL(G, J) - c_0 \log N /N \geq 4 d - c_0 \log N / N \geq 2d.
$$ 
\hl{Since} $card(\mathcal{M}) \leq (N+1)^{g^2}$, this entails 
\begin{equation} \label{S0bound}
S(\mathcal{M}) \leq (N+1)^{g^2} \; e^{- 2 N d } \leq e^{- N d }
\end{equation} 
for $N > N_5 = N_4 + (2 g^2/ d)^2.$ By definition of $\kappa(H,r,G)$, when  $r \in \mathcal{R}$, one must have  $\kappa < KL( G, J ) \leq 4d$, and \eqref{recall2} yields, for $N > N_5$,
$$
\abs{KL( G, J ) - KL( G, I)}  \leq  e^{\kappa / a}(c_0^2 / N) \leq c_3 / N \
$$
with $c_3 = c_0^2  e^{4d/a}.$  This forces $-\frac{1}{N}\log{U_N} = \tau(H,r,G) + o_1(1)$, with $\abs{o_1(1) } \leq \abs{o(1) }+ c_3 / N \leq c_4  \log(N) / N$, with $c_4 = c_0 + c_3$.  By definition of $S(\mcal{R}),$ this implies, for $N > N_5$, 
\begin{equation}\label{S3}
S(\mathcal{R}) = \sum_{r  \in  \mcal{R} }  \exp (- N  [ \tau(H,r,G) + o_1(1) ] ).
\end{equation}

\hl{By definition} of $C(H,G)$, one has $\tau(H,r,G) \geq C(H,G)$, so that \eqref{S3} yields
\begin{equation} \label{SS+}
S(\mathcal{R})  \leq  (N+1)^{g^2} \exp (- N \,C(H, G) + c_4 \log(N)) \leq \exp (- N \, C(H, G) + c_5 \log N )
\end{equation} 
for $N > N_5$ and  $c_5 = c_4 +g^2 +1$. Equation~\eqref{S0bound} then yields $ S(\mathcal{M}) \leq e^{- N d} \leq e^{- N C(H,G)} $, since $C(H,G) \leq d$. Combining this with \eqref{SS+} we obtain, for $N > N_5$ and $c_6 = 3 c_5$, 
$$
\mathcal{Q}(H, G) = S(\mathcal{M}) + S(\mathcal{R}) \leq ( 1 + N^c_5 ) e^{- N \,C(H,G)}  \leq \exp ( - N \,C(H,G) + c_6 \log{N}).
$$
 \hl{This yields} the large-deviation upper bound: 
\begin{equation} \label{upLD}
\frac{1}{N} \log \mathcal{Q}(H, G) \leq - C(H,G) + c_6 \log N /N. 
\end{equation}

\hl{By definition} of $C(H,G)$, there exists a matrix $w= w(H,G,N) \in K(H)$ such that 
\begin{equation}\label{wHGN}
C(H,G) \leq \tau(H,w,G) \leq C(H,G) + 1/ \sqrt{N}.
\end{equation}

\hl{Lemma} \ref{density} shows that for $N > N_6 =N_5 + \frac{ g^2}{a F_1}$, there is an $N$-rational matrix $s \in K_N(H)$ such that $\supp(s) = \supp(w)$ and $\vnorm{ s - w } \leq g/N.$   Equation~\eqref{wHGN} implies $\tau(H,w,G) \leq d + 1$ so that $\tau(H,s,G)$ must be finite since $ \supp(s) = \supp(w)$. A fortiori, $KL(G, \Psi(H,s))$ is finite, implying $\supp(\Psi(H,s)) \subset \supp(G)$. However, since $s \in K_N(H)$, Equation~\eqref{recall1} forces $\supp(J(H,s)) = \supp(\Psi(H,s))$, so that $KL(G, J(H,s))$ must also be finite.  Therefore, one has $s \in A_N$.  The H\"{o}lder continuity in $r$ for the composite cost function $\tau(H,r,G)$ was established above by  \eqref{eqholdertau}; we now apply it to get $c_7 = \eta(d+1, a, \calP)$ such that
\begin{equation}\label{tau.sw}
\abs{\tau(H,s,G) - \tau(H,w,G)} \leq c_7 \vnorm{ s-w }^{1/2} \leq c_7 \sqrt{g} / \sqrt{N}.
\end{equation}

\hl{Set} $c_8 = 1 + c_7 \sqrt{g}$ and $N_7 = N_6 + c_8^2$. Using \eqref{wHGN}, we get, for $N>N_7$,
\begin{equation} \label{tauHsG}
\tau(H,s,G) \leq \tau(H,w,G) + c_7 \sqrt{g} ) / \sqrt{N} \leq C(H,G) + c_8 / \sqrt{N} \leq d+1.
\end{equation}
\hl{Then}, $\kappa(H, s, G) \leq KL(G,\Psi(H,s)) \leq \tau(H,s,G)) \leq d + 1$.  For $N> N_7$, since $s \in A_N$, apply \eqref{recall2} to the triple $(H,s,G)$ to get 
\begin{equation} \label{***} 
\abs{KL( G, J(H,s) ) - KL( G, \Psi(H,s) ) }  \leq  c_8 /N   
\end{equation}
with $c_8 = c_0^2 e^{(d+1)/a}.$ Since $d \geq 1$, this implies
$$
KL( G, J(H,s) ) \leq KL( G, \Psi(H,s) ) + c_8 /N \leq d + 1 + c_8 /N < d + 2 < 4 d 
$$ 
for $N > N_8 = N_7 + 1/c_8.$ Hence, for $N > N_8$, the matrix $s \in A_N$ must belong to $\mathcal{R}$, and \eqref{S3} provides the lower bound :
\begin{equation}\label{S4}
S(\mathcal{R}) = \sum_{r \in \mathcal{R} } \exp (  - N  [ \tau(H,r,G) + o_1(N) ]  ) \geq \exp ( - N  [ \tau(H,s,G) + o_1(1) ] ). 
\end{equation}

\hl{Due to the bound} on $o_1(1)$ and \eqref{tauHsG}, 
$$
\frac{1}{N} \log S(\mathcal{R}) \geq - \tau(H,s,G) +o_1(1) \geq - C(H,G) - c_8 (1/ \sqrt{N} + \log N /N )
$$
for $N > N_8.$  Finally, since $\mathcal{Q}(H, G) \geq S(\mathcal{R})$, setting $c_9 = 2 c_8$ yields
$$
\frac{1}{N} \log \mathcal{Q}(H, G) \geq - C(H,G) - c_9 / \sqrt{N}.
$$
\hl{In view } of \eqref{upLD}, we now obtain, for $N > N_8,$
\begin{equation} \label{trans}
- C(H,G) - c_9 / \sqrt{N} \leq  \frac{1}{N} \log \mathcal{Q}(H, G)  \leq - C(H,G) + c_6 / \sqrt{N}.
\end{equation}

\hl{This conclude}s the proof of Theorem~\ref{transLD}.
\end{proof}
%%%%%%%%%%%%%%%%

%
The  H\"{o}lder continuity property of the composite cost $\tau(H,r,G)$ proved above can be used  to derive uniform holder continuity for the one-step cost function $C(H,G)$, as outlined in the next theorem.

%%%%%%%%%%%%%%%%%
\begin{Theorem} \label{holderCost}
Fix $0< a < 1, d >0,$ and the parameter set $\calP$ (Definition~\ref{P}). Consider any histograms $H, G, H',G'$ in the compact set $\mcal{H}(a)$ defined by  \eqref{HIST(a)}, and verifying 
\begin{equation} \label{conditions}
    \supp(G') = \supp(G), \;\; \supp(H) = \supp(H'), \;\;
    \vnorm{ H' - H } \leq a/F_g \;, \;\;  C(H,G) \leq d.
\end{equation}
\hl{Then, there} exists a constant $c \equiv c(d, a, \calP)$ such that 
\begin{equation} \label{continuityCHG}
\abs{ C(H',G') - C(H,G) } \leq c  \left(\vnorm{ H' - H }^{1/2} + \vnorm{ G' - G }\right). 
\end{equation}
\end{Theorem}
%%%%%%%%%%%%%%%%%
\begin{proof}
This proof is quite technical and is displayed in Appendix \ref{sec:ap3}.
\end{proof}

%%%%%%%%%%%%%%%%%%%%%%%%%%%%%%%%%%%%%%%%%%%%%%%%%%%%%%%%%
\subsection{Computation of the One-Step Transition Cost}
\label{subs3.5}
One of the main goals of this paper was to develop a computationally implementable set of rigorous large-deviation formulas in order to compute numerically the most likely evolutionary paths connecting any known initial histogram $H$ to any histogram $G$ observed at some later time $T$. In fact, for many experimental observations of $\{H,G\}$ the time $T$ may not even be known precisely. To this end, we will provide an explicitly computable expression for the one-step cost $C(H,G)$ defined by the theoretical formula \eqref{costHG}. This formula requires minimizing the explicit convex function $\tau(H,r,G)$ over all $g \times g$ matrices $r \in K(H)$. We discovered that for small mutation rates $m \leq 10^{-6}$ and large population sizes $N\geq 10^6$, the minimization of $\tau(H,r,G)$ can be achieved explicitly using an asymptotic expansion, and this minimization procedure provides a very accurate approximation of $C(H,G)$ with error terms of order $m^2 < 10^{-12}$.

The explicit computation of  precise approximations for $C(H,G)$ is here deliberately restricted to the simpler case when both $H$ and $G$ are interior histograms, i.e., verify $H(j)>0$ and $G(j) >0$ for all $j$. This enables a first-order Taylor expansion of $C(H,G)$ in the mutation rate $m$. For more generic $H,G$ having some coordinates equal to 0, good approximations of $C(H,G)$ are  feasible but are more complicated because Taylor formulas in $m$ are no longer valid (see remarks below). 

\revision{Recall that $F_j$ are growth factors for genotype $j$. Therefore, $\Phi_j = F_j H(j) / \langle F,H\rangle $ approximates the histogram after growth, but before the mutation phase.}

 \begin{Definition}
\label{interiohist}
A histogram $H \in \mathcal{H}$ is called {interior} if $H(j) > 0$ for all $1 \leq j \leq g$. We denote the open set of interior histograms as $\interior{\mathcal{H}}$. 
\end{Definition}
We will prove below that $C(H,G)$ is  smooth in  $H,G$ whenever 
$H,G \in \interior{\mathcal{H}}$.
Recall (see Equation~\eqref{Kj}) that matrices $r$ in the convex set $K(H)$ satisfy
\begin{align} 
r_{j,k} & \geq 0  \text{ for all } (j,k), \label{cons0} \\
\sum_k  r_{j,k} & <  F_j H(j) \text{ for all } j,  \label{cons1} \\
r_{j,k} & =  0   \text{ when } q_{j,k} = 0. \label{cons2}
\end{align}

\hl{The interior} $\interior{K}(H)$ of $K(H)$ is the set of all $r$ verifying \eqref{cons1},  \eqref{cons2}, and  $r_{j,k}> 0$ whenever $q_{j,k} >0$.  The following theorem gives an explicit computation of the one-step cost function $C(H,G)$.

\revision{Restrictions \eqref{cons0}--\eqref{cons2} arise from the constraint that the number of mutants $j \to k$ cannot be larger than the number of cells of the genotype $j$.}
%%%%%%%%%%%%%%%%%%
\begin{Theorem}[Explicit first-order expansion of the one-step cost C(H, G)]
\label{thCHG}
Fix the vector  $F=\{F_1,\dots,F_g\}$ of growth factors,  the transition matrix $Q= (q_{j,k})$ for mutations, and any $0 < a < 1$. Let $\Gamma(a) \subset \interior{\mathcal{H}}$ be the set of all histograms $H$ such that $H(j) > a$ for all $j= 1,\dots, g$. 
There is then a constant $c \equiv c(a,F,Q) > 0$ such that for all $H, G \in \Gamma(a)$, and all mutation rates $0 \leq m < c$, the transition cost $C(H,G)$ is a finite $C^{\infty}$ function of $(m, H, G)$. Moreover, $C(H,G)$ has an explicit first-order expansion in $m$ given by
\begin{equation}\label{formulaCHG}
C(H, G) =  KL(G,\Phi(H)) + m \sum_{j,k}  F_j H(j) q_{j,k} [1 - U_k/U_j ] + O(m^2),
\end{equation}
where  
\begingroup\makeatletter\def\f@size{9}\check@mathfonts
\def\maketag@@@#1{\hbox{\m@th\normalsize\normalfont#1}}
\[
\Phi_j(H) = F_j H(j) / \langle F,H\rangle \, , \quad KL(G,\Phi(H)) = \sum_j  G(j) \log(G(j) / \Phi_j(H) \, , \quad
U_j =  \exp \left( \frac{G(j)}{F_j H(j)}\right).
\]
\endgroup%

\revision{\hl{The cost fun}ction $C(H,G)$ in \eqref{formulaCHG} consists of two parts. The Kullback--Leibler divergence represents the cost associated with the multinomial selection step. It measures the ``distance'' between the population histogram after growth and the target histogram $G$. The second term represents the cost associated with mutations. It quantifies how unlikely it is that mutations occurring during the growth phase transform the post-growth population histogram into the target histogram $G$.}

\end{Theorem}
%%%%%%%%%%%%%%%%%%
\begin{proof}
Denote $f_{j,k} = q_{j,k} F_j H(j)$. For $r \in K(H)$ recall that $\tau(H,r,G) = mut(H,r) + KL(G, \Psi(H,r))$, with 
\begin{align} 
\Psi_i = \Psi_i(H,r) & =  \frac{1}{\langle F, H\rangle} \left( F_i H(i) - \sum_k  r_{i,k} + \sum_k r_{k,i} \right),  \label{eqpsi}\\
KL(G, \Psi(H,r)) & = \sum_i  G(i) \log ( G(i) / \Psi_i ),  \label{eqKL}\\
mut(H,r) & = \sum_{(j,k) \in \supp(Q)} \left[ mf_{j,k} + r_{j,k} \log r_{j,k} - r_{j,k} - r_{j,k} \log (m f_{j,k})\right],  \label{eqmut}
\end{align}

\hl{To minimize} $\tau(H,r,G)$, we only need to consider $r \in K(H)$ with $\tau(H,r,G)$ finite, which holds iff $\Psi_i(H,r) > 0$ for all $i$. 
Due to convexity in $r$ of $\tau(H,r,G)$ we first seek to minimize $\tau(H,r,G)$ over $r \in \interior{K}(H)$. Such a minimizing  $r$ must verify 
\begin{equation} \label{sys}
\partial_ {r_{j,k}} \tau(H,r,G) = 0 \qquad \,\text{whenever}\;\;q_{j,k} >0.
\end{equation}

\hl{The set} $\{ (j,k) \;|\; q(j,k) >0$ is by definition the  support $\supp(Q)$ of the fixed transition matrix $Q = (q_{j,k})$. For $(j,k) \in \supp(Q),$ basic derivations and algebra reduce \eqref{sys} to the system 
\begin{equation}\label{Dtau}
0 =  \log\frac{r_{j,k}}{m f_{j,k}} - \sum_i \; \frac{G(i)}{\Psi_i} \partial_{r_{j,k}} \Psi_i  
\end{equation}

\hl{The coefficients} $p(i,j,k) = \partial_{r_{j,k}} \Psi_i $ are constants given by
\begin{equation}\label{DPsi}
\langle F,H \rangle  p(i,j,k) = 
\begin{cases} 
-1, &  i =j  \neq k,\\
1, &   j \neq k = i, \\
0 &  \text{otherwise.}
\end{cases}
\end{equation}

\hl{Let} $d = card(\supp(Q)).$ Define the vector $ x \in (\mbb{R}^+)^d $ by $x_{j,k} = r_{j,k}/(m f_{j,k})$ for $q_{j,k} >0$.
Note the double-index notation for the coordinates $x_{j,k}$ of the vector $x$. 
Substitute \eqref{DPsi} into \eqref{Dtau} to obtain
\begin{equation}\label{xsolv}
\log x_{j,k} = - G(j)/A_j + G(k)/A_k,
\end{equation}
where $A_i \equiv A_i(x)$ is given by 
\begin{equation*}
A_i(x)  = F_i H(i) + \sum_k  m \left[ -  f_{i,k}\, x_{i,k} +  f_{k,i}\,x_{k,i} \right]. 
\end{equation*}

\hl{To enforce}  \eqref{cons1}, the vector $x$  must verify, for all $1 \leq i \leq g,$ the set of strict linear constraints 
\begin{equation}\label{consx}
F_i H(i) - m \sum_k f_{i,k} x_{i,k} > 0.
\end{equation}

\hl{Call} $\mathcal{O}$ the open set of all $x \in (\mbb{R}^+)^s$ verifying \eqref{consx}. Rewrite \eqref{xsolv} as an implicit equation $W(m,x) = 0$ for $x \in \mathcal{O}$, where the function $W(m,x)$ is given by
$$
W_{j,k}(m,x) = x_{j,k} - \exp\left( - G(j)/A_j + G(k)/A_k \right) 
$$
for all $(j,k) \in \supp(Q)$, $x\in \mcal{O}, $ and $m \geq 0.$ 
We can treat $W(m,x)$ as a vector-valued function $W(m,x): (\mbb{R}^+)^d \to \mbb{R}^d$.
Then, $W(m,x)$ is of class $C^{\infty}$ in $(m,x)$, as well as in $H,G \in \mathcal{H}^\circ$.
For $m=0$ and $(j,k) \in \supp(Q)$, the system $W(m,x) = 0$ has a unique solution $\wt{x} \in \interior{K}$ given by
\begin{equation} \label{xsolution}
\wt{x}_{j,k} = \exp\left(- G(j)/[ F_j H(j) ] + G(k)/ [ F_k H(k) ]\right).
\end{equation}

\hl{For all} $(s,t) \in \supp(Q)$ and all $i$, the derivatives $a(i,s,t) = \partial_{x_{s,t}} A_i(m,x)$ verify $a(i,s,t) =m b(i,s,t)$, where the only non-zero terms of $b(i,s,t)$ are $b(i,i,t) = -f_{i,t}$ and  $b(i,s,i) = f_{s,i}$.  The Jacobian $JAC(m,x)$ of $W(m,x)$ has coefficients $\partial_{x_{s,t}} W_{j,k}$ given for all $(j,k)$ and $(s,t)$ in $\supp(Q)$ by 
$$
 \partial_{x_{s,t}} W_{j,k} = \delta_{(j,k), (s,t)} - m Y_{(j,k),(s,t)} \exp\left(- G(j)/A_j + G(k)/A_k\right) 
$$
where $Y_{(j,k),(s,t)} = b(j,s,t) G(j) / A_j^2 - b(k,s,t) G(k) /A_k^2,$ 
and $\delta_{(j,k), (s,t)} \equiv \delta_{j,s} \delta_{k,t}$ is the Kronecker delta function.  
At the point $(0,\wt{x}),$ the Jacobian  $JAC(0,\wt{x})$ is hence equal to the identity matrix and thus invertible. The classical implicit function theorem then applies to $W(m,x) = 0$ and provides $c=c(a,F,Q) > 0$ such that the equation $W(m,x) = 0$ has a unique solution $x^*(m) \in \mathcal{O}$ for $m \leq c$  and $H,G \in \Gamma(a)$. The same classical theorem implies that $x^*(m) = x(m,H,G)$ is of class  $C^{\infty}$ in $(m,H,G) $.
Define $r^*(m)$ by
\begin{equation*}
r^*_{j,k}(m) =
\begin{cases}
m\, x^*_{j,k}(m),&  \;\text{when}\;\; q_{j,k} > 0,\\
0,&  \text{otherwise.}
\end{cases}
\end{equation*}
\hl{Then}, $r^*(m)$ is a solution of \eqref{sys} and inherits from $x^*(m)$ the $C^{\infty}$ smoothness in $(m,H,G)$. Moreover, $r^*(m)$ verifies the constraints \eqref{cons0}--\eqref{cons2} since $x^*(m) \in \mathcal{O}$. The positivity of all coordinates of $x^*(m)$ implies $r^*_{j,k}(m) >0$ for $q_{j,k} >0$. 
Hence, $r^*(m) \in \interior{K}(H)$ for $m< c$ and solves \eqref{sys}. The strict convexity of $\tau(H,r,G)$ on the open convex set $\interior{K}(H)$ forces $r^*(m)$ to be the unique minimizer of $\tau(H,r,G)$ for $r \in K(H)$ so that $C(H,G) = \tau(H,r(m),G)$.  The function $(H,r,G) \to \tau(H,r,G)$ is $C^{\infty}$ for $H,G \in \mathcal{H}^\circ$ and $r \in \interior{K}(H)$. Hence, the function $(m, H,G) \to C(H,G)$ is also $C^{\infty}$ for $m < c$ and $H,G \in \Gamma(a)$.

To obtain the first-order expansion for the solution $r^*(m)$
and the cost function $C(H,G)$, define 
\begin{equation} \label{UjEjk}
U_j  =  \exp \left( G(j)/ [ F_j H(j) ]\right) \quad \text{and} \quad  E_{j,k} =  U_k / U_j
\text{  for all  } 1\le j,k \le g, 
\end{equation}
so that $\wt{x}_{j,k} = f_{j,k} E_{j,k} = F_j H(j) q_{j,k} U_k/U_j$. To simplify notation, write $u \simeq v$ whenever $u(m) = v(m) + O(m^2)$ as $m \to 0$.  As $m \to 0,$ the differentiability of $x^*(m)$ gives   
\begin{equation}\label{rsolv}
r^*_{j,k}(m) =  m\, x^*_{j,k}(m) \simeq m\, \wt{x}_{j,k} =  m\,  f_{j,k} E_{j,k}.
\end{equation}
\hl{Inserting} \eqref{rsolv} into \eqref{eqmut} yields $mut(H,r^*(m)) \simeq  m\, \mu$, with 
$
\mu = \sum_{j,k}  f_{j,k} \left(1 - E_{j,k} + E_{j,k}\log E_{j,k} \right). 
$
Since $\log(E_{j,k}) = - G(j)/ [ F_jH(j) ] + G(k)/ [ F_k H(k)]$, we have 
\begin{equation} \label{eqmu}
\mu = \sum_{j,k} f_{j,k}( 1 - E_{j,k}) + \sum_{j,k}  q_{j,k} E_{j,k}\left(- G(j)  + G(k) \frac{F_j H(j)}{F_k H(k)}\right). 
\end{equation}

\hl{Substituting} \eqref{rsolv} into \eqref{eqpsi}, we obtain
\begin{align*}
\Psi_j &=  \Psi_j(H,r^*(m)) \simeq  (F_j H(j) - m t_j) / \langle F,H\rangle; \\
t_j &= \sum_k  F_j H(j) q_{j,k} E_{j,k} -  \sum_k  F_k H(k) q_{k,j}E_{k,j}.
\end{align*}
\hl{This yields} $\log(\Psi_j) \simeq  \log\left(F_j H(j) /\langle F,H\rangle\right) -  m \frac{t_j}{F_j H(j)}.$  Recall that $\Phi_j(H) = F_j H(j)/\langle F,H\rangle$ is the population histogram at the end of the daily deterministic growth starting with histogram $H$.  Substituting the expansion of $\log(\Psi_j)$ into \eqref{KL}, we obtain
\begin{align*}
&KL(G, \Psi(H, r^*(m)) \simeq  \kappa + m \eta;\\
\kappa = KL(G, \Phi) =  \sum_j  G(j) &\log\left(\frac{G(j)}{\Phi_j}\right)  > 0; \qquad  \qquad \eta  = \sum_j  \frac{ G(j) t_j}{F_j H(j)}.
\end{align*} 
\hl{This implies} $\eta  = \sum_{j,k}  \left[G(j) q_{j,k} E_{j,k} - \frac{G(j)}{F_jH(j)} F_k H(k) q_{k,j}\right] E_{k,j}$.  Exchange $j$ and $k$ in the second term  of the previous sum to get 
\begin{equation}
\label{eta1}
\eta = \sum_{j,k} q_{j,k} E_{j,k}\left[G(j) - G(k) \frac{F_j H(j)}{F_k H(k)}\right].
\end{equation}

\hl{Hence}, $C(H,G) = \tau(H,r^*(m),G) = \kappa + m\, (\eta +\mu) + O(m^2)$ is the first-order expansion of $C(H,G)$.  Combine \eqref{eqmu} and \eqref{eta1} to obtain
$
\eta +\mu = \sum_{j,k} F_j H(j) q_{j,k} \left( 1 - U_k/U_j \right),
$
which concludes the proof of \eqref{formulaCHG} and of the theorem.
\end{proof}
%%%%%%%%%%%%%%%%%%
%%%%%%%%%%%%%%%%%%

\revision{
\subsection{One-Step Transition Cost $C(H,G)$  for Boundary Histograms $H,G$} 
\label{sec:boundary}
In this section we provide some preliminary results on 
extending Theorem~\ref{thCHG} to pairs $(H,G)$ in which some coordinates are equal to zero. This is definitely a more complicated situation and we will address it in more detail in a future paper. In this section, we briefly outline how to handle boundary histograms in one key example,  where the initial and target histograms $H$  and $G$ verify  
\begin{equation}
H(g) = 0 \; ; \quad   H(j) >0 \; \text{~for all~} \; j \neq g; \quad  G(k)  >0 \; \text{~for all~} \; k.
\end{equation}
\hl{Situations } where $H$ has several coordinates equal to zero can be handled similarly. To simplify formulas, assume that mutants' transition probabilities satisfy $q_{j,k}>0$ for all $j\neq k$.

Recall that the unknown  rates $r = (r_{j,k})$ must verify, for all $j,k$,
$$ 
r_{j,k} \ge 0, \quad r_{j,j} = 0, \quad \sum_k r_{j,k} < F_jH(j)
$$
and \hl{denote}  $f_{j,k}= F_j H(j) q_{j,k}$, so that   $f_{g,k} = 0$ for all $k$. 
  
The composite one-step cost associated with $r$ for fixed histograms $H,G$ is 
$$
\tau(r) = mut(H, r) + KL(G, \Psi(r))
$$
with the previously defined  (see Theorem~\ref{thCHG})
\begin{align*}
mut(r) &= \sum_{ j | H(j)  >0}\sum_k \left[ m f_{j,k} + r_{j,k} \log(r_{j,k}) - r_{j,k} - r_{j,k} \log(m f_{j,k} ) \right], \\
\Psi_j(r) &= \frac{1}{\langle F, H\rangle} \left( F_j H(j) + \sum_k (r_{k,j} - r_{j,k}) \right), \\
KL(G, \Psi(r)) &=  \sum_ j G(j) \log\left( \frac{G(j)}{\Psi(r)} \right). 
\end{align*}

\hl{Since} $H(g) =0$, one has
$$ 
\Psi_g(r) =\frac{1}{\langle F, H\rangle} \sum_{j\neq g} r_{j,g}    \quad \text{and} \quad r_{g,k} = 0 \;\; \text{~for all~} k.
$$ 
\hl{The one-step transition} cost $C(H,G)$ is the minimum of $\tau(r)$ over all $r = (r_{j,k})$, as above.
In addition, the one-step cost
$\tau(r)$ is infinite unless $r$ verifies the restriction $ \sum_k  r_{k,g}  >0 $, which we now impose since we want to minimize $\tau(r)$; this restriction simply states that there are cells that can mutate into the genotype $g$ in one step.

Define the set of variables $x_{j,k} \ge 0$ as     
$$
x_{g,k}=0, \quad 
x_{j,j}= 0 \; \text{~for all~} \; j, \quad 
x_{j,k}= \frac{r_{j,k}}{m f_{j,k}} \; \text{~for all~} \; j \neq g. 
$$

\hl{To derive} the one-step cost, we consider fixed histograms
$H$ and $G$, and the mutation rate $m <10^{-6}$ is very small.
When  $r= r(m,H,G)$ minimizes $\tau(r)$,  the partial derivatives of $\tau$ with respect to each  positive $r_{j,k}$ must be zero, and this yields, as in Theorem~\ref{thCHG} above,
\begin{equation} \label{eq1}
\log x_{j,k} = - \frac{G(j)}{A_j(x)} + \frac{G(k)}{A_k(x)}  \; \text{~when~}  \; x_{j,k} >0,
\end{equation}
where
$$ 
A_j(x) = F_j H(j) + m  \sum_k  (f_{k,j} \, x_{k,j} - f_{j,k} \, x_{j,k}).
$$

\hl{We need} to study  the limits 
$$
y_{j,k} = \lim_{m \to 0} x_{j,k}(m,H,G).
$$ 

For $j \neq g$ and $k \neq g$, i.e., for  $H(j) >0$ and $H(k) >0$, we obtain, as in Theorem~\ref{thCHG},
$$
\log y_{j,k} = - \frac{G(j)}{F_j H(j)} + \frac{G(k)}{F_k H(k)}  
$$ 
and then 
$$
y_{j,k} = U_j/U_k  \text{~with~} U_j= \exp\left(- \frac{G(j)}{F_j H(j)} \right) 
$$ 
which yields 
$$
r_{j,k}= m f_{j,k} x_{j,k} \sim  m f_{j,k} \frac{U_j}{U_k} \text{~for~} j \neq G \text{~and~} k\neq j. 
$$

For $j=g$ and any $k$ we have $y_{g,k}=0$ since $x_{j,k}=0$, as well as $r_{j,k}=0$.
For  $j\neq g$, Equation~\eqref{eq1} becomes 
$$
\log(x_{j,g}) = - \frac{G(j)}{F_jH(j)} + \frac{G(g)}{\sum_{k \neq g}  (f_{k,g} \, x_{k,g}) }  \;\;\text {for}\; j \neq g
$$
since  $H(g)=0$ and $A_g(x) = m \sum_{k \neq g}  (f_{k,g} \, x_{k,g} ) $.
The variables $X_j = x_{j,g}/U_j$ thus  verify for all $j \neq g$ that
$$
\log(X_j)= \frac{G(g)}{m \sum_{k \neq g}  (U_k \, f_{k,g} \, X_k )}.
$$
\hl{The right-hand side} does not depend on $j \neq g$ and hence we can write $X_j= z  >0$ for all $j \neq g$. The preceding equation  shows that $z$ is a solution of 
\begin{equation}
\label{eqz}
z \log(z) = \frac{G(g)}{m \, u} \quad \text{with} \quad u = \sum_{k \neq g}  (U_k f_{k,g} ).
\end{equation}

\hl{Note} that $u = u(H,G) >0$ does not depend on $m$. This implies that $z >1$. 
The solution of Equation \eqref{eqz} can be expressed as
$z = L \, G(g) \, (m \, u)^{-1}$, where for any $t>0$
$$ 
L(t) \ge 1 \;\; \text{denotes the unique solution of} \;\; L(t)\log(L(t)) = t.
$$ 
\hl{Thus},  we obtain
$$
x_{j,g}= U_j L\left(\frac{G(g)}{m \, u}\right) \;\;\text{for all} \;\; j \neq g.
$$
\hl{Hence}, for fixed histograms $H$ and $G$ with $G(g)>0$ and $H(g)=0$, we obtain $y_{j,g}= \lim_{m\to 0} x_{j,g}= +\infty$, since for large $t$ one has $L(t) \approx {t}/{\log(t)}$. More precisely,  for all $j\neq g$ we obtain
$$
r_{j,g}= m f_{j,g} x_{j,g} = U_j f_{j,g} m L\left(\frac{G(g)}{m \, u}\right) \sim  \frac{f_{j,g} G(g)}{u} \frac{1}{|\log(m)|} .
$$

\hl{Substituting} into $mut(r)$ and $\Psi(r)$ the expression above for $r=(r_{j,k})$, which minimizes $\tau(r)$, we  obtain a fully  explicit expression of $C(H,G)$ for fixed $m, H, G$.
However, the Taylor expansion of $C(H,G)$  as   $m \to 0$ obtained in Theorem~\ref{thCHG}
is no longer valid, since  when  $m\to 0$, the transition cost  $C(H,G)$ tends to $+\infty$ at the very slow rate 
$$
C(H,G) \sim  C \log(|\log(m)|) 
$$
where $C>0$ is a constant that can be computed explicitly.

The formulas just outlined for boundary histograms are numerically applicable to Lenski-type experiments when $10^{-8} \le m \le 10^{-6}$. However,  for even smaller  mutation rates, $m \ll 10^{-9}$,   the only ``observable'' one-step transitions from $H$ to $G$ when $H(g) =0$ and $G(g)>0$ actually require having $G(g) = b m$ with some constant $b>0$. When this is the case, the formulas given above for $r_{j,k}$ still provide an explicit leading-order expansion of the one-step transition cost $C(H,G)$. }

%%%%%%%%%%%%%%%%%%%%%%%%%%%%%%%%%%%%
\section{Large Deviations for Evolutionary Trajectories}
\label{LDPathTheory}
For fixed large population size $N$, and fixed small mutation rate $m$, the stochastic genetic evolution of the daily populations $pop_n$ across discrete time $n$ is modeled  here by the random histogram path $\{ H_1\,\cdots  H_T \}$, where $H_n$ is the histogram of $pop_n$. The discrete Markov chain $\{H_n, \, n=1,2,\dots\}$ belongs to the compact convex state space $\mcal{H}$, and has Markov transition kernel $\mcal{Q}_N$. 
The large-deviation analysis of daily transitions discussed in Section~\ref{LDcycles}, provides an explicit formula for the one-step rate function $C(H,G) \ge 0$, and for large $N$, the precise  large-deviation approximation  
\[
\frac1N \log{Q_N(H,G)} \approx - \, C(H,G)
\]
that roughly says that $Q_N(H,G) \approx  \exp(- N\, C(H,G))$.
In this section, we rigorously extend this one-step large-deviation approximation to the ``path space'' of random histogram trajectories for the Markov chain describing the time-evolution of histograms.

\subsection{Large Deviations for a Single Trajectory}
\label{sec:LDTtraj}

\vspace{6pt}
\begin{Definition}
\label{def:OmegaT}
For any fixed \emph{time horizon} $T \geq 2$, denote by $\Omega_T \equiv \mathcal{H}^T $ the \emph{path space} of all possible histogram trajectories $\mbold{H} = [ H_1\,\cdots H_T ]$ with all $H_n \in \mathcal{H}.$
\end{Definition}
\begin{Definition}
The \emph{distance} between trajectories $\mbold{H}$ and $\mbold{H}'$ is defined  by \\
$\vnorm{ \mbold{H} - \mbold{H}' } = \max_{n =1\cdots T}  \vnorm{ H_n - H'_n }$. 
\end{Definition}
Intuitively a ball $B(\mbold{H}) \in \Omega_T$ of small radius can then be viewed as a ``thin tube'' of trajectories around the tube axis $\mbold{H}$.
Since the random population histograms $H_n$ are $N$-rational, our Markov chain actual trajectories will always belong to the set of \emph{N-rational}   trajectories $[H_1 \, H_2\cdots H_T]$, such that $H_n \in \mcal{H}_N$ for $1\leq n\leq T$.

Next, we construct an explicit rate functional $\lambda(\mbold{H}) \ge 0$ defined for all $\mbold{H} \in \Omega_T$, which extends the one-step rate function $C(H,G)$ to a $T$-step rate function defined on the path space. Furthermore, we will prove that $\exp(-N\, \lambda(\mbold{H}))$ roughly approximates the probability that the trajectory of the underlying Markov chain remains within any given very thin tube of trajectories centered around $\mbold{H}$.
%%%%%%%%%%%%%%%%%%
\begin{Definition}
\label{def:lambda}
For any trajectory $\mbold{H} \in \Omega_T$, we define the \emph{large-deviation rate function} $\lambda : \Omega_T \to [ 0, \infty ]$ by
\begin{equation} \label{lambda}
\lambda(\mbold{H}) = \sum_{n=1}^{T-1} C(H_n, H_{n+1}).
\end{equation}

\hl{We then} define the \emph{large-deviation set functional} $\Lambda(\Gamma) \in [0, \infty]$ for any $\Gamma \subset \Omega_T$ by 
\begin{equation} \label{Lambda}
\Lambda(\Gamma) = \inf_{ \mbold{H} \in \Gamma} \; \lambda(\mbold{H}).
\end{equation}
\end{Definition}
%%%%%%%%

 \begin{Definition}
 \label{def:bboldH}
 For any path $\mbold{H} \in \Omega_T$, define its \emph{essential minimum}  $b(\mbold{H}) = \min_{n =1\cdots T}  b(H_n)$,
 where $b(H_n)$ is the essential minimum of $H_n$ defined in Definition~\ref{def:suppmin}.
 \end{Definition}

\revision{Note that Definition~\ref{def:bboldH} restricts trajectories to consist only of interior histograms, in which all genotypes are present in the population. However, some genotype frequencies may be very small, for example, $N^{-1}$.}

The next theorem justifies calling $\lambda$ the rate function of the path space.  
 This is a uniform large-deviation result for single $N$-rational trajectories, and will be a key point in obtaining large-deviation estimates for subsets of the path space (see Theorem~\ref{thsetLD}).

%%%%%%%%%%%%%%%%%%
\begin{Theorem} \label{thtrajLD}
Fix the path length $T \geq 2$, two positive constants  $a,d >0$, and parameters $\calP$. 
Denote by $\mbold{H} =[H_1,\dots, H_T] \in \Omega_T$ the random path of the Markov chain describing daily transitions. Then, $(d,a,\calP)$ determine positive constants $c$ and $N_0$ such that the following holds. For any given fixed $N$-rational path $\mbold{h}= [h_1 \, h_2\, \dots \,h_T] \in \Omega_T$ such that $\lambda(\mbold{h}) \leq d$ and $b(\mbold{h}) \geq a$, one has 
\begin{equation} \label{trajsingle}
\frac{1}{N} \log P( \mbold{H} = \mbold{h} \; | \; H_1 = h_1) = - \lambda(\mbold{h}) + o(1) \;\text{~for all~}\; N > N_0
\end{equation} 
 with $\abs{ o(1) } \leq c\, T / \sqrt{N}.$
\end{Theorem}
\revision{Theorem~\ref{thtrajLD} shows that the large-deviation rate function can be computed as a sum of cost functions associated with the daily step; see \eqref{lambda}.}
\begin{proof}
For all $n \leq T-1$, we must have $C(h_n, h_{n+1})\leq d$ since $\lambda(\mbold{h}) \leq d$. Hence, {Theorem}~\ref{transLD} provides constants $N_0$ and $c$, determined by $(d, a, \calP)$, such that the Markov transition kernel $\mathcal{Q}$ verifies 
$
\frac{1}{N} \log \mathcal{Q}(h_n , h_{n+1}) = - C(h_n , h_{n+1}) + \eta_n(N) 
$, 
with $\abs{ \eta_n(N) } \leq c / \sqrt{N}$ for all $\mbold{h}$ as above, $n \leq T-1$, and $N > N_0.$  The Markov property yields 
$
\log P( \mbold{H} = \mbold{h} \; | \; H_1 = h_1) = \sum_{n =1}^{T-1} \log \mathcal{Q}(h_n , h_{n+1}). 
$
Since $\lambda(\mbold{h}) = \sum_{n =1}^{T-1} C(h_n , h_{n+1}),$ we obtain
\begin{equation} \label{trajh}
\frac{1}{N} \log P( \mbold{H} = h \; | \; H_1=h_1) = - \lambda(\mbold{h}) + o(1) 
\end{equation} 
with $\abs{o(1)} \leq \sum_{n =1}^{T-1} \abs{ \eta_n(N) } \leq cT / \sqrt{N}$ for $N > N_1$.
\end{proof}
%

%%%%%%%%%%%%%%%%%%%%%%%%%%%%
Next, we prove that the rate function $\lambda(\mbold{H})$ has a uniform  H\"{o}lder property.
%%%%%%%%%%%%%%%%%%
\begin{Theorem} \label{holderlambda}
Fix two constants $d>0$, $a>0$ and the process parameters $\calP$. Consider any (deterministic and fixed) path $\mbold{h}\in \Omega_T $ such that $ b(\mbold{h}) \geq a$ and $\lambda(\mbold{h}) \leq d$. Let $\mbold{H} \in \Omega_T $ be any path such that $\supp(H_n) = \supp(h_n)$ for all $n=1,\dots, T$ and verifying 
\begin{equation} \label{lambda1}
b(\mbold{H}) \geq a \;\text{~and~}\; \vnorm{\mbold{H} - \mbold{h}} \leq a/F_g.
\end{equation}
\hl{Then, there} is a constant $c \equiv c(d, a, \calP)$ such that for all $T$, one has
\begin{equation} \label{lambda2}
\abs{ \lambda(\mbold{H}) - \lambda(\mbold{h})} \leq c T\vnorm{ \mbold{H} - \mbold{h}}^{1/2}.
\end{equation}
\end{Theorem}
\begin{proof}
By Theorem~\ref{holderCost}, there is a constant $c_1(d,a,\calP)$ such that, for all $ n=1,\cdots, T-1$.
$$
\abs{ C(H_n, H_{n+1}) - C(h_n, h_{n+1}) } \leq 
c_1 \left( \vnorm{ H_n - h_n }^{1/2} + \vnorm{ H_{n+1} - h_{n+1}} \right)
$$
and hence $\abs{ C(H_n, H_{n+1}) - C(h_n, h_{n+1}) } \leq 2 c_1 \vnorm{ \mbold{H} - \mbold{h}}^{1/2} $. By definition of $\lambda$, this yields
$$
\abs{\lambda(\mbold{H}) - \lambda(\mbold{h})} \leq \sum_{n=1}^{T-1}\, \abs{ C(H_n, H_{n+1}) - C(h_n, h_{n+1}) } \leq cT \vnorm{\mbold{H} - \mbold{h}}^{1/2},
$$
which proves \eqref{lambda2} with $c = 2 c_1$. 
\end{proof}
%

%%%%%%%%%%%%%%%%%%%
%%%%%%%%%%%%%%%%%%%

%%%%%%%%%%%%%%%%%%%%%%%%%%%%%%%%%%%%
\subsection{Large Deviations for Sets of Trajectories}

Theorem~\ref{thtrajLD} now sets the stage for a large-deviation result for sets of trajectories. In particular, we will show that for large $A$, the probability of observing random paths $\mbold{H} = \{H_1,\cdots, H_n\}$ of population histograms such that $\lambda(\mbold{H}) > A$ is bounded above by $e^{-c A N}$ for some constant $c>0$.  This will naturally lead to the main result of this section, given by Theorem~\ref{thsetLD}.  We first need to define an open neighborhood of paths that is analogous to Definition~\ref{balls} for histograms. 
%%%%%%%%%%%%%%%%%%
\begin{Definition} \label{VNBNpaths}
For any $\Gamma \subset \Omega_T$, define the open neighborhood $V_N(\Gamma)$ as the union of all balls $V_N(\mbold{H})$ with radius $\frac{2}{3 N}$ and center $\mbold{H} \in \Gamma$.  Denote by $B_N(\Gamma)$  the (finite) set of all $N$-rational paths in $V_N(\Gamma)$. Define also $b(\Gamma) = \inf_{\mbold{H} \in \Gamma} b(\mbold{H}).$
\end{Definition}

%%%%%%%%%%%%%%%%%%
\begin{Theorem} \label{thlargecost}
Fix $0 < a < 1$ and  the model parameters $\calP$ (Definition~\ref{P}). Then, there is a constant $c \equiv c(a,\calP) >0$ such that for all $A > c T$ and  for any  initial histogram $Z$ such that $b(Z) \geq a$, the random paths $\mbold{H}$ starting at $H_1=Z$ will verify  
\begin{equation} \label{largecost}
P( \lambda(\mbold{H}) > A \; \text{~and~} \; b(\mbold{H}) > a \; | \; H_1=Z) \leq e^{- \frac{a}{2 T} N A}
\end{equation}
for all $N > c A/T.$
\end{Theorem}
\begin{proof}

%%%%%%%%%%%%%%%%%%%
Recall the notation used to prove Theorem~\ref{transLD}.  Consider any $N$-rational \revision{interior histograms} $H,G \in \mcal{H}(a)$ and $r \in K_N(H)$.  Given $(H_n= H , R_n/N = r) $, the population histogram $J_n$ after growth and mutations is a function  $J_n = J(H, r)$ of $H,r$ (see \eqref{WWn}), with $J(H,r)$ very close to $I = \Psi(H,r)$ (see Equation~\eqref{Jn}), and we defined the set  $A_N$ of acceptable mutation matrices $r=R_n/N$ by
\begin{equation} \label{recallA_N}
A_N = A_N(H,G) = \vset{ r \in K_N(H) \; | \; KL( G, J )< \infty \text{ and }  KL(G, I)< \infty}.
\end{equation}

\hl{We also defined} the approximate transition cost  $\apT(H,r,G) = mut(r, H) + KL( G, J(H,r) )$ in Theorem~\ref{transLD}. The proof of \hl{Theorem}~\ref{transLD} yielded
Equations~\eqref{trans4} and \eqref{recall2}, which provided constants $c_0(a,\calP)$ and $N_0(a,\calP)$ such that for $N > N_0,$
\begin{align} 
&\supp(J(H,r)) = \supp(I(H,r)) \quad \text{with} \quad \vnorm{ J - I }  \leq  c_0 / N, \;\text{~when~}\; r\in A_N ,   \label{JI1}\\
&\mathcal{Q}(H, G)  =  \sum_{r  \in  A_N}  \exp (- N\,  [\, \apT(H,r,G) + o(1) \,] ),   \label{QHG1}
\end{align} 
where $\abs{o(1)} \leq c_0  \log N  / N$. 

 Suppose that $C(H,G) > A$, so that $\tau(H,r,G) \geq C(H,G) > A $ for $r \in K(H)$.  From \eqref{mutbound}, we have a constant $c \equiv c(a,\calP)$ such that $mut(H,r) \leq c$ whenever $b(H) \geq a$ and $r \in K(H)$.  We then have $KL(G,I) = \tau(H,r,G) - mut(H,r) > A - c$.  Hence, \eqref{betabound} yields $\beta(G,I) = \max_{k \in \supp(G)}  {1/I_k(H,r)} \geq KL(G,I) > A - c$.  Therefore, for some $j \in \supp(G),$ one has $I_j(H,r) < 1 / (A-c)$ so that $J_j(H,r) < 1 / (A -c) + c_0 / N$ by \eqref{JI1}. This yields $1/J_j(H,r) > \frac{3}{4} (A - c)$ provided $N >  N_0 + 8 c_0 A $.
A fortiori, we get $\beta(G, J) = \max_{k \in \supp(G)}  {1/J_k(H,r)} > \frac{3}{4} (A - c)$.  
For $N > N_0 + 8 c_0 A $ and $r \in A_N$,
\begin{equation} \label{L(A)}
\apT(H,r,G) \geq KL(G,J) > b(G) \beta(G,J) - \log{g} > \frac{3 a}{4} (A - c) - \log(g) = \frac{3 a A}{4} - c_1 
\end{equation} 
by \eqref{betabound} with  $c_1 = \log(g) + 3 a c /4$. Since $card(E_N) \leq [ (N+1) F_g ]^{g^2}$, combining \eqref{L(A)} and \eqref{QHG1} yields 
\begin{equation} \label{Q1}
\begin{aligned}
\mathcal{Q}(H, G)  & \leq [ (N+1) F_g ]^{g^2} \exp \left( - N \left[ \frac{3 a A}{4} - c_1 \right] + c_0  \log N\right) \\
& \leq \exp \left( N \left[- \frac{3 a A}{4} + c_2 \right]\right)   
\end{aligned}
\end{equation} 
for $N > N_0 + 8 c_0 A$ with $ c_2  =  c_1 + c_0 +g^2 (\log F_g +1)$.

Let $\Gamma(A) = \vset{ \mbold{H}\in \Omega_T \; | \; H_1=Z, \; \lambda(\mbold{H}) > A, \, b(\mbold{H})> a}$. Let $\Gamma_N(A)$ be the set of $N$-rational paths in $ \Gamma(A)$.  Then, $P( \mbold{H} \in \Gamma_N(A) \; | \; H_1 = Z) = P( \mbold{H} \in \Gamma(A) \; | \; H_1= Z)$ since all random paths $\mbold{H}$ are $N$-rational. For $\mbold{H} \in \Gamma_N(A)$, the relation $\lambda(\mbold{H}) > A$ provides at least one time step $\nu = \nu(\mbold{H})  \leq T$ such that $C(H_\nu, H_{\nu+1}) > A/T$. Apply \eqref{Q1} to $(H_\nu, H_{\nu+1})$ to~get
$$
P(\mbold{H} = \mbold{h} \; | \; H_1= Z) \leq P(H_{\nu+1} = h_{\nu+1} \; | \; H_\nu = h_\nu) \leq \exp \left( N \left[- \frac{3 a}{4 T} A + c_2\right]\right) 
$$
for $N > N_0 + 8 c_0 A / T$ and $\mbold{H} \in \Gamma_N(A)$  
and any fixed path $\mbold{h} \in \Gamma_N(A)$.
Since $card(\Gamma_N(A)) \leq (N+1)^{T g}$, this yields
\begin{align*}
P(\mbold{H} \in \Gamma(A) \; | \; H_1=Z) &= P( \mbold{H} \in \Gamma_N(A)\; | \; H_1=Z) \\
&\leq (N+1)^{T g} \exp\left(N\left[- \frac{3 a}{4 T} A + c_2\right]\right) \\
&\leq \exp\left( N\left[- \frac{3 a}{4 T} A + c_3\right]\right)
\end{align*}
for $N > N_0 + 8 c_0 A / T$ with $c_3 = c_2 + T g$. Now, impose $A/T  > 4 c_3 $ to get 
$$ 
P( \mbold{H} \in \Gamma(A) \; | \; H_1=Z)  \leq \exp\left(- \frac{a}{2T}  N A \right) 
$$
{for $N > N_0 + 8 c_0 A/T.$ Set $c = 4 c_3 + N_0 /8 c_0 +16 c_0$, so that the simpler constraints $A/T > c$ and $N > c A/T$ force $A/T > c_3$ and $N > N_0 + 8 c_0 A/T$. This completes the proof.}
\end{proof}
%

%%%%%%%%%%%%%%%%%%
\revision{Theorem~\ref{thsetLD} shows that the probability of observing a trajectory in an open ``tube'' $V_N(\Gamma)$ centered at the set  of trajectories $\Gamma$ is asymptotically equivalent to 
$\exp(-N\Lambda (\Gamma))$. Consequently, the trajectories in $\Gamma$ minimizing $\Lambda$ are the most likely trajectories.}
\begin{Theorem} \label{thsetLD}
Let $\Lambda$ be the large-deviation set functional of the Markov chain of population histograms defined by \eqref{lambda} and \eqref{Lambda}. Fix the process parameters $\calP$, the path length $T,$ and $0<a < 1$.  Denote by $\mbold{H}$ a random histogram trajectory of our Markov chain, starting at some deterministic histogram  $Z$ such that $b(Z) \geq a$. Let $\Gamma \subset \Omega_T$ be any set of histogram trajectories starting at $Z$. Assume that $\Lambda(\Gamma) \leq  L $ for some fixed finite $L>0$, and that $b(\Gamma) \geq a$  (see Definition~\ref{VNBNpaths}).

Let  $V_N(\Gamma)$ be the open neighborhood of $\Gamma$  with (small) radius $\frac{2}{3N}$ (see Definition~\ref{VNBNpaths}). Then, $(T,a, \calP,L)$ determine constants $c$ and $N_0$ such that for all $N > N_0$, one has the  uniform large-deviation result 
\begin{equation} \label{setLD}
\frac{1}{N} \log P( \mbold{H} \in V_N(\Gamma) \; | \;  H_1 = Z) = - \Lambda(\Gamma) + o(1)
\end{equation}
that holds with $\abs{o(1)} \leq c / \sqrt{N}.$ This yields the asymptotic large-deviation limit 
$$
\lim_{N \to \infty}  \frac{1}{N} \log P( \mbold{H} \in V_N(\Gamma) \; | \; H_1 = Z) = - \Lambda(\Gamma).
$$
\end{Theorem}
\begin{proof}
Let $\Gamma$ be any subset of $\Omega_T$ of paths starting at $Z$ with $b(\Gamma) \geq a$ and  finite $\Lambda(\Gamma) = L > 0$. For any $A>0$, let $E(A)$ be the set of all $N$-rational paths $\mbold{h} \in \Omega_T$ with $h_1 = Z$, such that $\lambda(\mbold{h}) > A $ and $b(\mbold{h})> a/2$. Theorem~\ref{thlargecost} provides $c = c(a,\calP)$, such that  
$
P( \mbold{H} \in E(A) \; | \; H_1 = Z) \leq e^{- \frac{a}{2 T} N A } 
$
for all $A > c T$ and $N > c A/T.$ Set $d = T (c + 4 L / a) + L$. Set $A=d$, which forces $A >cT$ and $\frac{a}{2 T} A > 2 L$, so that 
\begin{equation} \label{PED}
P( \omega \in E(d) \; | \; H_1 = Z) \leq e^{- 2 N L } 
\end{equation} 
for all $N > N_1 = c d/T.$ Let $W_N(\Gamma) \subset B_N(\Gamma) \subset V_N(\Gamma)$ be the set of $N$-rational paths $\mbold{h} \in V_N(\Gamma)$ such that $\lambda(\mbold{h}) \leq d$ and $h_1= Z$.  Then, $B_N(\Gamma) \subset W_N(\Gamma) \cup E(d)$, so that \eqref{PED} yields
\begin{equation} \label{set0}
P( \mbold{H} \in B_N(\Gamma) \; | \; H_1 = Z) \leq  P( \mbold{H} \in W_N(\Gamma)  \; | \; H_1 = Z) +  e^{-2 L N}
\end{equation} 
for $N > N_1.$ Due to Theorem~\ref{holderlambda}, there is a constant $c_1 > 0$ determined by $(d,a,\calP)$, and thus by $(T, L, a, \calP),$ such that for any paths $\mbold{h}, \mbold{h}' \in \Omega_T,$ the inequality
\begin{equation} \label{lambda11}
\abs{\lambda(\mbold{h}') - \lambda(\mbold{h})} \leq c_1 \vnorm{\mbold{h}' - \mbold{h}}^{1/2}
\end{equation}
must hold whenever the following holds for all integers $1 \leq n \leq T$:
\begin{align} 
&\min\vset{\lambda(\mbold{h}), \lambda(\mbold{h}')} \leq d,    \label{lambda1a} \\
b(\mbold{h}') \geq a&, \qquad \vnorm{\mbold{h}' - \mbold{h}} \leq a/F_g, \qquad \supp(h'_n) = \supp(h_n) \label{lambda2a}.
\end{align}

\hl{Since} $b(\Gamma) \geq a$, Lemma \ref{basicBN} implies $b(V_N(\Gamma)) \geq a/2$, so that  all $\mbold{h} \in W_N(\Gamma)$ verify $\lambda(\mbold{h}) \leq d$ and $b(\mbold{h}) > a/2$. Therefore, Theorem~\ref{thtrajLD} provides $c_2$ and $N_2 > N_1$ determined by $(d, a, \calP),$ and thus by $(T, L, a, \calP)$, such that
\begin{equation} \label{set1}
\frac{1}{N} \log P( \mbold{H} = \mbold{h} \; | \; H_1 = Z) = - \lambda(\mbold{h}) + o_1(1)
\end{equation}
for all $N > N_2$ and any $\mbold{h} \in W_N(\Gamma)$ with $\abs{o_1(1)} \leq T c_2 / \sqrt{N}.$ For each $\mbold{h} \in W_N(\Gamma) \subset V_N(\Gamma),$ one can select a path $\hat{\mbold{h}} \in \Gamma$ such that $\vnorm{\hat{\mbold{h}}- \mbold{h}} \leq 1/N$. 
By Lemma \ref{basicBN}, the paths $\mbold{h}$ and $\mbold{h}' =\hat{\mbold{h}}$ verify both \eqref{lambda1a} and \eqref{lambda2a} so that \eqref{lambda11} applies and yields 
$
\abs{ \lambda(\hat{\mbold{h}}) - \lambda(\mbold{h})} \leq c_1 \vnorm{ \hat{\mbold{h}} - \mbold{h}}\leq c_1 / \sqrt{N}
$ 
provided $N > F_g / a.$  This implies $ - \lambda(\mbold{h}) +  o_1(1) = - \lambda(\hat{\mbold{h}}) + o_2(1)$ with $\abs{ o_2(1)} \leq (c_1 + T c_2) / \sqrt{N}$ for $N >N_3 =N_2 + F_g /a$.  For $N>N_3$ and $\mbold{h} \in W_N(\Gamma)$, this yields
$$
P( \mbold{H} = \mbold{h} \; | \; H_1 = Z) = \exp (- N \lambda(\hat{\mbold{h}}) + N o_2(1)) \leq \exp (- N L + c_3 \sqrt{N})
$$
with $c_3= c_1 + T c_2$ due to \eqref{set1} and $\lambda(\hat{\mbold{h}}) \geq \Lambda(\Gamma) = L.$ For $N > N_3$, sum over $\mbold{h} \in W_N(\Gamma)$
to get 
$$
\begin{aligned}
P( \mbold{H} \in W_N(\Gamma) \; | \; H_1 = Z) & \leq \exp [- N L + c_3 \sqrt{N} + gT  \log(N+1)] \\
& \leq \exp (- N L + (c_3 + gT) \sqrt{N})
\end{aligned}
$$
since $card(W_N(\Gamma)) \leq (N+1)^{gT}.$  Due to \eqref{set0}, this gives 
\begin{align}\label{PBNup}
P( \mbold{H} \in B_N(\Gamma) \; | \; H_1 = Z) &\leq \exp\left[- N L + (c_3 + T g) \sqrt{N}\right] + e^{-2 N L} \nonumber \\
&\leq 2 \exp\left[- N L + (c_3 + gT) \sqrt{N}\right]
\end{align}
for $N > N_3.$ One has $P( \mbold{H} \in V_N(\Gamma) \; | \; H_1 = Z) = P( \mbold{H} \in B_N(\Gamma) \; | \; H_1 = Z)$ since the random paths $\mbold{H}$ are always $N$-rational. Hence, for $N > N_3$, \eqref{PBNup} yields the upper bound
\begin{equation} \label{logPVNup}
\frac{1}{N} \log P( \mbold{H} \in V_N(\Gamma) \; | \; H_1 = Z) \leq - L + c_4 / \sqrt{N}
\end{equation}
with  $c_4 =c_3 + gT +1.$ Set $c_5 = \left[\frac{1+c_1}{T (c + 4 L/a)}\right]^2$ and $N_4 = N_3 + c_5$. This yields $L +(1 + c_1) / \sqrt{N} \leq d = L +T (c + 4 L / a)$ for $N > N_4.$  For each $N,$ select $\mbold{h}'' \equiv \mbold{h}''(N) \in \Gamma$ such that $L \leq \lambda(\mbold{h}'') \leq L + 1 / \sqrt{N}$. For $N > N_4$, one has $\lambda(\mbold{h}'') \leq d $, and \eqref{lambda11} applies to the pair of trajectories $\mbold{h}''$  and $\mbold{h}'$ for any $N$-rational $\mbold{h}' \in V_N(\mbold{h}'')$ to give
$
\abs{\lambda(\mbold{h}') - \lambda(\mbold{h}'') }  \leq c_1 / \sqrt{N}.
$
Hence, $\lambda(\mbold{h}') \leq L + (1+c_1) / \sqrt{N} \leq d$, which shows that $\mbold{h}' \in W_N(\Gamma)$.  Apply \eqref{set1} to $\mbold{h}' \in W_N(\Gamma)$ to obtain
\begin{equation}\label{lowPh'}
\frac{1}{N} \log P( \mbold{H} = \mbold{h}' \; | \; H_1 = Z) = - \lambda(\mbold{h}') + o_2(1) \geq - L - (1 + c_1 +T c_2) / \sqrt{N}
\end{equation}
for $N> N_4.$ Since $\mbold{h}' \in W_N(\Gamma) \subset V_N(\Gamma)$, one has 
$$
P( \mbold{H} \in V_N(\Gamma) \; | \; H_1 = Z) \geq P( \mbold{H} \in W_N(\Gamma) \; | \; H_1 = Z) \geq P( \mbold{H} = \mbold{h}' \; | \; H_1 = Z). 
$$
\hl{Combining this} with \eqref{lowPh'} yields the lower bound
\begin{equation} \label{logPVNlow}
\frac{1}{N} \log P( \mbold{H} \in V_N(\Gamma) \; | \; H_1 = Z) \geq - L - T (1 + c_1 + T c_2) / \sqrt{N}
\end{equation}
for $N>N_4.$ Combining \eqref{logPVNup} and \eqref{logPVNlow} concludes the proof.
\end{proof}

Theorem~\ref{thsetLD} is crucial for computing the most likely random evolutionary path  
$\mbold{h}^*$ 
of a population starting with a known initial histogram $H$ and ending with another known histogram $G$ after an arbitrary number of ``days''. Theorem~\ref{thsetLD} indicates that the likelihood of observing the random histogram trajectory $\mbold{H}$ within any fixed set of paths $\Gamma \subset \Omega_T$ is roughly of the order of $\exp(- N\, \Lambda(\Gamma))$, which is an extremely small exponential for a large $N$ whenever $\Lambda(\Gamma) >0$.
Therefore, if $\Lambda(\Gamma)>0$, then observing a random trajectory $\mbold{H}\in \Gamma$ 
is a rare event. The actual calculation of $\Lambda(\Gamma)$ involves minimizing the rate function $\lambda(\mbold{h})$ over all paths $\mbold{h}\in \Gamma.$  Let $\Gamma^* \subset \Gamma$ be the set of all such minimizing paths for the rate function $\lambda(\mbold{h})$.  Whenever the rare event $\mbold{H} \in \Gamma$ is realized by some observed random histogram trajectory $\mbold{H}$, then with a conditional probability extremely close to 1, the trajectory $\mbold{H}$ must follow very closely one of the minimizing paths $\mbold{h} \in \Gamma^*$. We explicitly quantify this argument in Section~\ref{mostlikely}.

%%%%%%%%%%%%%%%%%%%%%%%%%%%%%%%%%%%%
\subsection{Mean Evolution and Zero-Cost Trajectories}\label{zerotraj}

When $\Gamma$ is a thin open tube around a zero-cost path $\mbold{h}$, i.e., a path such that $\lambda(\mbold{h}) = 0$, then $\Lambda(\Gamma) =0$;  and our next theorem implies that $P(\mbold{H} \in \Gamma)$ will approach $1$ at exponential speed for large $N$. This motivates our explicit analysis and construction of all zero-cost paths. In particular, we now show that given an initial histogram $h_1$, the unique zero-cost trajectory starting at $h_1$ is recursively computable by an explicit formula.

%%%%%%%%%%%%%%%%%%
\begin{Theorem} \label{thzerotraj}
Fix a path length $T$. A histogram path $\mbold{h}= [ h_1\, \dots, h_T] \in \Omega_T$  satisfies $\lambda(\mbold{h}) = 0$ if and only if  $h_{n+1} = \zeta(h_n)$ for $1 \leq n \leq T-1$, where the histogram-valued  function $H \to \zeta(H) \in \mathcal{H}$ is defined for all $H \in \mathcal{H}$ by
\begin{equation} \label{zeta}
\zeta_j(H)  = \frac{1}{\langle F, H \rangle} \left( F_j H(j) - m \sum_k  q_{j,k} F_j H(j) + m \sum_k  q_{k,j} F_k H(k) \right), \quad j=1,\ldots,g.
\end{equation}
\hl{Hence, a zero-cost} path $\mbold{h}$ is uniquely determined by its starting point $h_1$. Consider any zero-cost path $\mbold{h}$. Let $V_N(\mbold{h})$ be the ball with a center at $\mbold{h}$ and radius $\frac{2}{3 N}$  in $\Omega_T$. Then, the initial histogram $h_1$ and the process parameters $\calP$ determine constants $c$ and $N_0$ such that the random path $\mbold{H} \in \Omega_T$ verifies
\begin{equation}  \label{V(h)}
1 \geq P( \mbold{H} \in V_N(\mbold{h})  \; | \; H_1 = h_1 ) \geq 1 - e^{- c  \sqrt{N}} 
\end{equation}
for all $N > N_0.$
\end{Theorem}
\revision{Given the histogram $h_1$ in the zero-cost path, formula \eqref{zeta} can be used to recursively compute the sequence $h_1 \to h_2 \to\cdots \to h_n$. The frequency of the $j$-th genotype in the zero-cost path is obtained from the growth term $F_j H(j) / \langle F, H \rangle$, minus the expected number of emigrants from genotype $j$, plus the expected number of immigrants into genotype $j$.}

\begin{proof}
Given any two histograms $H,G\in \mcal{H}$ with finite cost $C(H,G)$,  Lemma \ref{finitecost} proves the existence of a matrix $r \in \overline{K(H)}$ such that $C(H,G) = \tau(H, r, G)$. Hence, $C(H,G) = 0$ if and only if $\tau(H,r,G) = 0,$ which is equivalent to $mut(H,r) = KL(G, \Psi(H,r) ) = 0$.  From Proposition~\ref{mutproperties}, 
one has $mut(H,r) = 0$ iff 
\begin{equation} \label{r.optimal}
r_{j,k} = m q_{j,k} F_j H(j) = m\, q_{j,k} F_j H(j) 
\end{equation}
for all $j,k.$ By definition of the Kullback--Leibler divergence, one has  $KL(G, \Psi(H,r) ) = 0$ if and only if $\Psi(H,r) = G$. Combine this relation with \eqref{psi} and \eqref{r.optimal} to conclude that $C(H,G) = 0$ if and only if one has
$$
G(j) = \frac{1}{\langle F, H \rangle}\left( F_j H(j) - \sum_k  m\, q_{j,k} F_j H(j) + \sum_k  m\, q_{k,j} F_k H(k) \right)
$$
for all $j.$ Hence, $C(H,G) = 0$ if and only if $G= \zeta(H)$. Now, for any path $\mbold{h} \in \Omega_T$, the relation $\lambda(\mbold{h}) = 0$ holds if and only if $C(h_n, h_{n+1} ) = 0$ for all $1 \leq n \leq T-1$, which is equivalent to $h_{n+1}= \zeta(h_n)$ for all $1 \leq n \leq T-1.$  Finally, the bound in \eqref{V(h)} is an immediate consequence of \eqref{setLD} applied to the set $F = \vset{\mbold{h} \;|\; \lambda(\mbold{h})= 0}.$
\end{proof}

Next, we demonstrate that as $n \to \infty,$ any infinite zero-cost path $\mbold{h} = [h_n]_{n=1}^{\infty}$ achieves near fixation of some explicitly determined genotype.  
%%%%%%%%%%%%%%%%%%
\begin{Definition}\label{reachgen}
For any non-empty set $S$ of genotypes, define the set $R(S)$ of \emph{genotypes reachable from S} as the set of all genotypes such that there is a genotype sequence $k_1,\cdots, k_T$ of arbitrary length $T$ such that $k_1 \in S, k_T = j,$ and all $q_{k_t, k_{t+1}} > 0.$
\end{Definition}
%
%
%%%%%%%%%%%%%%%%%%
\begin{Theorem} \label{hinfty} 
Let $\mbold{h}$ be any zero-cost histogram path of infinite length starting at $h_1 \in \mathcal{H}$. Denote $S_n =\supp (h_n)$. Then, the set $S_n$ increases with $n$ and stabilizes after a finite time $n_0$ so that $S_n= \Ster$ for all $n > n_0$. Moreover, $\Ster = R(S_1)$, where $R(S_1)$ is the set of all genotypes reachable from $S_1$, as in Definition~\ref{reachgen}. Let $s \in \Ster$ be the fittest genotype within $\Ster$, i.e., such that $F_s = \max_{j\in \Ster} F_j$. Then, the initial histogram $h_1$ and the parameters $\calP$ determine $m_0 >0$ such that for all mutation rates $m \leq m_0$,    
\begin{equation}\label{fixation}
\lim_{n \to \infty}  h_n = \hter 
\end{equation}
where $\hter$ is the unique solution of $ \hter = \zeta(\hter)$, has support $\supp(\hter) = \Ster,$ and is a $C^{\infty}$ function of $m$ with the first-order expansion
\begin{align} \label{fixedpoint}
\hter(i) &\simeq m  q_{s, i} / (F_s - F_i) \quad \text{for} \quad i \in (\Ster \setminus s), 
\\
\hter(s) &\simeq 1- m\sum_{ i \in (\Ster \setminus s)} q_{s, i} / (F_s - F_i).
\end{align}
\end{Theorem}
\revision{If all genotypes are reachable from the initial histogram $h_1$ along the zero-cost path, then Theorem~\ref{hinfty} and formula \eqref{fixedpoint} determine the limiting histogram as time tends to infinity. In particular, formula \eqref{fixedpoint} describes the balance between selection through genotype fitness and mutation.}
\begin{proof}
The detailed proof is displayed in Appendix \ref{sec:ap4}.
\end{proof}

%%%%%%%%%%%%%%%%%%%%%%%%%%%%%%%%%%%%
\section{Most Likely Evolution from Initial to Terminal Histograms} \label{mostlikely}
Theorem~\ref{thsetLD} indicates that for any  set $\Gamma \subset \Omega_T$ of histogram paths,  the probability $P( \mbold{H} \in \Gamma)$ that the random histogram trajectory $\mbold{H}$ belongs to $\Gamma$ is \textit{roughly} estimated by $\exp(- N\, \Lambda(\Gamma))$.
We now apply this estimate to the set $\Gamma(H,G)$ of all histogram trajectories $\mbold{H} =\{H_1,\dots,H_T\}$ such that $H_1= H$ and $H_T =G$, where  the initial and terminal histograms $H$ and $G$ are known and fixed. For large $N$ and any trajectory $\mbold{h}\in\Gamma(H,G)$, we have shown that $\frac{1}{N} \log ( P( \mbold{H} = \mbold{h}) )$ tends to $ - \, \lambda(\mbold{h}) $. This indicates that the most likely evolutionary path $\mbold{h}^*$ ``followed'' by a random histogram trajectory $\mbold{H}^*$ connecting $H$ to $G$ in $T$ steps should be a minimizer of  $\lambda(\mbold{h})$ over all paths $\mbold{h}\in \Gamma(H,G).$ This is a natural question of interest for laboratory experiments on bacterial genetic evolution, since actual detailed observations of $H_n$  may be costly and hence may occur only every $T$ steps. In long experiments, the time $T$ itself may not even be known precisely, so the reconstruction of the whole most likely evolutionary path connecting $H$ to $G$ given only $H$ and $G$ should provide insights into key genotype changes necessary to realize a rare evolution from $H$ to $G$. We also outline how to compute the most likely number of days  $T$ enabling a rare transition from $H$ to $G$ in $T$ steps.

%%%%%%%%%%%%%%%%%%%%%%%%%%%%%%%%%%%%
\subsection{Interior Paths in the Space of Histogram Trajectories} \label{interior}

Note that for interior histograms $H,G \in \interior{\mathcal{H}}$, the transition cost $C(H,G)$ is always finite. A path $\mbold{H} = [H_1\cdots H_T ] $ will be called an \emph{interior path} if $H_n\in \interior{\mcal{H}}$ are interior histograms for all $1 \leq n \leq T$.  To develop explicit computational schemes, 
we focus on interior paths and histograms for the remainder of this paper.

%%%%%%%%%%%%%%%%%%%%%%%%%%%%%%%%%%%%
\subsection{Sets of Thin Tubes Realizing Rare Events}

\vspace{6pt}
\begin{Definition} \label{U^eta}
For any  small $\eta>0$ and any set of paths  $\Gamma \subset  \Omega_T$, we define the $\eta$-neighborhood $U^{\eta}(\Gamma)$ of $\Gamma$ as the union of all open balls of radius $\eta$ centered at arbitrary $\mbold{h}$ belonging to $\Gamma$. For small $\eta$, the set  $U^{\eta}(\Gamma)$ is a set of thin tubes of paths with ``axes'' $\mbold{h} \in \Gamma$.
\end{Definition}
Let $\mbold{H} \in \Omega_T $ be the random trajectory of population histograms.  For any closed set of interior paths $E \subset \Omega_T$, one can easily prove that $a = b(E) >0$. Assume,  moreover, that $\Lambda(E) > 0$. Then, for any $0 <A < \Lambda(E)$ and large enough $N > N_0(a, A, \calP)$, Theorem~\ref{thsetLD} provides the fast-vanishing bound
$$
P(\mbold{H} \in E) \leq P(\mbold{H} \in V_N(E)) \leq 2 e^{- N A}
$$ 
 so that the random events $\{\mbold{H} \in E\}$ and  $\{\mbold{H} \in V_N(E)\}$ are both \emph{\hl{rare events.}}

We now show (Theorem~\ref{conditioning}) that for large population size $N$, rare evolutionary events $\{\mbold{H} \in E\}$ with finite $\Lambda(E) > 0$ can only be realized by random population evolutions following a very thin tube of radius $1/N$ around one of the (possibly several) paths $\mbold{h}^* = \arg\min_{\mbold{h} \in E} \lambda(\mbold{h})$.

More precisely, let $\mcal{T} = B_N(\mbold{H})$ be the thin tube of trajectories centered around the random trajectory $\mbold{H}$ and having vanishingly small radius $1/N$. Fix a set of paths $E$ such that $\Lambda(E)>0$, so that the event $\mbold{H} \in E$ is a rare event.  The next theorem proves that whenever  we know that the thin random tube $\mcal{T}$ contains a path belonging to $E$, then with conditional probability extremely close to 1,  the random trajectory $\mbold{H}$
itself must be arbitrarily close to a trajectory $\mbold{h}^*$ that minimizes 
the rate function  over all paths in $E$.
Of course, this result is of greater interest when there is only one such minimizer $\mbold{h}^*$, but the uniqueness 
of $\mbold{h}^*$ cannot be guaranteed in general.
Computing such minimizing paths requires efficient numerical strategies, as discussed in Section~\ref{s:numtraj}.

\begin{Theorem} \label{conditioning} Fix the path length $T$ and an initial interior histogram $H$. Denote $P_H$ the probability distribution of random histogram paths $\mbold{H} \in \Omega_T$ starting at H. Let $E \subset \Omega_T$ be any closed set of interior paths starting at $H$ and satisfying $0 < \Lambda(E) < \infty $. Let $E^*\subset E$ be the set of all paths $\mbold{h}^*$ minimizing the rate function $\lambda(\mbold{h})$ over all $\mbold{h} \in E$. Then, $E^*$ is a closed subset of $E$. For any fixed $ \eta > 0$, the $\eta$-neighborhood $ U^{\eta}(E^*)$ of $E^*$ verifies 
\[
\lim_{N \to \infty}  P_H \left(\mbold{H} \in U^{\eta}(E^*) \; | \; \mbold{H} \in V_N(E) \right)   =  1
\]
with exponential speed of convergence.
\end{Theorem}
\begin{proof}
Set $a = b(E) > 0 $ and $L = \Lambda(E) > 0$. Let  $E(c) = \vset{ \mbold{h} \in E \; | \; \lambda(\mbold{h}) \leq c }$. Then, the function $\lambda(\mbold{h})$ is continuous on $E(2 L)$ due to Theorem~\ref{holderlambda} applied to interior paths. Hence, $E(2 L)$ is closed and must contain any path minimizing  $\lambda(\mbold{h})$ over $E(2 L).$ However, the two sets of minimizers of $\lambda$ over $E$ and over $E(2 L)$ are obviously identical. This proves $E^* \subset E(2 L) \subset E$.

Theorem~\ref{holderlambda} applied to interior paths provides a constant $c = c (T,L,a,\calP)$ such that for all $\mbold{h} \in E^*$ and $\mbold{h}' \in E$ with $\vnorm{ \mbold{h}' - \mbold{h} } < a/F_g$, one has $\abs{ \lambda(\mbold{h}') - \lambda(\mbold{h})} \leq c \vnorm{\mbold{h}'-\mbold{h}}^{1/2}.$ Fix any $0 < \eta < a/F_g$. Let $K \subset E$ be the open $\eta$-neighborhood  of $E^*$ within $E,$  and set $W= E - K.$ For each $\mbold{h}' \in K$, there is one $\mbold{h} \in E^*$ with $\vnorm{\mbold{h}' - \mbold{h}} < \eta.$ Therefore, Theorem~\ref{holderlambda} implies $\lambda(\mbold{h}') < L + c \sqrt{\eta}.$  This forces $\lambda(\hat{\mbold{h}}) \geq L + c \sqrt{\eta}$ for all $\hat{\mbold{h}} \in W$ so that $\Lambda(W) \geq L + c \sqrt{\eta}$.  Apply Theorem~\ref{thsetLD} to get $N_0$ and $c_0$, determined by $(T,L,a,\calP,\eta)$, such that for $N >N_0$, 
\begin{align*}
\frac{1}{N}  \log P_H (\mbold{H} \in V_N(W) ) & =  - L - c \sqrt{\eta}+o_1(1), \\
\frac{1}{N}  \log P_H (\mbold{H} \in V_N(E) )  & =  - L +o_2(1), 
\end{align*} 
with  $\abs{o_1(1)} < c_0 / \sqrt{N}$ and $\abs{o_2(1)} < c_0 / \sqrt{N}.$ These  results yield
\begin{align}
P_H (\mbold{H} \in V_N(W)) & \leq \exp(-N L - N c \sqrt{\eta} +c_0 \sqrt{N}),\\
P_H (\mbold{H} \in V_N(E)) & \geq \exp( -N L  - c_0 \sqrt{N}),
\end{align}
so that
$$
\frac{P_H (\mbold{H} \in V_N(W))}{P_H (\mbold{H} \in V_N(E))}  \leq \exp(- N c \sqrt{\eta} + 2 c_0 \sqrt{N}).
$$

\hl{Impose} $N > N_0 + {16 c_0^2}{c^2 \eta}$ to force  $N (c/2) \sqrt{\eta} > 2 c_0 \sqrt{N}$ to give 
$$
\frac{P_H (\mbold{H} \in V_N(W))}{P_H (\mbold{H} \in V_N(E))} \leq \exp(- N (c/2) \sqrt{\eta}).
$$

\hl{Since} $E = W \cup K,$ one has $V_N(E) \subset  V_N(W) \cup V_N(K).$ Therefore, 
$$
P_H (\mbold{H} \in V_N(K)) \geq P_H (\mbold{H} \in V_N(E)) - P_H (\mbold{H} \in V_N(W)).
$$

\hl{From} $ P_H (\mbold{H} \in V_N(K) \; | \; \mbold{H} \in V_N(E)) = \frac{P_H (\mbold{H} \in V_N(K))}{P_H (\mbold{H} \in V_N(E))}$, we now obtain
$$
1 \geq P_H (\mbold{H} \in V_N(k) \; | \; \mbold{H} \in V_N(E) ) \geq 1 - \frac{P_H (\mbold{H} \in V_N(K))}{P_H (\mbold{H} \in V_N(E))} \geq 1 - \exp\left(- N (c/2) \sqrt{\eta}\right)
$$ 
for $N > N_1.$ For $ N > N_1 +1/\eta$, one has $V_N(K) \subset U^{2 \eta}(E^*) \cap  V_N(E).$ Therefore,
$$
1 \geq P_H (\mbold{H} \in U^{2 \eta}\; | \; \mbold{H} \in V_N(E) ) \geq P_H (\mbold{H} \in V_N(K) \; | \;  \mbold{H} \in V_N(E)) \geq 1 - \exp\left(- N (c/2) \sqrt{\eta}\right),
$$
which concludes the proof.
\end{proof}

%%%%%%%%%%%%%%%%%%%%%%%%%%%%%%%%%%%%
\subsection{Most Likely Path Connecting Two Histograms} \label{mostlikelypath}

\vspace{6pt}
\begin{Definition} \label{reachable}
Let $Q = (q_{i,j} )$ be the mutation transfer matrix and $m$ the mutation rate. As in Definition~\ref{reachgen}, for any subset $S $ of the set of genotypes $\{ 1,\cdots , g\}$ we denote by $R(S)$ the set of genotypes reachable from $S$ in some finite number of steps.  For any given histograms $H$ and $G$, we say that {G is reachable from H} if $\supp(G) \subset R(\supp(H))$.
\end{Definition}
Note that if there exists some power $p$ such that all coefficients of $Q^p$ are positive, one has $R(S) = \vset{ 1, \dots, g }$ for any non-empty set of genotypes $S$ and hence any  $G$ is reachable from any $H$.  
The next theorem answers an important question for bacterial genetic evolution: how can one reconstruct the \hl{most likely}
evolutionary path starting at a known initial histogram $H$ and reaching a known terminal histogram $G$ after a fixed number $T$ of daily cycles.

\begin{Theorem}\label{tube.fixedT} Fix a path length $T$ and a small $0 < a <1$. Fix any interior histograms $H$ and $G$ such that $G$ is reachable from $H$, and verifying $b(H) >a$, $b(G) >a$.  Let $B= B(G,\eta)$ be an open ball of center $G$ and a small radius $\eta <a/2$. Define the following compact set $E_T$ of paths  starting at $H$ and reaching $B$ for the first time at time T, by
\begin{equation} \label{E_T}
E_T = \vset{\mbold{h} \in \Omega_T \,|\, b(\mbold{h}) \geq a ;  \, h_1 = H  ; \,  h_T \in \bar{B} ; \, h_n \neq B(G) \text{ for } 1 \leq n < T}.
\end{equation}
\hl{Then, for} $\eta$ small enough,  $\Lambda (E_T) = \inf_{\mbold{h} \in E_T} \, \lambda(\mbold{h})$ must be finite. The set $E_T^*$ of all paths $\mbold{h} \in E_T$ such that $ \lambda(\mbold{h}) = \Lambda (E_T)$ is then non-empty and contains only interior paths.

When $\Lambda (E_T) >0$, then the event $\mcal{H} \in E_T$ is a rare event, and  for any open neighborhood U of $E_T^*$, one has
$$
\lim_{N \to \infty}  \frac{1}{N} \log P (\mbold{H} \in U\; | \; \mbold{H} \in E_T) = 1
$$
with convergence at exponential speed.

When $\Lambda (E_T) =0$, then there is a unique zero-cost trajectory $\mbold{h}^* \in E$ determined by \hl{Theorem~}
\ref{thzerotraj},
with $h_1 =H$ and $h_T\in\bar{B}$. Then, for any open neighborhood U of $\mbold{h}^* $, one has
$$
\lim_{N \to \infty}  P (\mbold{H} \in U) =  \lim_{N \to \infty}  P (\mbold{H} \in E_T) = 1
$$
with convergence at an exponential speed.

\end{Theorem}
\begin{proof}
We omit the detailed proof which is a direct application of the preceding Theorem~\ref{conditioning}, and of our results on zero-cost trajectories.
\end{proof}
\textbf{\hl{Intuitive interpretation of theorem}
~\ref{tube.fixedT}: } For a large $N$ and a very small $\eta$, the preceding theorem indicates that a stochastic population evolution connecting interior histograms $H$ and $G$ in exactly $T$ days must follow a very thin tube of paths centered around an  ``optimal'' path $\mbold{h}^*$ that minimizes $\lambda(\mbold{h})$ over all $\mbold{h}$ such that $h_1=H$, $h_T =G$, and $h_j\neq G$ for $j <T$. We discuss this minimization problem in the next section.

For laboratory experiments observing bacterial populations over a very long time (e.g., several years), a plausible situation may involve two histograms $H$ and $G$ where $H$ is observed on day $n=1$, but the exact observation time $T$ of $G$ is only known to be smaller than some given $T_{max}$. The following theorem indicates how one can estimate $T$.

\begin{Theorem}\label{tube}
We use here the same notation $a,H,G,\eta, E_T, B=B(G,\eta)$ and hypotheses as in the preceding Theorem~\ref{tube.fixedT}, but $T$ can now vary in a fixed interval $[0,T^+]$. For any path $\mbold{h} \in \Omega_{T^+}$, denote by $\theta(\mbold{h})$ the first hitting time of the small ball $B = B(G,\eta)$.

Denote by $E = \bigcup_{T\leq T^+} E_T$ the set of all paths starting at $H$ with the hitting time $\theta\leq T^+$.
Then, $\Lambda (E) = \inf_{T\leq T^+} \, \Lambda(E_T)$ is  finite and there is at least one  time $T^*<T^+$ such that $\Lambda (E) = \Lambda(E_{T^*})$. If $T^*$ is unique and $\Lambda(E_{T^*}) >0$, then one has 
$$
\lim_{N \to \infty}  \frac{1}{N} \log P (\theta(\mbold{H}) =T^* \; | \; \mbold{H} \in E ) = 1
$$
with convergence at an exponential speed.
\end{Theorem}
\begin{proof}
 This result is a direct corollary of Theorem~\ref{conditioning}, since a finite union of rare events $E_T$ can essentially ``only be realized'' if the $E_T$ with the smallest $\Lambda(E_T)$ is realized.
\end{proof}
%
%
%%%%%%%%%%%%%%%%%%%%%%%%%%%%%%%%%
%%%%%%%%%%%%%%%%%%%%%%%%%%%%%%%%%%%%%%
\section{Recursive Computation of Cost-Minimizing Trajectories Connecting Two Histograms}
\label{s:numtraj}
%\\
The previous section on large deviations demonstrates that given two interior histograms $H$ and $G$, the most likely bacterial population evolution from $H$ to $G$ in a given number of days, $T$, must be a very thin tube of histograms trajectories centered around a deterministic path $\mbold{h}^* = [h_1^* \, h_2^* \dots h_T^*]$ minimizing the rate function $\lambda(\mbold{h})$ over all $\mbold{h} \in \Omega_T$ such that $h_1 = H$ and $h_T = G$. 
By analogy with control problems and the definition of energy minimizing geodesics $\gamma(t)$ on Riemannian manifolds, we call any such $\mbold{h}^*$ a \emph{geodesic} from $H$ to $G$ whenever $\lambda(\mbold{h}^*)$ is finite. 
Indeed, in this discrete version of geodesics the one-step cost $C(h_n, h_{n+1})$ at the discrete time $n$ is analogous to the Riemannian energy $||\frac{d \gamma}{dt} ||^2$ at a continuous time $t$ (see our companion paper~\cite{azencott2023rare} where this point of view is explored further). When $h_n^* \in \mcal{H}^\circ$ for all $1 \leq n \leq T$, we call $\mbold{h}^*$ an \emph{interior geodesic}. The numerical computation of geodesics presents multiple numerical and mathematical challenges in our discrete context. We now develop a key recursive equation in reverse time that must be verified by all geodesics. This equation is an essential practical tool for the efficient numerical computation of geodesics connecting two histograms $H$ and $G$.

\begin{Theorem}\label{MinPathInt}
Fix any small $a>0$ and the process parameters $\mcal{P}$. If the mutation rate $m$ is small enough, then for any pair $y, z \in \mcal{H}^\circ$ of interior histograms verifying $b(y) \ge a$ and $ b(z)\ge a$, there is a unique interior histogram 
$x=\chi(m,y,z) $ such that $[x,y,z]$ is an interior geodesic connecting $x$ to $z$
in two steps. There is a positive constant  $m_0 \equiv m_0(a,\mcal{P})$ such that, for all $0 < m <m_0$, the following hold:
\begin{enumerate}
    \item The histogram-valued function $x= \chi(m,y,z)$ is $C^{\infty}$ in $m,y,z$;  
    \item The one-step costs $C(x,y)$ (from $x$ to $y$) and $C(y,z)$  (from $y$ to $z$) are both differentiable in $y \in \mcal{H}^\circ$;
    \item $x=\chi(m,y,z) $ is the unique histogram solving $grad_y \left[ C(x,y) + C(y,z) \right] = 0 $, where the gradient $grad_y$
is computed in the convex set of interior histograms $x,y \in \mcal{H}^\circ$; 
    \item For any $T \ge 2$, any  interior geodesic $\mbold{h}= [h_1 \, h_2 \dots h_T]$ such that $b(\mbold{h}) >a$ must verify the reverse recurrence relation 
\begin{equation}\label{reverse}
h_n = \chi(m, h_{n+1}, h_{n+2}), \qquad \text{for} \;  1 \leq n \leq T-2;
\end{equation}
    \item The geodesic $\mbold{h}$ is completely determined by its last two histograms, $h_T$ and $h_{T-1}$;
    \item 
{The first-order Taylor expansion in $m$ of $x = \chi(m,y,z)$  is given by  
$x = \wh{x} + m \left<\wh{x} , w \right> + O(m^2)$, 
where the interior histogram $\wh{x}$, given by \eqref{ZeroOrder}, and the vector $w \in R^g$ both  depend only on the histograms $(y, z)$,  and are explicitly computed below by the succession of the five formulas~\eqref{X}--\eqref{betas}.}
\end{enumerate}
\end{Theorem}
\begin{proof}
Any sub-segment $[ h_n = x , h_{n+1} = y , h_{n+2} = z ]$ of an interior geodesic $\mbold{h}$ is also a geodesic connecting  $x =h_n $ to $z = h_{n+2}$ in two steps. For the  two-step interior geodesic $(x, y, z)$, the two-step cost function $u(x,y,z) = C(x,y) + C(y,z)$ is minimized in $y$ by $y = h_{n+1}$. For any three histograms $x,y,z \in \mcal{H}^\circ,$ both $C(x,y)$ and $C(y,z)$ are finite and differentiable in $y$ by Theorem~\ref{thCHG} and convex in $y$ by Lemma \ref{finitecost}. Hence, for fixed $x,z \in \mathcal{H}^\circ$, the function $f(y) = C(x,y) + C(y,z)$ is finite, convex, and differentiable for all $y$ in the open convex set $\mathcal{H}^\circ$. If $y \in \mathcal{H}^\circ$ is a minimizer of $f(y)$ over all $\mathcal{H}^\circ$, any vanishingly small $\Delta y$ ($\Delta y$ is a histogram in this context) so that $y + \Delta y \in \mathcal{H}^\circ$ must verify $\sum_j \Delta y(j) =0$ and $f(y+\Delta y) \geq f(y)$. Hence, the gradient of $\left[f(y) + \mu  \sum_j  y(j) \right]$ must be 0 for some Lagrange multiplier $\mu$. For each $s\in \vset{1, \dots, g}$, this yields the system
\begin{equation}\label{twosteps}
\frac{\partial}{\partial y_s} f(y) + \mu= 0.
\end{equation}

\hl{Denote} $u(x,y,z) = f(y) = C(x,y) + C(y,z)$. For given  $y, z \in \mathcal{H}^\circ$, extend \eqref{twosteps} into the following system of $(1+g)$ equations to be solved for a histogram $x$ and a Lagrange multiplier $\mu$:
\begin{align}
\frac{\partial}{\partial y_s} u(x,y,z) + \mu &= 0, \qquad  1 \leq s \leq g,   \label{s1} \\
\sum_k  x_k  &=1.  \label{s2}
\end{align}

This provides $(g+1)$ equations for $(g+1)$ unknowns $(x,\mu)$. We now show that for $m=0,$ this system has a unique explicit solution $(\wh{x}, \nu)$. By Theorem~\ref{thCHG}, for $m < m_0$, with $m_0 > 0$ being small enough, the function $u(x,y,z)$ is $C^{\infty}$ in $(m, x,y,z)$, with the explicit first-order Taylor expansion in $m$ given by \eqref{formulaCHG}. Recall our earlier notation,
\begin{equation*}
E_{s,k}(x,y) = \exp \left( -\frac{y_s}{F_s x_s} + \frac{y_k}{F_k x_k} \right),  \qquad
E_{s,k}(y,z) = \exp \left( -\frac{z_s}{F_s y_s} + \frac{z_k}{F_k y_k} \right). 
\end{equation*}

\hl{Taking derivatives} in $y$ of the first-order expansions \eqref{formulaCHG} of $C(x,y)$ and $C(y,z)$ readily yields the following first-order Taylor expansions in $m$, with remainders of order $m^2 \leq 10^{-12}$:  
\begin{align} 
\frac{\partial}{\partial y_s} C(x,y) &\simeq A_s (x,y) + m \wh{A}_s(x,y), \label{DyCxy1}\\ 
A_s(x,y) &= 1 + \log\frac{y_s}{F_s} - \log x_s, \label{DyCxy2}\\
\wh{A}_s(x,y) &= \sum_k \left( Q_{s,k} E_{s,k}(x,y) - \frac{F_k x_k}{F_s x_s} Q_{k,s}E_{k,s} \right), \label{DyCxy3}\\
\frac{\partial}{\partial y_s} C(y,z) &\simeq D_s(y,z) + m \wh{D}_s(y,z),  \label{DyCyx1} \\
D_s(y,z)  &= \frac{F_s}{\langle F,y \rangle} - \frac{z_s}{y_s}, \label{DyCyx2} \\
\wh{D}_s(y,z) &= F_s \sum_k Q_{s,k} - \left( F_s + \frac{z_s}{y_s} \right) \sum_k  E_{s,k}(y,z) Q_{s,k} - \frac{z_s}{F_s y_s^2} \sum_k  F_k y_k Q_{k,s} E_{k,s}(y,z). \label{DyCyx3}
\end{align}

\hl{For} $m= 0$,  the system in \eqref{s1} becomes $A_s(x,y) +D_s(y,z) + \mu =0$ for each $s$, which 
yields
\begin{equation*}
1 + \log\frac{y_s}{F_s x_s} + \frac{F_s}{\langle F, y \rangle} - \frac{z_s}{y_s} + \mu = 0.
\end{equation*}
\hl{Hence, we} have $x = e ^{1+\mu} X$, where the vector  $X$ is given by 
\begin{equation} \label{X}
X_s = \frac{y_s}{F_s}  \exp \left( \frac{F_s}{\langle F, y \rangle} - \frac{z_s}{y_s} \right) >0.
\end{equation}

\hl{The constraint} given by \eqref{s2} gives $e^{1+\mu} = \frac{1}{\sum_t  X_t}$. Therefore, for $m=0$ and all $s,$ the unique solution $\wh{x} = x(0,y,z)$ and $\wh{\mu}= \mu(0,y,z)$ of the system given by \eqref{s1} and  \eqref{s2} is
\begin{equation}\label{ZeroOrder}
\wh{x}_s = \frac{X_s}{\sum_t  X_t}, \qquad
X_s = \frac{y_s}{F_s}  \exp \left( \frac{F_s}{\langle F, y \rangle} - \frac{z_s}{y_s} \right), \qquad 
1+ \wh{\mu} = - \log\left(\sum_t X_t\right).
\end{equation}

\hl{Denote} $U(s,t) = \frac{\partial}{\partial x_t}  \frac{\partial}{\partial y_s}  u(x,y,z) + \mu$. For $m = 0,$ one has 
$$
U(s,t) = \frac{\partial}{\partial x_t}  \frac{\partial}{\partial y_s} [ A_s(x,y) +D_s(y,z) ] =  - 1_{\vset{s = t}}  \frac{1}{\wh{x}_s}.
$$ 
\hl{The} $g \times g$ matrix $U$ is thus diagonal with non-zero diagonal terms, and therefore  $U$ is invertible. Therefore, the implicit function theorem applies to the system given by \eqref{s1} and \eqref{s2}. Hence, for some fixed $m_1 >0$, there is a unique solution $(x(m, y,z), \mu(m,y,z))$ to the system \eqref{s1}--\eqref{s2}, and the functions  $(x,\mu)$ are $C^{\infty}$ in $(m, y,z)$.

Let $x \simeq \wh{x} + m V$ and $\mu \simeq \wh{\mu} + m \mu_1$ be the first-order Taylor expansions of $(x, \mu)$ in $m$.  Denote $V_s = \wh{x}_s w_s$ so that $x_s \simeq \wh{x}_s (1 + m w_s)$. The constraint \eqref{s2} then implies $\langle \wh{x}, w \rangle = 0$.  The first-order expansion of \eqref{s1} becomes
\begin{equation}\label{s3}
A_s(\wh{x} + mV , y) + D_s(y,z) + \wh{\mu} + m (\wh{A}_s(\wh{x} , y) +  \wh{D}_s(y,z) + \mu_1) \simeq 0.
\end{equation}
Since $A_s(\wh{x} + m V , y) \simeq A_s(\wh{x},y) - m w_s$, the zero-order term in \eqref{s3} vanishes due to the values of $\wh{x}$ and $\wh{\mu}$. The first-order term must vanish as well, which gives for all $s$
$$
w_s = \wh{A}_s(\wh{x} , y) +  \wh{D}_s(y,z) + \mu_1.
$$
\hl{Since} $ \langle \wh{x}, w \rangle = 0$, this yields 
$$
\mu_1 = - \sum_t \; \wh{x}_t [ \wh{A}_t(\wh{x} , y) +  \wh{D}_t(y,z)]
$$
\hl{and we then have}
$$
w_s = \wh{A}_s(\wh{x} , y) +  \wh{D}_s(y,z) - \sum_t   \wh{x}_t [\wh{A}_t(\wh{x} , y) +  \wh{D}_t(y,z) ]. 
$$

\hl{Define vectors} $\alpha$ and $\beta$  by $\alpha_s = \wh{A}_s(\wh{x},y)$ and $\beta_s = \wh{D}_s(y,z)$ so that
\begin{equation} \label{w}
w = \alpha + \beta - \langle \wh{x} , \alpha +\beta \rangle.
\end{equation}
\hl{Then the expression}s for $\wh{A}_s, \wh{D}_s,$ and $E_{s,k}$ given above directly yield  
\begin{equation} \label{alphas}
\alpha_s =  \sum_k  \left( Q_{s,k} e_{s,k} - \frac{F_k X_k}{F_s X_s} Q_{k,s} e_{k,s} \right)
\end{equation}
with $ e_{s,k}=1/e_{k,s} = \exp \left[ -y_s/(F_s \wh{x}_s) + y_k/(F_k \wh{x}_k)\right]$ and
\begin{equation} \label{betas}
\beta_s = F_s \sum_k Q_{s,k} - \left( F_s + \frac{z_s}{y_s} \right) \sum_k  f_{s,k} Q_{s,k} - \frac{z_s}{F_s y_s^2} \sum_k  F_k y_k Q_{k,s} f_{k,s}
\end{equation}
with $f_{s,k}=1/f_{k,s} = \exp \left[-z_s/(F_s y_s) + z_k/(F_k y_k)\right]$.
Note that $\wh{x}, \alpha, \beta,$ and $w$ depend only on $y$ and $z$. The preceding formulas provides the explicit first-order expansion $x_s = \wh{x}_s + m \wh{x}_s w_s$, concluding the proof.
\end{proof}

The reverse recurrence relation in Theorem~\ref{MinPathInt} introduces several challenges regarding numerical implementation. Once a target histogram $h^*_T = G$ is fixed, each geodesic is completely determined by its penultimate histogram $h^*_{T-1} =Z$. However, the penultimate histogram $Z$ is not known beforehand. Moreover, in concrete applications, the time horizon $T$ is often not known precisely. Performing a straightforward exhaustive numerical search for the optimal penultimate histogram, even for a very small number $g$ of genotypes, can be inefficient. We performed a preliminary numerical investigation of geodesics computations for a small number of genotypes $g = 3$ in~\cite{geiger_thesis}. Our empirical preliminary exploration indicates that geodesics connecting two interior histograms $H$ and $G$ seem to contain only interior histograms and are interior geodesics. In addition, we discuss the numerical implementation of the optimization algorithm in more detail in Ph.D. thesis~\cite{su_thesis} and a companion paper~\cite{azencott2023rare},
where we extensively explore and test efficient strategies to numerically compute geodesics for a larger number of genotypes, up to $g=10$.

Another potential research direction is to extend the mathematical formalism and the numerical computation of geodesics starting at a boundary histogram
$H$. For a boundary histogram, the number of $j$-cells is zero for some $j$ (same as $H(j) = 0$). Our analysis can be extended to this case, but this would require substantial research effort.

%%%%%%%%%%%%%%%%%%%%%%%%%%%%%%%%%%%%
\section{Conclusions}

We developed a large-deviation analysis of rare genetic events for \revision{interior histograms} in locked-box-type stochastic models of the genetic evolution of large bacterial populations. We consider discrete time-homogeneous Markov chains that model daily cycles of growth, mutations, and selection, where each daily selection randomly reduces population size to a fixed large size $N >10^5$. The daily population $pop_n$ then has fixed size $N$ and is characterized by its histogram $H_n \in \mcal{H} \subset R^g$, where  $H_n(j)$ is the frequency of genotype  $j$ in $pop_n$. The state space of the Markov chain $H_n$ is the compact convex set $\mcal{H}$ of all histograms in $R^g$. \revision{In addition,} mutations occur at a fixed small rate $m < 10^{-6}$. 
\revision{Restricting our analysis to interior histograms allows us to develop explicit formulas for the cost function \eqref{formulaCHG} and the backward optimization algorithm (Theorem~\ref{MinPathInt}), accurate up to order $m^2$. Since mutations occur at the small rate $m \le 10^{-6}$, these formulas are sufficiently accurate for many practical situations. Moreover, the validity of our asymptotic expansion for the one-step cost function was verified in~\cite{su_thesis} using the importance sampling algorithm.}
Other parameters (e.g., the mutation matrix $q_{i,j}$) in the Markov chain model can be fitted to emulate typical laboratory experiments on bacterial genetic evolution.

For the one-step Markov transition kernel $\mcal{Q}(H,G)$, we first proved the key uniform large-deviation estimate $\frac{1}{N} \log [ \mcal{Q}(H,G) ] \approx  - C(H,G)$ valid for large $N$ and arbitrary interior histograms $H,G$, where the one-step cost function $C(H,G) \ge 0$ is computed explicitly up to terms of order $m^2$.
This led us to obtain a uniform large-deviation estimate for the random histogram trajectories $\mbold{H} = \left[ H_1, \dots, H_T \right]$ of the Markov chain. For a fixed trajectory  $\mbold{h} = \left[ h_1, \dots, h_T \right]$, we demonstrate that  $Prob\{||\mbold{H} -\mbold{h}||< 1/N \}$ is roughly 
$\approx \exp(-N\,\lambda(\mbold{h}) )$, with a rate function $\lambda(\mbold{h}) $ explicitly given by the sum $C(h_1,h_2) + \dots + C(h_{T-1},h_T)$.
For any set $\Gamma$ of histogram trajectories of fixed duration $T$, let $\Gamma(N)$ be the $\Gamma$-neighborhood of radius $1/N$. We demonstrate that $Prob\{ \mbold{H}  \in  \Gamma(N) \}$ is of order $\exp(- N\, \Lambda(\Gamma) )$, with the rate functional defined as $\Lambda(\Gamma) = \inf_{\mbold{h} \in \Gamma} \, \lambda(\mbold{h})$.

For a fixed $T$, we applied our large-deviation theoretical results to develop a computational optimization algorithm for numerically computing the most likely histogram trajectory $\mbold{h}^*$ connecting a given initial histogram $h^*_1 =H\in \mcal{H}^\circ$ to a given terminal histogram $h^*_T =G\in \mcal{H}^\circ$. This problem is of particular interest when $G$ is not completely concentrated on the fittest genotype $g$, since in that case the conditional probability $P(H_T=G\; | \; H_1= H)$ vanishes at an exponential speed when $N$ becomes large.
Indeed, when $\Gamma(H,G)$ is the set of all histogram paths $\mbold{h}$ such that $h_1 =H$ and $h_T= G$, and if  $\Lambda(\Gamma(H,G)) \ne 0$, the random event $\mcal{E}= \{ \mbold{H} \in \Gamma(H,G) \}$ must have an exponential vanishing probability 
for large $N$. Our results then demonstrate that whenever $\mcal{E}$ actually occurs, any random histogram trajectory $\mbold{H}$ realizing $\mcal{E}$ must lie in a very thin tube around the path $\mbold{h}^*$ minimizing $\lambda(\mbold{h})$ over all $\mbold{h}\in \Gamma(H,G)$. In laboratory experiments, knowing only $H_1 =H$ and the observed terminal histogram $H_T= G$, there is a practical interest to reconstructing the whole unobserved random genetic evolution  $\mbold{H} \approx \mbold{h}^*$ between times $1$ and $T$.

In this paper
we developed (Theorem~\ref{MinPathInt}) a concrete theoretical and numerical approach for computing a minimizer $\mbold{h}^*$ of $\lambda(\mbold{h})$ over all $\mbold{h} \in \Gamma(H,G)$. By analogy with the geodesics minimizing Riemannian energy on a Riemannian manifold, we call any such minimizer $\mbold{h}^*$ a large-deviation ``geodesic'' connecting $H$ to $G$ in $T$ steps. We demonstrate that any geodesic $\mbold{h}^*$ is completely determined by its last two histograms, the known $h^*_T= G$ and the unknown ``penultimate histogram'' $h_{T-1} =Z$. This is due to the recursive relation in reverse time, which expresses  $h^*_{n+2}$ as a function $ \chi(m,h^*_{n+1}, h^*_{n})$. 
We explicitly compute the histogram-valued function $\chi$ for a small mutation rate $m$. This reduces the search for $\mbold{h}^*$ (given $H$ and $G$) to a search for the optimal penultimate histogram $Z$. In our companion paper~\cite{azencott2023rare}, we explored this strategy to develop an efficient numerical algorithm for numerically computing  the large-deviation geodesics for a relatively large number of genotypes, up to $g \approx 10$.

\revision{Although most of our analysis is restricted to interior histograms, it applies to histograms in which all genotype frequencies are strictly positive but may be very small. Therefore, our explicit formulas can provide valuable insight into the mechanisms describing the emergence of non-dominant genotypes in large populations. 
In particular, one can study the influence of deleterious and advantageous mutations by appropriately modifying the mutation matrix $Q$.
In addition, we present preliminary results for the analysis of the cost function for the emergence (from zero) of just one genotype in Section~\ref{sec:boundary}.}
Thus, our \revision{analytical} results open the way to interesting further applied work on the biological mechanisms that can trigger concrete rare genetic events, such as the emergence or fixation of non-dominant genotypes. Our explicit cost functions formulas can help quantify how the occurrence of specific rare genetic events is impacted by changes in the mutation rate $m$, or by the emergence of deleterious mutations.

\revision{Although most of our analysis is restricted to interior histograms, it applies to histograms in which all genotype frequencies are strictly positive but may be very small. Therefore, our explicit formulas can provide valuable insight into the mechanisms underlying the appearance of non-dominant genotypes in large populations. In particular, one can study the influence of deleterious and advantageous mutations by appropriately modifying the mutation matrix $Q$.
In addition, in Section~\ref{sec:boundary}, we present preliminary results on the cost function for the emergence of a single genotype from zero frequency.

Thus, our analytical results open the way to further applied work on the biological mechanisms that can trigger concrete rare genetic events, such as the emergence or fixation of non-dominant genotypes. Our explicit formulas for the cost function can help quantify how the occurrence of specific rare genetic events is affected by changes in the mutation rate $m$ or by the presence of deleterious mutations. These formulas also provide a variational framework for comparing competing rare-event scenarios by identifying the most likely paths leading to a prescribed genetic outcome.}
\vspace{6pt}

\subsection*{Acknowledgments}
This research was partially supported by the U.S. National Science Foundation grant DMS-1412927. 
I.T. was also partially supported by the U.S. National Science Foundation grant DMS-1903270.

The authors thank Profs. R. Azevedo and T. Cooper for helpful discussions on this subject.

%%%%%%%%%%%%%%%%%%%%%%%%%%%%%%%%%%%%%%%%%%%%%

\appendix

\section{Summary of Notation} Below, we summarize the notation used in this manuscript.
\label{sec:ap1}

\begin{itemize}
\itemsep0em
\item
$g \in \bN$: Number of genotypes in the bacterial population
\item
$N \in \bN$: Number of cells in the bacterial population
\item
$[F_1, F_2,\cdots, F_g]$: Ordered genotype growth factors,
with $F_1 <F_2 <\cdots <F_g$
\item
$m \in \bR$: Mutation rate
\item $Q \in \bR^g \times \bR^g$, with entries $q_{j,k} \ge 0$-conditional probabilities of mutation
(see Definition~\ref{Qjk})
\item $M=mQ$: Mutation matrix (see Definition~\ref{mutM}) 
\item Model parameters $\calP=\{N, g, F_1,\cdots, F_g, m, Q\}$ (see Definition~\ref{P})
\item $R_n(j,k)$: Number of $j$-cells mutating into $k$-cells on day $n$
(see Definition~\ref{Rn}) 
\item
$H \in \bR^g$: Population histogram of bacterial frequencies frequencies; $H = [H(1), H(2),$ $\ldots, H(g)]$, with $H(j) \ge 0$
for $j=1,\ldots,g$ and
$\sum_{j=1}^g H(j) = 1$ 
 \item
$H(j) \in [0,1]$: Frequency of genotype $j$ in the population
\item
$H_n \in \bR^g$: Population histogram of bacterial frequencies on day $n$
\item
$H_n(j) \in [0,1]$: Frequency of genotype $j$ in the population on day $n$
\item
$\mathcal{H} \subset \bR^g$: Space of all histograms
(see Definition~\ref{calH})
\item
{ $\mathcal{H}^\circ =\{H: H \in \mathcal{H} \text{ and } H(j) > 0 \text{ for all } j=1,\ldots,g\}$ is the 
set of all interior histograms}
\item
$\mathbf{H}= \{H_1,\cdots, H_T\}$: Time-dependent 
path of length $T$ (days) in the space $\mathcal{H}$
connecting the initial histogram $H_1$ and the final histogram $H_T$
\item
$\Omega_T$: Space of all paths of length $T$ (see Definition~\ref{def:OmegaT})
\item
$\Omega_T^+$: Space of all interior paths of length $T$; 
if $\mathbf{H} \in \Omega_T^+$, then $H_n(j) \ge \eps$
for all $1 \le t \le T$ and $1 \le j \le g$ and for some
$\eps \ll 1$
\item 
{$\supp (H)$: Support of a histogram 
$\supp (H) = \vset{j \, | \,H(j)>0}$ for any $H \in \mathcal{H}$
(see Definition~\ref{def:suppmin})}
\item 
$b(H)$: The essential minimum for a histogram
$b(H) = \min_{j \in \supp (H)} \; H(j)$ (see Definition~\ref{def:suppmin})
\item $b(\mbold{H})$: The essential minimum for a path of histograms $b(\mbold{H}) = \min_{n =1\cdots T}  b(H_n)$
(see Definition \ref{def:bboldH})
\item$ \ceil{u}$ is
the smallest integer greater than or equal to $u$ for $u\geq 0$
\item 
$\Phi_j(H) = F_j H(j) / \langle F, H \rangle$: Approximation of the genetic histogram after the growth phase
(see Equation~\eqref{phi})
\item 
$\Psi_j(H, r) = \left(F_j H(j) - \sum_k  r_{j,k} + \sum_k  r_{k,j}\right) / \langle F, H \rangle$: Approximation of the genetic histogram after the mutation phase
(see Equation~\eqref{psi})
\item
A matrix $A$ is \emph{N-rational} if $N A$ has non-negative integer coefficients 
\item 
The set $\calZ$ contains all $N$-rational $g \times g$ matrices
\item 
$K(j, H)$ is the set of $g \times g$ matrices with non-negative coefficients with constraints \eqref{Kj} 
(see Definition \ref{KNconstr})
\item 
$K(H)  = \bigcap_{j=1}^g \; K(j, H)$ and $K_N(H) = \calZ \cap K(H)$
\item 
$N r_{j,r} = R_{j,k}$: Total number of $j$-cells mutating into $k$-cells
\item $\mathcal{Q}_N (H, G) = P(H_{n+1} = G \; | \; H_n = H)$: Markov transition kernel (see Definition~\ref{kernel})
\item 
$\tau (H, r, G) = mut(r, H) + KL(G, \Psi(H, r))$: Composite transition rate
\item $KL(G,J)$: Kullback--Leibler divergence (see Definition~\ref{def:kl})
\item $mut(H,r)$: Rate function for mutations (see Definition~\ref{def:mut})
\item $C(H,G) = \min_{r \in K(H)} \tau(H,r,G)$: One-step cost function 
(see Definition~\ref{def:onestep})
\item
$\Omega_T (a) \subset \Omega_T$: The compact set of histogram trajectories such that $\mbold{H} \in \Omega_T$ and $b(\mbold{H}) \geq a$
\item The distance between trajectories: $\vnorm{ \mbold{H} - \mbold{H}' }= \max_{n =1\cdots T}  \vnorm{ H_n - H'_n }$
\item $\lambda(\mbold{H})$: Large-deviation rate function for the path 
$\mbold{H}$; 
$\lambda(\mbold{H}) = \sum_{n=1}^{T-1} C(H_n, H_{n+1})$
\mbox{(Equation~\eqref{lambda})}
\item
$\Lambda(F)$: Large-deviation set functional; 
$\Lambda(F) = \inf_{ \mbold{H} \in F} \; \lambda(\mbold{H})$ for any 
$F \subset \Omega_T$
\mbox{(Equation~\eqref{Lambda})}
\item $V_N (\Gamma)$: Open neighborhood of paths for any $\Gamma \subset \Omega_T$ (Definition~\ref{VNBNpaths})
\item $b(\Gamma)$: The essential minimum for a set of paths; 
$b(\Gamma) = \inf_{\mbold{H} \in \Gamma} b(\mbold{H})$ for any 
$\Gamma \subset \Omega_T$ (Definition~\ref{VNBNpaths})
\item $B_N(\Gamma)$: The (finite) set of all $N$-rational paths in $V_N(\Gamma)$
\end{itemize}

%%%%%%%%%%%%%%%%%%%%%%%%%%%%%
\section{Proofs of Theorems in Section~\ref{stochmodel}}
\label{sec:ap2}

\subsection*{Proof of Lemma \ref{density}}

In this proof only, $\floor{z}$ denotes the largest integer less than or equal to $z$.  Fix any $j \in \supp(H)$, and select any $k^*= k^*(j)$ such that $r_{j, k^*} = \max_ k \; r_{j,k}.$  For any $k$, if  $r_{j,k*} < g/N$, define $s_{j,k}$ by
\begin{equation} \label{defsjk}
s_{j,k} =
\begin{cases}
0, & r_{j,k} = 0; \\
1 / N, & 0 < r_{j,k} < 1/N; \\
 \floor{N r_{j,k}} / N, & 1/N \leq r_{j,k}.
\end{cases} 
\end{equation}

\hl{If } $g/N \leq r_{j,k*}$, define 
\begin{equation} \label{defsjk*}
s_{j,k*} = [ N \; r_{j,k*} ] / N - (g-1 )/N.
\end{equation}

\hl{Equations} \eqref{defsjk} and \eqref{defsjk*} imply that 
$\vnorm{s - r} \leq g/N$ and $\supp(s) = \supp(r)$.
Let $S(j) = \sum_k \, s_{j,k} $. To prove $s \in K(j, H)$, we only need to show that $S(j) < F_j H(j)$ for all $j \in \supp(H)$, which involves two cases.  Thus, let $j \in \supp(H).$\smallskip\\
{\bf \hl{Case 1:}
} Suppose $r_{j,k*} < g/N.$  Since $F_j H(j) > F_1 b(H) > F_1 a$, impose $N > g^2/a F_1$ to force $F_j H(j) > g^2/N $. Definition \eqref{defsjk} ensures that $ s_{j,k} < g/N$ for all $j, k$, and hence $S(j) < g^2/N < F_j H(j)$.\smallskip\\
{\bf \hl{Case 2:}} Suppose $g/N \leq r_{j,k*}$.  Let $u(j) \leq (g-1)$ be the number of indices $k$ such that $0 < r_{j,k} < 1/N$. Definitions \eqref{defsjk*} and \eqref{defsjk} imply $s_{j,k*}\leq r_{j,k*} - (g-1)/N$ and $\sum_{k \neq k*}  s_{j,k} \leq u(j)/N + \sum_{k \neq k*}  r_{j,k}.$  This yields $S(j) \leq \sum_k  r_{j,k} < F_j H(j)$ since $r \in K(H)$.\\
This concludes the proof.

%%%%%%%%%%%%%%%%%%%%%%%%%%%%%
\section{Proofs of Theorems in Section~\ref{LDcycles}\label{sec:ap3}}

First, we need the following technical lemma.
\begin{Lemma}
    For any fixed $0 < \alpha <1$ and for all $x,y\in [0,1/2],$ the following inequality holds: 
\begin{equation} \label{xlogxholder}
\abs{y \log y - x \log x } \leq \frac{4}{1-\alpha} \abs{x - y}^{\alpha}.
\end{equation}

\hl{Moreover, for} $x,y\in [0,A]$ with $A > 1$ the following inequality holds:
\begin{equation}\label{basic.holder}
\abs{x \log x  - y \log y } \leq \left(\frac{5 (1 + A ) \abs{\log(A)} }{1-\alpha}\right) \abs{x - y}^\alpha .
\end{equation}
\end{Lemma}
%%%%%%%%%%%%%%%%
\begin{proof}
Consider $0 \leq x \leq 1$. Then,
\begin{equation}\label{xlogx}
| x \log(x) | = x \log (1/x) \leq (x/(1-\alpha)) \log (1/x^{1-\alpha}) \leq x^{\alpha} / (1-\alpha).
\end{equation}
Also \hl{consider} $\abs{y \log y - x \log x }$ with $ \abs{y-x} < w \leq 1$ . Then, either $x, y \leq 2 w$ or $x,y \geq w$. When $x, y \leq 2 w$, Equation~\eqref{xlogx} gives 
$$
\abs{y \log y - x \log x } \leq \abs{x \log x} + \abs{y \log y} \leq \frac{4}{(1-\alpha) w^{\alpha}}. 
$$
\hl{When both} $x,y \geq w$, Taylor's formula and \eqref{xlogx} yield 
$$
\abs{y \log y - x \log x } \leq w ( 1 + \log(1/w) ) \leq \frac{4}{1-\alpha} w^{\alpha}. 
$$

% INDENT
\hl{This proves} Equation~\eqref{xlogxholder}.
For $x,y \in [0,A]$ with $A >1$, apply Equation~\eqref{xlogxholder} to $x/A$ and $y/A$ to obtain
Equation~\eqref{basic.holder}.
\end{proof}

%%%%%%%%%%%%%%%%%%%%%%%%%%%%
\subsection{Proof of Proposition~\ref{propholdermut}}
Consider $a$ and $\alpha$
in Proposition~\ref{propholdermut}.
Take $H, H' \in \mcal{H}(a)$ with $\supp(H') = \supp(H)$, $r \in K(H),$  $r' \in K(H')$. Fix $(j,k)$ with $H(j) q_{j,k} > 0$, which implies $H'(j) q_{j,k}> 0$. Then, $ x := r'_{j,k} $ and $y := r_{j,k}$ verify $\abs{ x- y } \leq \vnorm{ r' - r }$ and are bounded by $\max\vset{ F_j H'(j),  F_j H(j)} \leq F_g $.  By definition, 
\begin{align*}
L_{j,k}(H,r) &= y \log y - y \log(e F_j H(j)) + F_j H(j), \\
L_{j,k}(H',r') &= x \log x - x \log( e F_j H'(j)) + F_j H'(j),
\end{align*}
and \hl{set }$\hat{L} := \abs{ L_{j,k}(H,r) - L_{j,k}(H',r') },$ which satisfies the inequality 
\begin{equation} \label{dL}
\begin{aligned}
\hat{L} \leq & \abs{ y \log y  - x \log x } + \abs{y-x}  \abs{ \log(e F_j H(j)) } + \\
& \abs{x}  \abs{ \log ( H'(j)/H(j) ) } + F_j  \abs{ H'(j)- H(j) }. 
\end{aligned}
\end{equation}

Next,  \hl{reformulate} \eqref{dL} by writing $\hat{L} \leq U_1 + U_2 + U_3 + U_4$ 
where $U_i$ is the $i$-th term in the right-hand side of \eqref{dL}.
Clearly $U_3 + U_4 \leq 2 (F_g / a) \vnorm{ H' - H }$.  Then, \eqref{basic.holder} gives $U_1 \leq \frac{c_1}{1-\alpha}\vnorm{ r' - r }^\alpha$, with $c_1 =5 (1+F_g) \log(F_g) $. Since $F_g \geq F_j H(j) \geq a F_1$, one has $\abs{\log(e F_j H(j) )} \leq c_2 /a$, with $c_2 = 2 + \log(F_g/F_1)$, and hence, $ U_2 \leq c_2 \vnorm{r'- r}$. These bounds yield, for $H(j) q_{j,k} > 0,$
\begin{equation}\label{hatLjk}
\hat{L} \leq  \frac{\hat{c}}{a (1-\alpha) }\left( \vnorm{ r' - r }^{\alpha} + \vnorm{ H' - H }\right) 
\end{equation}
with $\hat{c} = c_1+ c_2 + 2 F_g \leq 14 F_g \log(F_g) $. This result still holds when $H(j) q_{j,k} = 0$ since then $\hat{L} = 0$. Summing \eqref{hatLjk} over all $j,k$, we get 
$$
\begin{aligned}
\abs{ mut(H',r') - mut(H,r) } & \leq \frac{\hat{c} g^2 }{a (1-\alpha)}\left(  \vnorm{ r' - r }^{\alpha} + \vnorm{ H' - H}\right)  \\
& \leq c\left(  \vnorm{ r' - r }^{\alpha} + \vnorm{ H' - H}\right),
\end{aligned}
$$
where $c$ is the constant stated in this proposition. This concludes the proof.

%%%%%%%%%%%%%%%%%%%%%%%%%%%%
\subsection{Proof of Proposition~\ref{LDmut}}
Take $H$ and $r$ as stated in Proposition~\ref{LDmut}. 
Given $H_n = H $, the coefficients $Z_n(j,k)$ of the companion matrix $Z_n$ are independent and have Poisson distributions $\pi(H,j,k)$ with respective means $N m\,q_{j,k} F_j H(j)$. For $j \in \supp(H)$ and any $k$, set $u = m\,q_{j,k} F_j H(j)$ and $v = r_{j,k}$. Since $N v$ is an integer, apply \eqref{XNv} to $X=Z_n(j,k)$ to obtain 
\begin{equation}\label{L0}
\frac{1}{N}\log P(Z_n(j,k) / N  =  r_{j,k} \; | \; H_n = H) = - L_{j,k}(H,r) + o_1(1)
\end{equation}
 with $\abs{ o_1(1) } \leq 2 \log N /N.$ This equation remains true for $j \notin \supp(H)$ and all $k$, since if $j \notin \supp(H)$, then $H(j) = 0$ and $Z_n(j,k) = 0$.

Conditional independence of the $Z_n$ coefficients then yields 
$$ 
P( Z_n/N = r \: | \; H_n = H)  =  \prod_{j , k} \; \; P ( Z_n(j,k)/N = r_{j,k} \; | \; H_n = H). 
$$

\hl{Equation} \eqref{L0} implies
\begin{equation}\label{L1}
\frac{1}{N} \log P( Z_n / N = r \; | \; H_n = H )  =  - \sum_{j, k} L_{j,k}(H,r) + o_2(1)  =  - mut(H,r) + o_2(1)
\end{equation}
with $\abs{ o_2(1) }  \leq g^2 \abs{ o_1(N) } \leq 2 g^2  \log N /N$.

Set $c = 2 + \log(g)$, and impose $N >c/a \geq \frac{c}{b(H)}$. Applying \eqref{control} yields
\begin{equation} \label{L3}
\begin{aligned}
& \left|\frac{1}{N} \log P (  R_n/N = r \; | \; H_n = H ) - \frac{1}{N}\log P ( Z_n/N = r \; | \; H_n = H )\right| \\
& \leq  \frac{1}{N} \log\left( 1 + 2 g \, d(H)^{N/2} \right). 
\end{aligned}
\end{equation}
\hl{The right-hand} side of \eqref{L3} is bounded above by $\frac{1}{N} \log(1 + 2 g)$, so that Equation~\eqref{L1} implies
\begin{equation}\label{L4}
\frac{1}{N} \log P ( R_n/N = r \; | \; H_n = H )  =  - mut(r, H) + o(1) 
\end{equation}
with $\abs{ o(1) }  \leq  \frac{1}{N} \left( \log(1 + 2 g) + 2 g^2 \log N  \right) \leq 4 g^2 \log(N)/N $. \\
This concludes the proof.

%%%%%%%%%%%%%%%%%%%%%%%%%%%%
\subsection{Proof of Proposition~\ref{LDmultinomial}}

The coordinates of $V= N G$ are non-negative integers that sum up to $N$. 
Therefore, \eqref{multinomial2} gives 
\begin{equation}\label{logmu}
\frac{1}{N} \log(\mu_{N,J}(N G)) = \frac{1}{N} \log N! - \sum_{j \in \supp(G)} \frac{1}{N} \log V(j)! + \sum_{j \in \supp(G)}  \frac{V(j)}{N} \log J(j). 
\end{equation}

\hl{For} $j \in \supp(G) \subset \supp(J)$, we apply Stirling's formula \eqref{stirling} to $V(j)!= [ N G(j)]!$ and obtain 
\begin{equation}\label{logV}
\frac{1}{N}\log V(j)! = G(j)\log N + G(j) \log G(j) - G(j) + o_j(1)
\end{equation}
with $\abs{ o_j(1) } \leq 2\abs{\log(N G(j)) / N} $.  Since $G$ is N-rational, then for each $j \in \supp(G)$ one has $ N G(j)>0$, and hence  $1 \leq N G(j) \leq N$. This yields $0 \leq \log( N G(j) ) \leq \log(N)$ and $\abs{o_j(1)}  \leq  2  \log N /N$.  
Since  $\supp(G) \subset \supp(J)$, we apply \eqref{stirling} to $\log(N!)$ along with \eqref{logV} and Equation~\eqref{logmu} to get
\begin{equation} \label{KL1} 
\frac{1}{N} \log(\mu_{N,J}(N G)) = \sum_{j \in \supp(G)} \left[- G(j) \log {G(j)} + G(j) \log J(j) \right] + o(1)
\end{equation}
with uniform remainder $\abs{o(1)} \leq 2 (g+1) \log N /N $. Notice that the sum in \eqref{KL1} is equal to the Kullback--Leibler divergence $- KL(G,J) < 0$  given by Equation~\eqref{KL}. \\
This concludes the proof.

%%%%%%%%%%%%%%%%%%%%%%%%%%%%
\subsection{Proof of Proposition~\ref{LDKL}}

Consider $G,G',$ as stated in Proposition~\ref{LDKL}. 
Then, $N > N^*$ forces $N > 2/b(G) $ so that $\supp(G') = \supp(G)$ by Lemma \ref{basicBN}.  
Let $J_n$ be the population histogram at the end of the mutation phase (Phase 2).
Suppose first that $KL(G,J_n)$ is finite so that $ \supp(G') =\supp(G) \subset \supp(J_n)$. 
By construction of the Markov chain $H_n$, one has
$
P ( H_{n+1} = G' \; | \; H_n ,  R_n )= P ( H_{n+1} = G' \; | \; J_n ).
$
Given $J_n$, the conditional distribution of $N H_{n+1}$ is the multinomial $ \mu_{N, J_n}$ given by \eqref{multinomial}. Since $G'$ is $N$-rational, $P (  H_{n+1} = G'  \; | \; J_n) = \mu_{N, J_n}(N G').$ Then, \eqref{KLlogmu} yields
\begin{equation}\label{h1}
\frac{1}{N} \log P(  H_{n+1} = G' \; | \; J_n ) = - KL(G', J_n) + o(1)
\end{equation}
with $\abs{ o(1) }  \leq  2 (g+1) \log N/N$.  From the KL derivative  \eqref{partialKL}, we also deduce
\begin{equation}\label{h2}
\abs{KL(G', J_n) - KL(G, J_n)}  \leq  \frac{2 g}{3 N} \left( 1 + \abs{ \log b(G) } + \abs{ \log b(J_n) }\right). 
\end{equation}

\hl{Condition} $b(J_n) \geq 1/(N F_g) $ holds due to \eqref{bJn}. 
Since $b(G) \geq 2/ N$, the right-hand side of \eqref{h2} is bounded above by $2 g \log N / N$ provided $N > N^*$. Equation~\eqref{h1} then implies 
\begin{equation}\label{h4}
\frac{1}{N} \log P(H_{n+1} = G' \; | \; J_n ) =  - KL(G, J_n) + o_1(1)
\end{equation}
with $\abs{ o_1(1) } \leq 5 g \log N /N$.  Hence, \eqref{LDselect1} is proved when $KL(G, J_n)$ is finite. \smallskip \\
When $KL(G, J_n) = + \infty$, we have $\supp(G)\not \in \supp(J_n)$. For $G' \in B_N(G)$, one has $\supp(G') = \supp(G) \not \in \supp(J_n)$ so that the transition from $J_n$ to $H_{n+1} = G'$ is impossible during Phase 3. Thus, both sides of \eqref{LDselect1} are equal to $- \infty$. \\
This concludes the proof.

%%%%%%%%%%%%%%%%%%%%%%%%%%%%
\subsection{Proof of Lemma \ref{lipKL}}

Define the entropy $\mscr{E}(G)$ as $0 \leq \mscr{E}(G) = - \sum_{j \in \supp(G)} G(j) \log(G(j)) \leq \log(g)$.  By definition of $KL(G,J)$, one has 
\begin{equation} \label{beta2}
KL(G, J) + \mscr{E}(G) = \sum_{j \in \supp(G)}  G(j) \log(1/J(j)).
\end{equation}

\hl{For each} $k \in \supp(G)$, this yields $b(G) \log{(1/J(k))} \leq KL(G, J) + \mscr{E}(G)$. Hence, the term $\beta(G,J)$ in \eqref{betabound}  verifies
$\log{\beta(G,J)} \leq [ KL(G, J) + \log(g) ] / b(G)$.
Since $G$ is a histogram, Equation~\eqref{beta2} implies 
$$ 
KL(G,J) \leq KL(G, J) + \mscr{E}(G) \leq \max_{j \in \supp(G)}  \log(1/J(j)) = \log(\beta(G,J)).
$$
\hl{This proves} Equation~\eqref{betabound}.
Now consider $G, J, J' \in \mcal{H}$ with $KL(G,J) \leq KL(G,J') < \infty$.  Formula \eqref{partialKL} implies $\abs{\partial_{J(k)} KL(G, J)} \leq \beta(G,J)$ for all $k \in \supp(G)$, with $\beta(G,J)$ as above. By Taylor's formula and \eqref{betabound}, we obtain
\begin{equation}
\abs{KL(G,J) - KL(G,J')}  \leq  g  \beta(G,J) \vnorm{ J - J' } \leq c  \vnorm{ J - J' }
\end{equation}
with $c = g  \exp ( [ K(G,J) + \log(g) ] / b(G) )$, proving Equation~\eqref{lipKLJI}.
Now take $a, G, G',$ and $J$ as stated in (ii). For $j \in \supp(G) \cup \supp(G')$, let $u(j) = [G'(j) - G(j) ]\log J(j)$ and $v(j) = G'(j) \log G'(j) - G(j) \log G(j).$ By definition of  Kullback--Leibler divergence, we have
\begin{equation} \label{sumKL}
 \abs{KL(G',J) - KL(G,J)} \leq \sum_{ j \in \supp(G) \cup \supp(G') } ( \abs{u(j)} + \abs{v(j)}).
\end{equation}

\hl{Equation} \eqref{xlogxholder} with $\alpha = 1/2$ implies
\begin{equation} \label{generic.v(j)}
| v(j) | \leq \; 8 \; ||\; G' - G \;|| ^{1/2}, \qquad j \in \supp(G) \cup \supp(G'). 
\end{equation}

\hl{For} $ j \in \supp(G) \cup \supp(G')$, one has $j\in \supp(J)$, and the
bound \eqref{betabound} gives 
\begin{align}
\abs{ u(j) } &\leq \vnorm{ G'- G }  
\max\left(\log  \beta(G,J)  , \log \beta(G,J)\right)  \nonumber \\ 
 &\leq \vnorm{ G'- G }\left( \frac{1}{a}\left[ \log g + 
 \max\left(KL(G,J), KL(G',J)\right)\right]\right),  \label{generic.u(j)}
\end{align}
and we can \hl{combine} \eqref{generic.u(j)}, \eqref{generic.v(j)}, and \eqref{sumKL} to get 
\begin{equation}
\abs{KL(G',J) - KL(G,J)} \leq c_1 \vnorm{ G' - G } ^{1/2} 
\end{equation} 
with $c_1 = \frac{g}{a} \left(8 + \log g + 
\max\left( KL(G,J) , KL(G',J)\right)\right)$, proving Equation~\eqref{holderKLG}.
For $G,G',J$ as above, consider the special case $\supp(G) =\supp(G') $. From \eqref{betabound}, one gets, for $j \in \supp(G) = \supp(G') $,
\begin{equation*} 
\abs{ u(j) } \leq \frac{1}{a} \vnorm{ G'- G }\left(\log g + \min \vset{KL(G,J), KL(G',J)}\right).
\end{equation*}
and Taylor's formula implies that
\begin{equation} 
\abs{v(j)} \leq \vnorm{ G' - G} \left[ 1 + \log (1/a)\right] .
\end{equation}
Therefore, \hl{these last two} bounds along with \eqref{sumKL} prove Equation~\eqref{lipKLG}. \\
This concludes the proof.

%%%%%%%%%%%%%%%%%%%%%%%%%%%%
\subsection{Proof of Proposition~\ref{holdertau}}

Consider $H',H,G, r',$ and $r$ verifying (i)--(iv) in Proposition~\ref{holdertau}.  Apply \eqref{holdermut} with H\"{o}lder coefficient $\alpha=1/2$ to obtain
\begin{equation}\label{deltamut}
\abs{mut(r', H') - mut(r, H) } \leq c_0  \left[  \vnorm{ r' - r }^{1/2} + \vnorm{ H'- H } \right ] ,
\end{equation}
where $c_0 = \frac{ 20}{a} F_g \log(F_g).$  From \eqref{lipPsi}, we also obtain
\begin{equation} \label{lipPsi2}
\vnorm{ \Psi(H',r') - \Psi(H,r) } \leq \frac{3 g F_g}{F_1} \left( \vnorm{ r' - r } + \vnorm{ H' - H } \right). 
\end{equation}

\hl{By definition} of $\tau(H,r,G)$ in \eqref{tauHrG} and item (iii) above, one has 
\begin{equation} \label{KL+2}
KL( G, \Psi(H,r) ) \leq \tau(H,r,G) \leq A,
\end{equation} 
and \hl{since} $\tau(H', r',G)$ and $\tau(H, r, G)$ are finite, both $\supp(\Psi(H,r))$ and $\supp(\Psi(H',r'))$ contain $\supp(G)$. Therefore, we can apply \eqref{KL+2} and \eqref{lipKLJI} to obtain 
$$
\abs{KL( G, \Psi(H',r')) - KL(G, \Psi(H,r))}  \leq  c_1  \vnorm{ \Psi(H',r') - \Psi(H, r) }
$$
with $c_1 = g^{1+1/a} e^{A/a}.$  In view of \eqref{lipPsi2}, this yields
$$
\abs{ KL( G, \Psi(H',r') ) - KL( G, \Psi(H,r) ) }  \leq  c_2  \left[  \vnorm{ r' - r } + \vnorm{ H'-H } \right]
$$
with $c_2 = 3 c_1 g F_g / F_1 $.  Combining this last result with \eqref{deltamut} and setting $\eta = c_0 + c_2$ yields
$$
\abs{\tau(H',r',G) - \tau(H,r,G) }\leq \eta \left[\vnorm{ r' - r }^{1/2} + \vnorm{ H'-H } \right].
$$
\hl{Noting that} $\eta$ here is precisely \eqref{eta} concludes the proof.

%%%%%%%%%%%%%%%%%%%%%%%%%%%%
\subsection{Proof of Theorem~\ref{holderCost}}

Consider $G, G', H\in \mcal{H}(a)$ verifying \eqref{conditions}. One can then select $r \in K(H)$ with $\tau(H,r,G) \leq d+1$, so that $I = \Psi(H,r)$ must verify $\supp(G') = \supp(G) \subset \supp(I)$ and $ KL(G, I) \leq d+1$ due to the definitions  of $C(H,G)$  and $\tau(H,r,G)$ by Equations~\eqref{costHG} and \eqref{tauHrG}. 
Apply \eqref{holderKLG} to get the Lipschitz bound: 
\begin{equation} \label{KLGG'}
\abs{KL(G',I) - KL(G,I)} \leq c_1  \vnorm{ G' - G }  
\end{equation} 
with $c_1 = \frac{g}{a} (3 + \log g  + d).$ This implies, by definition of $\tau(H,r,G)$,
$$
\abs{\tau(H,r,G') - \tau(H,r,G)} = \abs{KL(G', I) - KL(G, I)} \leq c_1 \vnorm{ G' - G }. 
$$
and then, since $ C(H,G') \leq \tau(H,r,G'),$
\begin{equation} \label{Ctau}
C(H,G') \leq \tau(H,r,G) + c_1 \vnorm{ G' - G } \leq d_1
\end{equation} 
with $d_1 = d + 1 + c_1$. In \eqref{Ctau}, take the infimum of the middle term over $r \in K(H)$ to obtain
\begin{equation} \label{CCGG'}
C(H,G') \leq C(H,G) + c_1  \vnorm{ G' - G }. 
\end{equation} 
\hl{Let} $c_2= \frac{g}{a} (3 + \log g + d_1)$. Since $ C(H,G') \leq d_1$, the generic result \eqref{CCGG'} can now also be rewritten by switching the roles of $G$ and $G'$ provided one also replaces $d$ with $d_1$ and $c_1$ with $c_2$.  This yields $C(H,G) \leq C(H,G') + c_2 \vnorm{ G' - G }$. Therefore, the bound in \eqref{CCGG'}~yields
\begin{equation} \label{lipCHGG'}
\abs{ C(H,G) - C(H,G') } \leq c_2  \vnorm{ G' - G }.
\end{equation}

\hl{Define a linear} mapping $\mcal{L}$ for all $\rho \in K(H)$, denoted $\hat{\rho} := \mcal{L}(\rho)$, given by
\begin{align} \label{biject}
\hat{\rho}_{j,k}  &= \rho_{j,k} H(j) / H'(j) \;\text{for $j \in \supp(H)$ and all $k$ }, \\
 \hat{\rho}_{j,k} &= 0 \;\;\; \text{for $j \notin \supp(H)$ and all $k$ }.
\end{align} 
Thus, by definition of $K(H')$ and $K(H)$, one readily verifies that $\hat{\rho} \in K(H)$. Therefore, $\mathcal{L}$ maps $K(H')$ into $K(H)$ bijectively with an inverse mapping $ \hat{\rho} \to \rho$ defined by exchanging $\rho$ and $\hat{\rho}$ as well as $H'$ and $H$ in \eqref{biject}.  Definition \eqref{biject} forces $\supp(\rho) = \supp(\hat{\rho})$. Since $\supp(H')= \supp(H)$, then \eqref{supppsiHh} yields $\supp(\Psi(H',\rho))= \supp(\Psi(H,\hat{\rho}))$. Hence, either these two supports contain $\supp(G)$ or neither of them do. This implies the equivalence
\begin{equation} \label{bothfinite}
\tau(H',\rho,G)< \infty \iff \tau(H,\hat{\rho},G) < \infty.
\end{equation}

\hl{For } $\rho \in K(H')$, the bound $\vnorm{ \rho } \leq F_g$ holds due to \eqref{Kj}.  Thus, we get $\abs{ 1- H(j)/H'(j) } \leq \vnorm{ H' - H } / a $ for $j \in \supp(H') = \supp(H)$. Hence, for all $\rho \in K(H'),$
\begin{equation} \label{rohatro}
\vnorm{ \rho - \hat{\rho} }  \leq \vnorm{ \rho}\, \vnorm{ H' - H } / a  \leq  \frac{F_g}{a}  \vnorm{ H' - H } \leq 1.
\end{equation}

\hl{Define} $U =\vset{ r \in K(H) \,|\, \tau(H,r,G) \leq d}$.  This implies $C(H,G) = \min_{r \in U} \; \tau(H,r,G) $ since $C(H,G) \leq d$. For $\hat{\rho} \in U$, the composite cost $\tau(H, \hat{\rho},G)$ is finite so that $\tau(H',\rho,G)$ is also finite due to \eqref{bothfinite}. Apply Proposition~\ref{holdertau} to get a constant $\eta= \eta(D,a,\mathcal{P})$ given by \eqref{eta}, and  such that 
$$
\abs{ \tau(H',\rho,G) - \tau(H,\hat{\rho},G) }  \leq  \eta \left(\vnorm{\rho - \hat{\rho}}^{1/2} + \vnorm{ H'-H }\right).
$$

\hl{For} $\hat{\rho} \in U$, the bound in \eqref{rohatro} yields
\begin{equation} \label{tauHh}
\abs{\tau(H',\rho,G) - \tau(H,\hat{\rho},G) }  \leq  \eta_0 \vnorm{ H'- H }^{1/2}
\end{equation} 
with $\eta_0 = 2 (F_g / a)^{1/2} \eta $.
Since $\mathcal{L}$ is a bijection and $U\subset K(H)$, any $r \in U$ is of the form $r = \hat{s}$ for some $s = s(r) \in K(H')$. Then, for any $r \in U$, applying \eqref{tauHh} to $\rho = s$ and $r = \hat{s}$ implies 
$$
C(H',G) \leq \tau(H', s ,G)  \leq  \tau(H, r ,G) + \eta_0  \vnorm{ H'- H }^{1/2},
$$
and we can
\hl{take the} infimum of the right-hand side over all $r \in U$ to get 
\begin{equation}\label{ChCH}
C(H',G) \leq  C(H,G) + \eta_0 \vnorm{ H'- H }^{1/2} \leq d + \eta_0.
\end{equation}

{\hl{Define} $V = \vset{\rho \in K(H) \, | \, \tau(H',\rho,G) \leq d + \eta_0}$. This implies $C(H',G) = \min_{\rho \in V} \; \tau(H',\rho,G)$.  For $\rho \in V$, the cost $\tau(H', \rho,G)$ is finite so that $\tau(H,\hat{\rho},G)$ is also finite. Apply again Proposition~\ref{holdertau} with $A= d + \eta_0$ to get $\eta_1 = \eta(d +\eta_0, a, \calP)$ such that, with \eqref{rohatro},} 
\begin{equation} \label{tauHh2}
\abs{\tau(H',\rho,G) - \tau(H,\hat{\rho},G)}  \leq  \eta_1  \left(\vnorm{ \rho - \hat{\rho} }^{1/2} + \vnorm{ H' -H }\right)  \leq  \eta_2 \vnorm{ H' - H }^{1/2}
\end{equation} 
with $\eta_2 = 2 (F_g / a)^{1/2} \eta_1 $.  This implies, for all $\rho \in V,$
$$
\tau(H',\rho,G) \geq \tau(H,\hat{\rho},G) - \eta_2 \vnorm{ H' - H }^{1/2} \geq C(H,G) - \eta_2 \vnorm{ H' - H }^{1/2},
$$
and 
\hl{take the infimum} of the left-hand side over $\rho \in V$ to get 
$
C(H',G) \geq C(H,G) - \eta_2 \vnorm{ H' - H }^{1/2}.
$
Combine this with \eqref{ChCH}, and set $\gamma = \max\vset{\eta_2, \eta_0}$ to obtain
\begin{equation} \label{holderH'HG}
\abs{C(H',G) - C(H,G)} \leq \gamma \vnorm{ H' - H }^{1/2}.
\end{equation}
\hl{This force}s $C(H',G) \leq d + \gamma.$ Apply then \eqref{lipCHGG'} to obtain $c_3 = c_2(d + \gamma, a, \calP)$ such that 
$$
\abs{C(H',G) - C(H',G')} \leq c_3  \vnorm{ G' - G }.
$$
\hl{Combining} the above expression with \eqref{holderH'HG} yields Equation~\eqref{continuityCHG} with the constant $c = c_3 + \gamma$ and concludes the proof.

%%%%%%%%%%%%%%%%%%%%%%%%%%%%%
\section{Proofs of Theorems in Section~\ref{LDPathTheory}}
\label{sec:ap4}

%%%%%%%%%%%%%%%%%%%%%%%%%%%%
\subsection*{Proof of Theorem~\ref{hinfty}}

Define  $c > 1$ and $0 < \eta <1 $ by $c = 1 + \max_j \left( \sum_k q_{j,k} \right)$ and $1- \eta = \max_{1\leq k\leq g-1}  F_k/ F_{k+1}.$   Fix temporarily a constant $0 < A <1/2$ to be selected later on. Impose $m c \leq A $ so that all $m_j = m  \sum_k Q_{k,j} \leq A.$  For $H \in \mathcal{H},$ set $U_j(H) = F_j H(j) (1 - m_j) + m \sum_k Q_{k,j}F_k H(k)$. Setting $G(j)= \zeta_j(H) = U_j(H) /<F, H>,$ we have that $G(j) = 0$ if and only if $H(j) = 0$ and $H(k) = 0$ for all $k$ such that $Q_{k,j}>0.$  This implies $\supp(H) \subset \supp(\zeta(H)).$ Consequently, since $h_{n+1} = \zeta(h_n)$, the set $S_n =\supp(h_n)$ must increase with $n$. Therefore, there exists a fixed (sub)set of genotypes $\Ster$ and an $n_0$ such that $S_n = \Ster$ for $n \geq n_0$. An easy recurrence based on the conditions above for $G(j) = \zeta_j(H) = 0$ shows that $\Ster$ is the set of all $j$ reachable by some finite sequence $k_1, k_2,\cdots, k_T = j$ with $k_1 \in \supp(h_1)$ and all $q_{k_t, k_{t+1}} >0$.
 
Now, call a genotype $q$ \emph{dominant} in  $H$ if $H(q) = \max_k \, H(k) $. 
Then, we have $\langle F, H \rangle  \leq  g F_g H(q)$, $1/g \leq H(q)$, and
\begin{equation}\label{boundLj}
(1 - A) F_j H(j)  \leq  U_j (H)  \leq  F_j H(j) + A \langle F, H\rangle  \leq  F_j H(j)+ g A F_g H(q)
\end{equation}
for all $j$.  
Since $G(j) / G(q) = U_j(H) / U_q(H)$, this implies
\begin{equation}\label{GjGi}
\left(\frac{(1 - A) F_j} {( F_q + g A F_g )}\right) \left(\frac{H(j)}{H(q)}\right)  \leq  \frac{G(j)}{G(q)} \leq \frac{F_j H(j) + g A F_g H(q)}{F_q H(q) (1 - A)}.
\end{equation}

\hl{For}  $j < q$, one has $H(j) \leq H(q)$ and $F_j \leq (1 - \eta) F_q$ by definition of $c$ and $\eta.$ Then, \eqref{GjGi} yields 
\begin{equation}\label{upGG}
\frac{G(j)}{G(q)} \leq  \left(\frac{1-\eta}{1-A}\right)\frac{H(j)}{H(q)} + \frac{g A  F_g}{(1-A) F_q} \leq  \frac{1-\eta + A F_g/F_1}{1-A}.  
\end{equation}

\hl{Impose} $A < A_1 =\frac{\eta}{1+F_g/F_1}$ to get  
$
G(j)/G(q) \leq 1-\eta/2 < 1 
$
for all $j < q.$ Hence, the dominant genotype $q(G)$ in $G=\zeta(H)$ verifies $q(G)  \geq q(H)$.  Consequently, the dominant genotype $q_n = q(h_n)$ of $h_n$  verifies $q_{n+1} \geq q_n$ for all $n$ so that there exists a finite $n_1 \geq n_0$ and a genotype $q^*$ such that $q_n = q^*$ for all $n > n_1$. From \eqref{boundLj}, we get $h_n(q_n) \geq 1/g$ and hence $q^* \in \Ster$.  Let $s = \max(\Ster)$ be the fittest genotype in $\Ster$ so that $q^* \leq s$.  Assume there is a $j \in \Ster$ verifying $ q^* < j \leq s$. We will proceed by contradiction to show that such a $j$ cannot exist. For $n > n_1$, one has $h_n(j)  \leq  h_n(q^*)$ and, by \eqref{GjGi},
\begin{equation} \label{recur}
\theta \left(\frac{h_n(j)}{h_n(q^*)}\right) \leq \frac{h_{n+1}(j)}{h_{n+1}(q^*)}  \leq  1
\end{equation}
with $\theta  =  (1 - A) F_j / ( F_{q^*} + g A F_g ) $.  Since $j > q^*,$ we get $F_j  \geq  F_{q^*}/ (1-\eta)$ by definition of $\eta$, so that
$$
\theta  \geq  \frac{1 - A}{(1-\eta) ( 1 + g A  F_g/F_{q^*})}  \geq  \frac{1 - A}{(1-\eta) ( 1 + A g F_g/F_1 )}.
$$

\hl{Impose} $A < A_2 =\frac{\eta^2  F_1}{2 g F_g}$ to  force $\theta > 1 +\eta$. By recurrence, \eqref{recur} implies for $s \geq j > q*$ and $n' > n > n_1$ that $(1+\eta)^{n' - n}\frac{h_n(j)}{h_n(q^*)} \leq  1.$  Fix $n >n_1$. Since $j \in \Ster$, one has $h_n(j) >0.$ Letting $n' \to \infty$ yields a contradiction. Hence, there is no $j \in \Ster$ with $ q^* < j \leq s$, so that $q^* = s$.

For $i \in \Ster$ with $i < s$ and $n > n_1$, the terms $y_n(i) = h_{n}(i) / h_n(s)$ verify $y_{n+1}(i) \leq  \mu y_n(i) + \nu$  by \eqref{upGG} and $q_n= s$ with $\mu = 1- \eta / 2$ and $\nu = 2 g A  F_g / F_1$.  Iterating this inequality gives $y_n(i) \; \leq \; \mu^{n-n_1} y_{n_1}(i) + \frac{\nu}{1-\mu}$.
Select $n_2 = n_2(A) > n_1$ to force $(1-\eta / 2)^{n-n_1} y_{n_1}(i) < A$ for all $n > n_2$ and $i \in (\Ster \setminus s).$ Then, $y_n(i)  \leq c_1 A$, with $c_1 = 1 + \frac{4 g F_g }{\eta  F_1}$.  This yields $h_n(i)  \leq  c_1 \; A$ for all $n > n_ 2(A)$ and $i \in (\Ster \setminus s)$.  Hence, $h_n(s)  \geq  1 -  c_1 (g-1) A$ for $n > n_2(A)$, since $\supp(h_n) = \Ster$.  We now fix  $A = A_3 = \min\vset{ A_2, \frac{1}{2 g c_1} }$ and set $n_3 = n_2(A_3)$ to get the fixed lower bound $h_n(s) \geq 1/2$ for all $n> n_3.$  For $n > n_3$, let $z_n$ be the vector of all the $z_n(i) = y_n(i)$ with $i \in (S - s)$. Since $\supp(h_n) = S$, the non-zero $U_i(h_n)$ \emph{only depend on} $z_n$ and can  be denoted $U_i(z_n)$.  We then have $z_{n+1} = f(z_n)$, with $f_i(z)  = U_i(z) / U_s(z)$.  Set $p = card(S) -1.$ The rational fraction $f(z)$ is well-defined on the set $\Gamma$ of  all $z \in[0, 1]^p$ such that $U_s(z) \geq 1/2$ for $0 \leq m \leq A_3/c$.  For the limit case $m=0$ of no mutations, the function $f(z)$ has the form $\hat{f}_i (z) = \frac{F_i}{F_s}z_i$ for all $i \in (S - s)$. The Jacobian matrix $D_z \hat{f}$ obviously verifies
\begin{equation} \label{Dzf}
\vnorm{ D_z \hat{f}} \leq \max_{i \in (S - s)} \frac{F_i}{F_s}  \leq 1 - \eta
\end{equation}
for all $z \in [0,1]^p.$  For $z \in \Gamma$,  the numerator and denominator of each rational fraction $f_i(z)$ are separately affine in $m \leq A_3 / c$ and $z$, with the denominator bounded below by $1/2$ and uniformly bounded coefficients. Elementary algebraic computations then prove that as $m \to 0$, the Jacobian matrix $D_z f (z)$ tends to $D_z \hat{f}(z)$ \emph{uniformly} over all $z \in \Gamma$.  Due to \eqref{Dzf}, this provides a constant $0 < A_4 < A_3/c$ such that $\vnorm{ D_z f } \leq (1 - \eta/2 )$ for all $m < A_4$ and all $z \in \Gamma$. This yields $\vnorm{ z_{n+1} - z_n } \leq (1-\eta/2) \vnorm{ z_n - z_{n-1} }$ for $m < A_4$ and $n > 1+n_3$. This contraction property classically shows that $z^*= \lim_{n \to \infty} z_n$ exists with $z^* \in \Gamma.$ Furthermore, $z^*$ is the \emph{unique} solution of $z^* - f(z^*) = 0$. Note that $z^*$ is an implicit  function of $m$ for $0 \leq  m < A_4$, with $z^*(0) = 0$.  Due to \eqref{Dzf}, the Jacobian $Id - D_z f$ of $z -f(z)$ is invertible for $m =0$. Since $f(z)$ is a $C^{\infty}$ function of $(m, z) \in [0, A_4) \times\Gamma$, the implicit function theorem applies to $z - f(z) =0$ and proves the existence of a constant $0 < A_5 < A_4$ such that $z^*$ is a $C^{\infty}$ function of $m \in [0, A_5).$  Since $( 1- h_n(s) ) / h_n(s) = \sum_{i \in (S-s)} z_n(i)$ and $h_n(i) = h_n(s) z_n(i)$ for $i \in (S-s)$, we see that $H= \lim_{n \to \infty} h_n$ exists with $\supp(H) \subset S$ and verifies $H(s) = 1 / [ 1 + \sum_{i \in (S-s)}  z^*(i)]$ and $H(i) = H(s) z^*(i).$  Hence, $H$ is a $C^{\infty}$ function of $m \in [0, A_5)$ and $\supp(H) = S$. For $m=0$, the solution of $z - \hat{f}(z) = 0$ is clearly $z^* = 0$, and the associated histogram $H$ verifies $H(s) = 1$ and $\supp(H) = \vset{s }$.  The first-order Taylor expansion of $H$ as a function of $m$ is of the form $H(i) \simeq \;  m v_i$ for $i \in (S - s)$ and $H(s) \simeq \; 1- u m$. Substitute this into $H = \zeta(H)$ to get  $v_i = Q_{s, i} / (F_s - F_i)$ and $u = \sum_i v_i$. This concludes the proof.

%\bibliography{refs}   

\begin{thebibliography}{10}

\bibitem{assaf2017wkb}
{\sc M.~Assaf and B.~Meerson}, {\em {WKB} theory of large deviations in
  stochastic populations}, Journal of Physics A: Mathematical and Theoretical,
  50 (2017), p.~263001.

\bibitem{azenLDT2012}
{\sc R.~Azencott, M.~I. Freidlin, and S.~Varhadan}, {\em Large Deviations at
  Saint-Flour}, Springer, 2012.

\bibitem{Barrick2013}
{\sc J.~E. Barrick and R.~E. Lenski}, {\em Genome dynamics during experimental
  evolution}, Nature Reviews. Genetics, 14 (2013), pp.~827--839.

\bibitem{Barrick2010}
{\sc J.~E. Barrick, C.~C. Strelioff, R.~E. Lenski, and M.~R. Kauth}, {\em
  \textit{Escherichia coli} {rpoB} mutants have increased evolvability in
  proportion to their fitness defects}, Molecular Biology and Evolution, 27
  (2010), pp.~1338--1347.

\bibitem{darwin2011}
{\sc F.~A. C.~C. Chalub and J.~F. Rodrigues}, eds., {\em The Mathematics of
  {Darwin’s} Legacy (Mathematics and Biosciences in Interaction)},
  Birkhäuser, 2011.

\bibitem{champ6}
{\sc N.~Champagnat}, {\em A microscopic interpretation for adaptive dynamics
  trait substitution sequence models}, Stoch. Proc. Appl., 116 (2006),
  pp.~1127--1160.

\bibitem{champ01}
{\sc N.~Champagnat, R.~Ferriere, and G.~B. Arous}, {\em The canonical equation
  of adaptive dynamics: A mathematical view}, Selection, 2 (2001), pp.~71--81.

\bibitem{champ8}
{\sc N.~Champagnat, R.~Ferri\`ere, and S.~M\'el\'eard}, {\em Unifying
  evolutionary dynamics: From individual stochastic processes to macroscopic
  models}, Theoretical Population Biology, 69 (2006), pp.~297 -- 321.
\newblock ESS Theory Now.

\bibitem{champ02}
{\sc N.~Champagnat and A.~Lambert}, {\em Evolution of discrete populations and
  the canonical diffusion of adaptive dynamics}, Ann. Appl. Prob., 17 (2007),
  pp.~102--155.

\bibitem{cooper1}
{\sc T.~F. Cooper, D.~E. Rozen, and R.~E. Lenski}, {\em Parallel changes in
  gene expression after 20,000 generations of evolution in {E. coli}}, Proc.
  Nat. Acad. Sci. USA, 100 (2003), pp.~1072--1077.

\bibitem{Cooper2001}
{\sc V.~S. Cooper, D.~Schneider, M.~Blot, and R.~E. Lenski}, {\em Mechanisms
  causing rapid and parallel losses of ribose catabolism in evolving
  populations of \textit{Escherichia coli}}, Journal of Bacteriology, 183
  (2001), pp.~2834--2841.

\bibitem{Gordo2013}
{\sc J.~A.~M. de~Sousa, P.~R.~A. Campos, and I.~Gordo}, {\em An {ABC} method
  for estimating the rate and distribution of effects of beneficial mutations},
  Genome Biology and Evolution, 5 (2013), pp.~794--806.

\bibitem{Deatherage2017}
{\sc D.~E. Deatherage, K.~L. Jamie, A.~F. Bennett, R.~E. Lenski, and J.~E.
  Barrick}, {\em Specificity of genome evolution in experimental populations of
  \textit{Escherichia coli} evolved at different temperatures}, Proc. Nat.
  Acad. Sci. USA, 114 (2017), pp.~E1904--E1912.

\bibitem{dembo}
{\sc A.~Dembo and O.~Zeitouni}, {\em Large deviations Techniques and
  applications}, Springer-verlag, 1998.

\bibitem{Desai2007}
{\sc M.~Desai and D.~Fisher}, {\em Beneficial mutation–selection balance and
  the effect of linkage on positive selection}, Genetics, 176 (2007),
  pp.~1759--1798.

\bibitem{durrett2008}
{\sc R.~Durrett}, {\em Probability Models for {DNA} Sequence Evolution
  (Probability and Its Applications)}, Springer, 2008.

\bibitem{Fox2015}
{\sc J.~W. Fox and R.~E. Lenski}, {\em From here to eternity - the theory and
  practice of a really long experiment}, {PLoS} Biology, 13 (2015),
  p.~e1002185.

\bibitem{geiger_thesis}
{\sc B.~J. Geiger}, {\em Large Deviations for Dynamical Systems with Small
  Noise}, PhD thesis, University of Houston, https://hdl.handle.net/10657/4801,
  2017.

\bibitem{gomez2023markovian}
{\sc A.~G{\'o}mez-Corral, M.~J. Lopez-Herrero, and D.~Taipe}, {\em A markovian
  epidemic model in a resource-limited environment}, Applied Mathematics and
  Computation, 458 (2023), p.~128252.

\bibitem{Gordo2012}
{\sc I.~Gordo, L.~Perfeito, and A.~Sousa}, {\em Fitness effects of mutations in
  bacteria}, Journal of Molecular Microbiology and Biotechnology, 21 (2012),
  pp.~20--35.

\bibitem{heger}
{\sc M.~Hegreness, N.~Shoresh, D.~Hartl, and R.~Kishony}, {\em An equivalence
  principle for the incorporation of favorable mutations in asexual
  populations}, Science, 311 (2006), pp.~1615--1617.

\bibitem{Illinworth2012}
{\sc C.~Illinworth and V.~Mustonen}, {\em A method to infer positive selection
  from marker dynamics in asexual population}, Bioinformatics, 28 (2012),
  pp.~831--837.

\bibitem{Levy2015}
{\sc S.~F. Levy, J.~R. Blundell, S.~Venkataram, D.~A. Petrov, D.~S. Fisher, and
  G.~Sherlock}, {\em Quantitative evolutionary dynamics using high-resolution
  lineage tracking}, Nature, 519(7542) (2015), pp.~181--186.

\bibitem{meleard2016}
{\sc S.~Meleard}, {\em Modèles aléatoires en Ecologie et Evolution
  (Mathématiques et Applications, 77) (French)}, Springer, 2016.

\bibitem{meleard2015}
{\sc S.~Meleard and V.~Bansaye}, {\em Stochastic Models for Structured
  Populations: Scaling Limits and Long Time Behavior (Mathematical Biosciences
  Institute Lecture Series, 1.4)}, Springer, 2015.

\bibitem{papageorgiou2025enhanced}
{\sc V.~E. Papageorgiou, G.~Vasiliadis, and G.~Tsaklidis}, {\em An enhanced
  epidemic susceptible-infected-hospitalized-recovered-deceased {(SIHRD)}
  stochastic model with emphasis on the impact of hospitalizations on epidemic
  evolution}, Methodology and Computing in Applied Probability, 27 (2025),
  p.~55.

\bibitem{Peng2017}
{\sc F.~Peng, S.~Widmann, A.~W\"unsche, K.~Duan, K.~A. Donovan, R.~C.~J.
  Dobson, R.~E. Lenski, and T.~F. Cooper}, {\em Effects of beneficial mutations
  in \textit{pykF} gene vary over time and across replicate populations in a
  long-term experiment with bacteria}, Molecular Biology and Evolution, 35
  (2018), pp.~202--210.

\bibitem{Plank1979}
{\sc L.~D. Plank and J.~D. Harvey}, {\em Generation time statistics of
  \textit{Escherichia coli} {B} measured by synchronous culture techniques},
  Journal of General Microbiology, 115 (1979), pp.~69--77.

\bibitem{rice}
{\sc S.~H. Rice}, {\em Evolutionary Theory}, Sinauer Assocites, 2004.

\bibitem{Ross1996}
{\sc S.~M. Ross}, {\em Stochastic Processes}, John Wiley $\&$ Sons, New York,
  1996.

\bibitem{schoneman_thesis}
{\sc J.~Schoneman}, {\em Stochastic Models for Genetic Evolution}, PhD thesis,
  University of Houston, http://hdl.handle.net/10657/3540, 2016.

\bibitem{Simon2013}
{\sc D.~Simon}, {\em Evolutionary Optimization algorithms}, John Wiley $\&$
  Sons, New York, 2013.

\bibitem{su_thesis}
{\sc Y.~Su}, {\em Rare Events Simulation In Bacterial Genetic Evolution
  Models}, PhD thesis, University of Houston,
  https://hdl.handle.net/10657/10734, 2022.

\bibitem{azencott2023rare}
{\sc Y.~Su, B.~Geiger, I.~Timofeyev, A.~Mang, and R.~Azencott}, {\em Rare
  events analysis and computation for stochastic evolution of bacterial
  populations}, Stochastic Analysis and Applications, 43(1) (2024), pp.~1--29.

\bibitem{len94a}
{\sc F.~Vasi, M.~Travisano, and R.~Lenski}, {\em Long-term experimental
  evolution in \textit{Escherichia coli}. ii. changes in life-history traits
  during adaptation to a seasonal environment}, The American Naturalist, 144
  (1994), pp.~432--456.

\bibitem{komarova2005}
{\sc D.~Wodarz and N.~Komarova}, {\em Computational Biology of Cancer: Lecture
  Notes and Mathematical Modeling}, World Scientific, 2005.

\bibitem{Woods2011}
{\sc R.~J. Woods, J.~E. Barrick, T.~F. Cooper, U.~Shrestha, M.~R. Kauth, and
  R.~E. Lenski}, {\em Second-order selection for evolvability in a large
  \textit{Escherichia coli} population}, Science, 331 (2011), pp.~1433--1436.

\bibitem{azencoop}
{\sc W.~Zhang, V.~Sehgal, D.~Dinh, R.~R. Azevedo, T.~Cooper, and R.~Azencott},
  {\em Estimation of the rate and effect of new beneficial mutations in asexual
  populations}, Theoretical Population Biology, 81 (2012), pp.~168--178.

\end{thebibliography}

\end{document}